\documentclass{amsart}

\usepackage[latin1]{inputenc}
\usepackage[english]{babel}
\usepackage{amsmath}
\usepackage{amssymb}
\usepackage{centernot}
\usepackage{caption}
\usepackage{color}
\usepackage{subcaption}
\usepackage{etoolbox}
\usepackage{upgreek}
\usepackage{mathabx}
\usepackage{titlesec}
\usepackage{array}
\usepackage{enumitem}
\usepackage{makecell}
\patchcmd{\section}{\scshape}{\bfseries}{}{}
\makeatletter
\renewcommand{\@secnumfont}{\bfseries}
\makeatother
\newcommand*{\justifyheading}{\raggedright}
\titleformat {\section}
  {\normalfont \Large \bfseries \justifyheading}{\thesection.}{1em}{}
\titleformat {\subsection}
  {\normalfont \large \bfseries \justifyheading}{\thesubsection.}{1em}{}
\patchcmd{\abstract}{\scshape\abstractname}{\textbf{\abstractname}}{}{}

\numberwithin{equation}{section}
\newtheorem{thm}{Theorem}[section]
\newtheorem{prop}[thm]{Proposition}
\newtheorem{lem}[thm]{Lemma}
\newtheorem{cor}[thm]{Corollary}

\newtheorem{theo}[thm]{Theorem}
\newtheorem{defn}[thm]{Definition}

\theoremstyle{definition}

\newtheorem{conv}[thm]{Convention}
\newtheorem{rem}[thm]{Remark}

\newtheorem{quest}[thm]{Question}
\newtheorem{exam}[thm]{Example}
\newenvironment{proof2}[1][Proof.]{\begin{trivlist}
\item[\hskip \labelsep {\bfseries \itshape #1}]}{\end{trivlist}}
\newenvironment{proof3}[1][Proof.]{\begin{trivlist}
\item[\hskip \labelsep {\bfseries \itshape #1}]}{\end{trivlist}}

\newcommand{\expe}[1]{\text{\bf{E$\left[\right.$}}#1\text{\bf{$\left.\right]$}}}
\newcommand{\proba}[1]{\text{\bf{Pr$\boldsymbol{\left[\right.}$}}#1\text{$\boldsymbol{\left.\right]}$}}

\newcommand{\Zon}[1]{\mathcal{Z}_{#1}^{\scriptscriptstyle\emph{central}}}

\newcommand{\R}{\mathbb{R}}
\newcommand{\QQ}{\mathbb{Q}}
\newcommand{\RNN}{\mathbb{R}_{\geq 0}}
\newcommand{\RP}{\mathbb{R}_{+}}
\newcommand{\PN}{\mathbb{P}}
\newcommand{\NN}{\mathbb{N}}
\newcommand{\kfield}{\mathbf{k}}
\newcommand{\indeg}[2]{\emph{indeg}_{\scriptscriptstyle (#1,#2)}}
\newcommand{\outdeg}[2]{\emph{outdeg}_{\scriptscriptstyle (#1,#2)}}
\newcommand{\nooutdeg}[2]{\emph{nod}_{\scriptscriptstyle (#1,#2)}}
\newcommand{\indegvec}[2]{\mathbf{\mathtt{indeg}}_{\scriptscriptstyle (#1,#2)}}
\newcommand{\outdegvec}[2]{\mathbf{\mathtt{outdeg}}_{\scriptscriptstyle (#1,#2)}}
\newcommand{\nooutdegvec}[2]{\mathbf{\mathtt{nod}}_{\scriptscriptstyle (#1,#2)}}
\newcommand{\Arr}[1]{\mathcal{A}_{#1}}

\newcommand{\Conv}[1]{\emph{conv}\left(#1\right)}

\newcommand{\rgr}[2]{G_{#1,#2}}

\newcommand{\upset}[1]{2^{#1}{\displaystyle\backslash}{\scriptstyle\{\emptyset\}}}

\newcommand{\midset}[1]{2^{#1}{\displaystyle\backslash}{\scriptstyle\{#1,\emptyset\}}}
\newcommand{\degvec}[1]{\mathbf{d}_{\scriptscriptstyle #1}}
\newcommand{\degnovec}[1]{d_{\scriptscriptstyle #1}}
\newcommand{\deginvec}[2]{\mathbf{\mathtt{inof}}_{\scriptscriptstyle (#1,#2)}}
\newcommand{\degoutvec}[2]{\mathbf{\mathtt{outof}}_{\scriptscriptstyle (#1,#2)}}
\newcommand{\degin}[2]{\emph{inof}_{\scriptscriptstyle (#1,#2)}}
\newcommand{\degout}[2]{\emph{outof}_{\scriptscriptstyle (#1,#2)}}
\newcommand{\CZon}[1]{\mathcal{Z}_{#1}}
\newcommand{\shull}[2]{\left( #1,#2 \right)}
\newcommand{\hull}[2]{\left[#1,#2\right]}
\newcommand{\aff}[1]{\emph{aff}\left(#1\right)}
\newcommand{\relint}[1]{\emph{relint}\left\langle#1 \right\rangle}
\newcommand{\inter}[1]{\emph{int}\left\langle#1 \right\rangle}
\newcommand{\bound}[1]{\partial \left\langle#1 \right\rangle}
\newcommand{\chat}[1]{\widehat{#1}}
\newcommand{\height}[2]{{\scriptstyle\emph{\bf height}}_{\scriptscriptstyle #1}^{\scriptscriptstyle #2}}
\newcommand{\depth}[2]{{\scriptstyle\emph{\bf depth}}_{\scriptscriptstyle #1}^{\scriptscriptstyle #2}}
\newcommand{\noncross}[1]{{\mathbf{NC}}({\scriptstyle #1})}
\newcommand{\precref}{\ensuremath{\prec_{\scriptscriptstyle\emph{ref}}\!\!\!\!\!\!\!\!\cdot\,} \ \,}

\usepackage{graphicx}
\usepackage{mathrsfs}
\usepackage[foot]{amsaddr}

\usepackage[round]{natbib}

\title[Acyclic orientations and spanning trees.]{Acyclic orientations and spanning trees.}
\author{Benjamin Iriarte}
\address{Department of Mathematics, Massachusetts Institute of Technology, Cambridge MA, 02139, USA}
\thanks{The author was supported by NSF grant DMS-1068625 during the entirety of this work.}
\email{biriarte@math.mit.edu}


\begin{document}
\begin{abstract}
We introduce polytopal cell complexes associated with partial acyclic orientations of a simple graph, 
which generalize acyclic orientations. Using the theory of cellular resolutions, two of these polytopal cell complexes are 
observed to minimally resolve
certain special combinatorial polynomial ideals related to acyclic orientations. These ideals are
explicitly found to be Alexander dual, which relative to comparable results in the 
literature, generalizes in a cleaner and more illuminating 
way the well-known duality between permutohedron and tree ideals. 
The combinatorics underlying these results naturally leads to a canonical way to represent rooted spanning forests of a labelled 
simple graph as non-crossing trees,
and these representations are observed to carry a plethora of information about generalized tree ideals and acyclic orientations of a graph, 
and about non-crossing partitions of a totally ordered set. A small sample of the enumerative and structural 
consequences of collecting and organizing this information are studied in detail. Applications of this combinatorial
miscellanea are then introduced and explored, namely: Stochastic processes on state space
equal to the set of all 
acyclic orientations of a simple graph, including
irreducible Markov chains, which exhibit stationary distributions ranging from 
linear extensions-based to uniform; a surprising formula for the expected number
of acyclic orientations of a random graph; and a purely algebraic presentation of the main problem in bootstrap percolation, likely making 
it tractable to explore the set of all percolating sets of a graph with a computer.      
\end{abstract}

\maketitle
\section{Introduction.}

This article is a sequel to~\cite{iriartegraphs}, focusing instead on the structural and enumerative properties of acyclic orientations. We 
introduce a number of novel perspectives, results and resources 
for the study and discovery of fundamental 
properties of acyclic orientations, 
and their generalization, partial
acyclic orientations, of a simple graph; 
these include polytopal cell complexes and polynomial ideals~{\bfseries [}\cite{millersturmfels}{\bfseries ]}, graphical 
zonotopes~{\bfseries [}\cite{postnikovpermutohedra},~\cite{beckbook}{\bfseries ]}, 
and 
Markov chains~{\bfseries [}\cite{lovaszrandom},~\cite{aldoustrees}{\bfseries ]}, among others. 
We adopt an original approach to the well-known connection between labelled trees, parking functions, non-crossing partitions, and 
graph orientations. This is the viewpoint of non-crossing trees, not
properly treated or even reported in the literature, and which we exploit to obtain
new results about these objects. Non-crossing trees are, in part,
motivated by the techniques of~\cite{finkbijections}, but owe their existence to the 
(subtleties of the) Alexander duality~{\bfseries [}\cite{milleralexander}{\bfseries ]} between two special polynomials ideals defined during 
Section~\ref{sec:twoideals} of the present writing. 
The perspectives presented here complement those of previous
key studies, including but not exhausting, those found in~\cite{chebikintrees},~\cite{postnikovtrees},~\cite{dochtermannsanyal},
~\cite{persiansminimal},~\cite{persiansdivisors},~\cite{stanleynoncross} and~\cite{stanleyparking}, and the references therein. 
The present work is, in fact, an evident seed for future research more than a conclusive exposition of the topic, and the number of
(sometimes quite provocative) open problems and directions for future research should gradually become clear. 
\par
In principle, an inconvenient aspect of acyclic orientations of a simple graph is their apparent but, nevertheless, artificial relation to bijective labellings
of the vertex set with a totally ordered set. This viewpoint was exploited during the author's previous article on this subject.  
Conceivably, adopting a perspective different to that of bijective labellings seems equally fated to
illuminate the study of acyclic orientations of a simple graph, and this is what we pursue in this writing. 
One example of how we apply ideas developed during this sequel is the construction of
a random walk on a certain simple connected regular graph with vertex set equal to the set of all acyclic orientations of any fixed simple graph, 
which therefore exhibits a unique uniform stationary distribution. The importance and applicability of such constructions
is evidently exemplified in~\cite{brodertrees},~\cite{aldousrandom} and~\cite{kelnertrees}, 
and formalized in~\cite{lovaszmixing} 
and~\cite{lovaszefficient}. Many other fine works
have made use of similar ideas to solve different combinatorial algorithmic problems.  
\par
Understanding acyclic orientations of a simple graph from their grounds usually entails making precise connections
of their theory with the theory of spanning trees, much better understood; this is also the case in the present 
article. The
particular connection between these sets of combinatorial objects that we choose to follow, developed here for the first time, is far 
from obvious and will be presented later in Section~\ref{sec:noncrosstrees}, 
where it sprouts naturally from the constructions of Sections~\ref{sec:polytopes}-\ref{sec:twoideals}. 
\par
Let us describe in fair detail the contents of the different sections of the article. 
\par
In Section~\ref{sec:polytopes} we introduce, again for the first time in the literature, 
an elegant inequality description of a well-known polytope related to the acyclic orientations of
a fixed simple (connected) graph on vertex set $[n]$, $n\in\PN$; it can be found in Subsection~\ref{subsec:zonotope}. 
The above description has the form predicted in~\cite{postnikovpermutohedra} for 
the generalized permutohedra. This polytope of 
\emph{partial acyclic orientations} has a Minkowski sum decomposition whose summands
appear also as summands in Postnikov's expression of the graph associahedron as a sum of simplices.  
A first step along this road from the polytope of partial acyclic orientations to the graph
associahedron of \emph{graph tubings}~{\bfseries [}\cite{devadoss}{\bfseries ]} leads us to consider, in the case of connected graphs, the 
Minkowski sum of the former polytope with an $(n-1)$-dimensional 
simplex. The construction of Cayley's trick applied to this case 
serves us to discover one more polytopal cell complex associated to the graph, a complex pivotal 
in the study of certain ``artinianizations'' of the ideals defined in Section~\ref{sec:twoideals} and (therefore) instrumental
in the search for minimal free resolutions of these ideals~{\bfseries [}\cite{bayercellular}{\bfseries ]}, and whose
combinatorial dual is precisely the totally non-negative part of the graphical arrangement~{\bfseries [}\cite{stanleyhyperplanes}{\bfseries ]}.  
\par
In Section~\ref{sec:twoideals}, this clean geometrical duality of polytopal cell complexes manifests itself as an algebraic duality between two
polynomial ideals associated to a fixed simple connected graph, defined therein; one of these ideals is motivated by the role of acyclic orientations 
in the graphical zonotope, and the other by the inequality description of the polytope of Subsection~\ref{subsec:zonotope}.
The proof
of this \emph{Alexander duality}, found early in the section, contains the stepping stones for Section~\ref{sec:noncrosstrees}. 
We regard
some of the results contained in this section as being ``close siblings'' to those found in~\cite{dochtermannsanyal},
~\cite{persiansminimal} and~\cite{persiansdivisors}, yet our {\it modus operandi} aims to fix a 
necessarily problematic (at least for our purposes)
aspect of these other works: The generalization of the duality between the (standard) permutohedron and tree ideals implicit in them is by no means
self-evident nor truly discussed, and it does not follow from a clean geometrical duality generalizing the picture of the permutahedron and
the barycentric subdivision of the simplex; as such, these other perspectives do not yield the algorithmic consequences that we need later on in
Sections~\ref{sec:noncrosstrees}-\ref{sec:applications}. 
\par
Section~\ref{sec:noncrosstrees} introduces non-crossing trees of a simple graph, certain pictorial representations of labelled rooted trees
reminiscent of~\cite{finkbijections}.
There is one non-crossing tree per each rooted spanning forest of the graph. In Subsection~\ref{subsec:standardmon}, we explain how 
each non-crossing tree naturally encodes a uniquely determined standard monomial of the generalized tree ideal, defined in Section~\ref{sec:twoideals}, and
(therefore) a uniquely determined orientation of the graph with no directed cycles. Among these orientations supported on non-crossing trees, we find
the acyclic orientations of the graph, which spring up, again naturally, from non-crossing trees satisfying a certain efficiency condition. 
In Subsection~\ref{subsec:noncrosspart}, we adopt ``the other'' point of view on non-crossing trees, and observe how we then obtain
chains of the non-crossing partitions lattice. These two points of view are combined to produce a coherent picture of the combinatorial objects 
involved in this work.
\par
Section~\ref{sec:applications} contains applications of the ideas developed in Sections~\ref{sec:polytopes}-\ref{sec:noncrosstrees} to
algorithmic/computational problems involving (mostly random) acyclic orientations.  
Subsection~\ref{subsec:randwalk} presents five different stochastic processes on state spaces equal to the set of all acyclic orientations
of a simple graph, and whose stationary distributions range from 
dependent on the number of linear extensions~{\bfseries [}as in~\cite{iriartegraphs}{\bfseries ]} 
to uniform. In order of appearance, these are the \emph{Card-Shuffling Markov chain}, the \emph{Edge-Label Reversal} and the
\emph{Sliding-$(n+1)$ stochastic processes}, the \emph{Cover-Reversal random walk}, and the \emph{Interval-Reversal random walk}.
The Card-Shuffling Markov chain had also been previously discovered as a hyperplane walk
in~\cite{athanasiadisdiaconis}, and the Cover-Reversal random walk is grounded on the work of~\cite{savageconnectivity} and of
Section~\ref{sec:polytopes} of the present writing. 
This subsection culminates with the presentation of the Interval-Reversal random walk, an irreducible
reversible Markov chain with uniform stationary distribution on the acyclic orientations of a simple graph, never presented
before in the literature, and motivated by a close inspection of Section~\ref{sec:polytopes} here. Subsection~\ref{subsec:randgraph}
presents a surprising expression for the expected number of acyclic orientations of an Erd\"os-R\'enyi random graph in terms
of \emph{parking functions}, a consequence of the study of non-crossing trees in Section~\ref{sec:noncrosstrees}. Subsection~\ref{subsec:kbootstrap}
introduces a commutative-algebraic approach to determining all percolating sets in \emph{$k$-bootstrap percolation} on
any simple graph~{\bfseries [}{\it e.g.}~\cite{bootstrap}{\bfseries ]}; this direction could yield good fruits if further explored in the future. 

\begin{proof3}[Acknowledgements:]
I would like to specially thank my advisor Richard P. Stanley and Jacob Fox, whose support and always useful advice
made it possible to write this work. 
\end{proof3}     

\section{Polytopal complexes for acyclic orientations.}\label{sec:polytopes}

\subsection{A Classical Polytope.}\label{subsec:zonotope}

\begin{defn}\label{defn:orienfunc}
Let $G=G(V,E)$ be a simple graph and let: 
$$\overleftrightarrow{E}:=\left\{(u,v)\in V^2:\{u,v\}\in E\right\}.$$ 
An \emph{orientation $O$ of $G$} is a function $O:E\rightarrow E\cup \overleftrightarrow{E}$ such that
for all $e=\{u,v\}\in E$, we have that $O(e)\in\left\{e,(u,v),(v,u)\right\}$. We will let
$O_{\scriptscriptstyle\emph{trivial}}$ be the identity map $E\rightarrow E$. 
\end{defn}

\begin{defn}\label{defn:part}
For a simple graph $G=G(V,E)$, a partition $\Upsigma$ of the set $V$ is said to be a \emph{connected partition of $G$} if
$G[\upsigma]$ is connected for all $\upsigma\in\Upsigma$, where $G[\upsigma]$ denotes the induced subgraph of $G$ on
vertex set $\sigma$.
\end{defn}

\begin{defn}\label{defn:partgraph}
Let $G=G(V,E)$ be a simple graph and $\Upsigma$ a connected partition of $G$. Then, the \emph{$\Upsigma$-partition graph \
$G^{\Upsigma}=G^{\Upsigma}(\Upsigma,E^{\Upsigma})$ of $G$} is 
the graph such that, for $\upsigma,\uprho\in\Upsigma$ with $\upsigma\neq \uprho$, $\{\upsigma,\uprho\}\in E^{\Upsigma}$ if and only if there exists
$u\in \upsigma$ and $v\in \uprho$ with $\{u,v\}\in E$.
\end{defn}

\begin{defn}\label{defn:pao1}
Let $G=G(V,E)$ be a simple graph. An orientation $O$ of $G$ is said to be a \emph{partial acyclic orientation (p.a.o.) of $G$} if $O$ can be obtained in the following way: 

There exists a connected partition $\Upsigma$ of $G$
and an acyclic orientation $O^{\Upsigma}$ of the $\Upsigma$-partition graph $G^{\Upsigma}$ of $G$ such that, for all $e=\{u,v\}\in E$:
\begin{enumerate}
\item If $e\subseteq\upsigma\in\Upsigma$, then $O(e) = e$.
\item\label{defn:pao1:it2} If $u\in \upsigma$ and $v\in\uprho$ for some $\upsigma,\uprho\in \Upsigma$ with $\upsigma\neq\uprho$, \
and if $O^{\Upsigma}(\{\upsigma,\uprho\})=(\upsigma,\uprho)$, then $O(e) = (u,v)$.
\end{enumerate}
\par
We will also consider two functions, $\dim_G$ and $J_G$, associated to the set of \emph{p.a.o.}'s of $G$. To define them, 
let $O$ be a \emph{p.a.o.} of $G$
with associated connected partition $\Upsigma$.
\par
The first function, $\dim_G$, maps from the set of all \emph{p.a.o.}'s of $G$ to $\NN$, and is given as:
\begin{equation*}
\dim_G(O)=|V|-|\Upsigma|.
\end{equation*}
\par
The second function, $J_G$, has also domain the \emph{p.a.o.}'s of $G$, but it maps to the set of finite distributive lattices:
\begin{align*}
&J_G(O)=J(O^{\Upsigma}), \\
&\text{where $J(O^{\Upsigma})$ is the poset of order ideals of $O^{\Upsigma}$.}
\end{align*}
\end{defn}

\begin{rem}\label{rem:pao}
For a \emph{p.a.o.} $O$ of $G$ as in Definition~\ref{defn:pao1}, we will often identify $O$ with its induced partially ordered set $(V,\leq_{O})$, where
for all $u,v\in V$ we have that $u<_{O} v$ if and only if $u\in \upsigma$ and $v\in\uprho$ for some $\upsigma,\uprho\in\Upsigma$ with $\upsigma\neq \uprho$, and there exist 
$\upsigma_0,\upsigma_1,\dots,\upsigma_k \in \Upsigma$ with $\upsigma_0=\upsigma$ and $\upsigma_k =\uprho$ 
such that $(\upsigma_{i-1},\upsigma_{i})\in O^{\Upsigma}\left[E^{\Upsigma}\right]$ for all
$i\in[k]$.  
\end{rem}

\begin{lem}\label{lem:easy}
In Definition~\ref{defn:pao1}, if $O$ is a \emph{p.a.o.} of $G$, then $\dim_G(O)$ is equal to 
$|V|-\mathbf{l}\left(J_G(O)\right)$, where $\mathbf{l}(\cdot)$ denotes the length function for graded posets. 
\end{lem}
\begin{proof2}
Let $\Upsigma$ be the connected partition of $G$ associated to $O$. The result follows 
since then $\mathbf{l}\left(J(O^{\Upsigma})\right) = |\Upsigma|$.
\qed\end{proof2}

\begin{lem}\label{lem:ordideal}
In Definition~\ref{defn:pao1}, consider a \emph{p.a.o.} $O$ of $G$, and for $I\in \underline{J_G}(O)$, let $I^u = \bigcup_{\upsigma\in I}\upsigma$. If we let
$P$ be the poset of all $I^u$ with $I\in \underline{J_G}(O)$, ordered by inclusion of sets, then $P\simeq J_G(O)$.
\end{lem}
\begin{proof2}
This is straightforward, since for $I_1,I_2\in \underline{J_G}(O)$, both $I_1\cap I_2\in \underline{J_G}(O)$ and $I_1\cup I_2\in \underline{J_G}(O)$.
\qed\end{proof2}
\begin{rem}\label{rem:notation}
Naturally, in Lemma~\ref{lem:ordideal} and in subsequent writing, for $O$ a \emph{p.a.o.} of $G$, $\underline{J_G}(O)$ denotes the ground set
of $J_G(O)$.  
\end{rem}
\begin{rem}\label{rem:ordideal}
In fact, following Lemma~\ref{lem:ordideal}, in Definition~\ref{defn:pao1} we will regard
$J_G(\cdot)$ as a collection of subsets of $V$ ordered by inclusion. 
\end{rem}
\begin{lem}\label{lem:injectiveord}
In Definition~\ref{defn:pao1} and Remark~\ref{rem:notation}, the map $\underline{J_G}$ is an injective map.
\end{lem}
\begin{proof2}
Let $O_1$ and $O_2$ be \emph{p.a.o.}'s of $G$ such that $\underline{J_G}(O_1) = \underline{J_G}(O_2)$. 
Hence, $J_G(O_1) = J_G(O_2)$. Considering a maximal chain $\emptyset = \upsigma_0 \subsetneq \cdots \subsetneq \upsigma_{k} = V$
of this poset, we observe that $\Upsigma = \left\{\upsigma_{i}\backslash\upsigma_{i-i}\right\}_{i\in[k]}$ is the connected partition
of $G$ associated to both $O_1$ and $O_2$. The poset of \emph{join-irreducibles} of $J_G(O_1) = J_G(O_2)$
determines a unique acyclic orientation $O^{\Upsigma}$ of the $\Upsigma$-partition graph $G^{\Upsigma}$, and so
both $O_1$ and $O_2$ are obtained from $O^{\Upsigma}$ as in Definition~\ref{defn:pao1}.\ref{defn:pao1:it2}. Clearly then
$O_1 = O_2$. 
\qed\end{proof2} 
\begin{defn}\label{defn:complexzon}
Consider a simple graph $G=G(V,E)$. We will define an abstract cell complex
$\mathscr{Z}_G=(\underline{\mathscr{Z}_G},\preceq_z,\dim_z)$, with underlying set
of faces $\underline{\mathscr{Z}_G}$ ordered by $\preceq_z$, and with dimension function $\dim_z$, in the following manner:
\begin{enumerate}
\item $\underline{\mathscr{Z}_G}$ is the set of \emph{p.a.o.}'s of $G$.
\item For $O_1,O_2$ \emph{p.a.o.}'s of $G$, $O_1\preceq_z O_2$ if and only if $\underline{J_G}(O_2)\subseteq \underline{J_G}(O_1)$.
\item For $O$ a \emph{p.a.o.} of $G$, $\dim_z(O) = \dim_G(O)$. 
\end{enumerate}
\end{defn}

\begin{figure}[ht]
\begin{tabular}{lll}
\begin{subfigure}{.24\textwidth}
  \raggedright
  \includegraphics[width=1.22\linewidth]{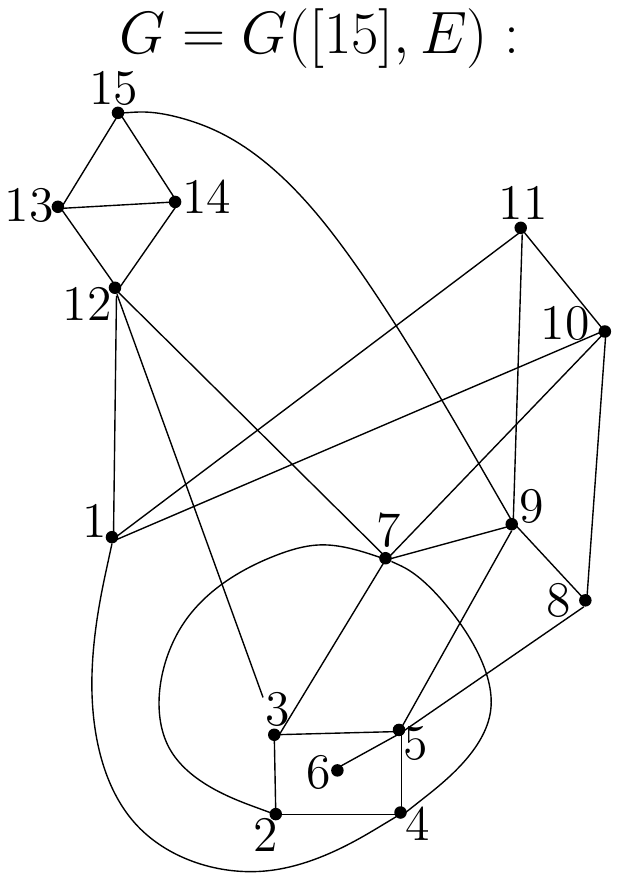}
  \caption{}
  \label{fig:defpao1}
\end{subfigure}
& 
\makecell{\begin{subfigure}{.38\textwidth}
  \raggedright
  \includegraphics[width=1.05\linewidth]{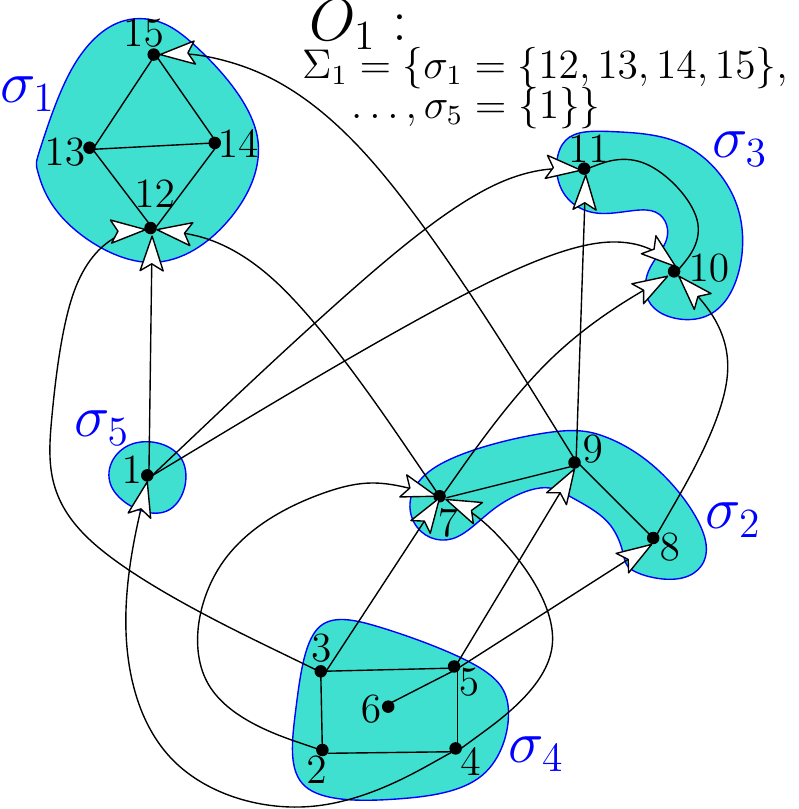}
  \subcaption{}\label{fig:defpao2}
\end{subfigure}
\\
\begin{subfigure}{.38\textwidth}
  \raggedleft
  \subcaption{}\label{fig:defpao3}
  \includegraphics[width=.5\linewidth]{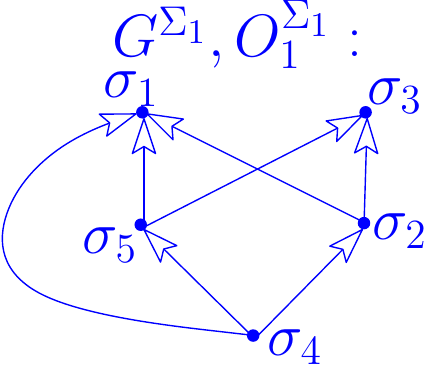}
\end{subfigure}}
&
\makecell{\begin{subfigure}{.38\textwidth}
  \raggedright
  \includegraphics[width=1.05\linewidth]{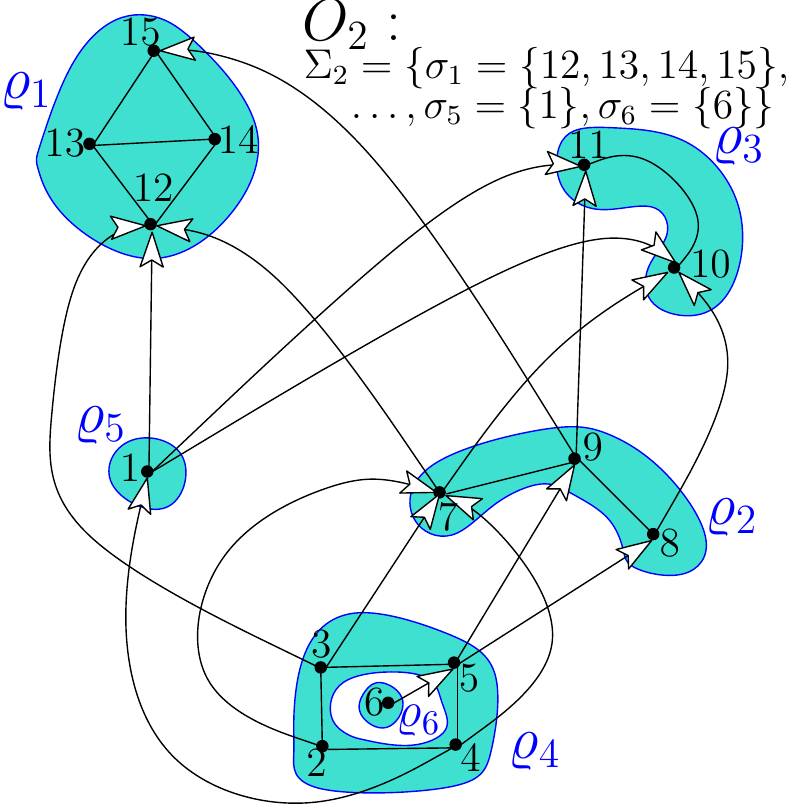}
  \subcaption{}\label{fig:defpao4}
\end{subfigure}
\\
\begin{subfigure}{.38\textwidth}
  \raggedleft
  \subcaption{}\label{fig:defpao5}
  \includegraphics[width=.5\linewidth]{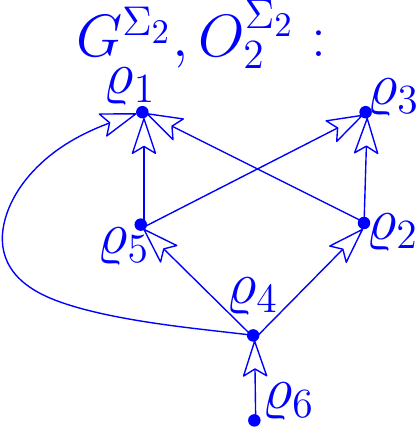}
\end{subfigure}}
\\
\multicolumn{3}{l}{
\begin{subfigure}{\textwidth}
  \raggedright
  \includegraphics[width=\linewidth]{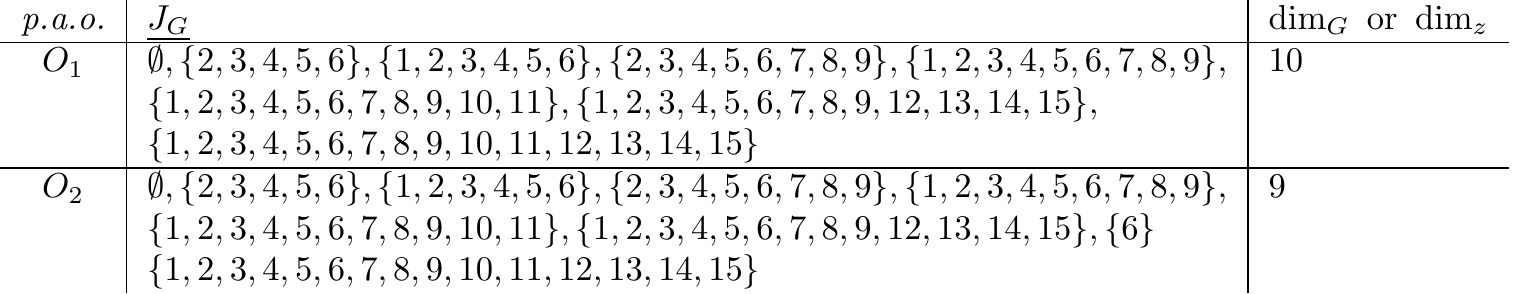}
  \subcaption{}
  \label{fig:defpao6}
\end{subfigure}
}
\end{tabular}

\caption{Examples of \emph{p.a.o.}'s and the order relation $\preceq_z $ of Definition~\ref{defn:complexzon}.} 
\label{fig:defpao}
\end{figure}

\begin{exam}\label{exam:defpao}
In Figure~\ref{fig:defpao}, we present two examples of \emph{p.a.o.}'s, $O_1$ and $O_2$, of a graph $G$ on vertex set
$[15]=\{1,2,\dots,15\}$, such that $O_2\preceq_z O_1$. Figure \ref{fig:defpao1} shows a
connected simple graph $G=G([15,E])$. Figure \ref{fig:defpao2} presents a particular \emph{p.a.o.} $O_1$ of $G$, with associated
connected partition $\Sigma_1$ (each of its blocks represented by closed blue regions), and Figure \ref{fig:defpao3} the $\Sigma_1$-partition graph
$G^{\Sigma_1}$ and its acyclic orientation $O_{1}^{\Sigma_1}$.  Similarly, Figure \ref{fig:defpao4} shows another \emph{p.a.o.} $O_2$ of $G$, with associated
connected partition $\Sigma_2$ (blocks represented by closed blue regions), and Figure \ref{fig:defpao5} the $\Sigma_2$-partition graph
$G^{\Sigma_2}$ and its acyclic orientation $O_{2}^{\Sigma_2}$. Table \ref{fig:defpao6} then offers complete calculations of
$\underline{J_G}(O_1)$, $\underline{J_G}(O_2)$, $\dim_G(O_1)=\dim_z(O_1)$ and $\dim_G(O_1)=\dim_z(O_2)$. Note that since
$\underline{J_G}(O_1)\subsetneq \underline{J_G}(O_2)$, then $O_2\prec_z O_1$. 
\end{exam}

\begin{lem}\label{lem:tech1}
Let $G=G([n],E)$ be a simple graph, and let $a,b\in \R$ and  $c\in\RP$. Consider the function
$F:2^{[n]}\rightarrow \R$ given by $F(\upsigma)=a+b|\upsigma|+c|E(G[\upsigma])|$, $\upsigma\in[n]$. 
Then, for all $\upsigma,\uprho\subseteq[n]$:
\begin{align*}
F(\upsigma)+F(\uprho)\leq F(\upsigma\cap\uprho) + F(\upsigma\cup\uprho).
\end{align*}
Equality holds if and only if $\upsigma\backslash\uprho$ and $\uprho\backslash\upsigma$ are \emph{completely non-adjacent
sets in $G$}, {\it i.e.} 
if and only if $\left\{\{i,j\}\in E:i\in \upsigma\backslash\uprho\text{ and }j\in\uprho\backslash\upsigma\right\}=\emptyset$. 
\end{lem}
\begin{rem}\label{rem:tech1}
In standard combinatorial theory terminology, in Lemma~\ref{lem:tech1}, we say that the function $F$ is \emph{lower semi-modular}
{\bfseries [}\cite{crapofoundations}{\bfseries ]}. 
\end{rem}
\begin{theo}\label{theo:polrealzon}
Let $G=G([n],E)$ be a simple graph with abstract cell complex $\mathscr{Z}_G$, as in 
Definition~\ref{defn:complexzon}. Then, the face complex of the polytope:
\begin{align}\label{eq:zon1}
\CZon{G} :=\\
& \left\{x\in \R^{[n]}:\sum_{i\in[n]}x_i = n+|E|\text{ and }\right.\nonumber\\
& \left.\sum_{i\in\upsigma} x_i\geq |\upsigma|+|E(G[\upsigma])|\text{ for all }\upsigma\subseteq [n]\right\}, \nonumber
\end{align}
is a polytopal complex realization of $\mathscr{Z}_G$.
\end{theo}
\begin{proof2}
Per Lemma~\ref{lem:injectiveord} and for the sake of clarity, in this proof we will think of elements of $\underline{\mathscr{Z}_G}$ as
their images under $\underline{J_G}$. 
 
To begin, an easy verification shows that the point $\frac{1}{2}\cdot\degvec{G}+1$ lives inside 
$\CZon{G}$, so $\CZon{G}$ is non-empty. Also, $\CZon{G}$ is bounded. 

Now, consider a (relatively open) non-empty face $F$ of $\CZon{G}$, and 
let $\mathcal{C}_F$ be the collection of all $\upsigma\subseteq[n]$
such that $\sum_{i\in\upsigma}y_i = |\upsigma| + |E(G[\upsigma])|$ if $y\in F$. 

A first key step in the proof will be to establish 
that $\mathcal{C}_F\in \underline{\mathscr{Z}_G}$. We will do this in a series of sub-steps. Let $\mathcal{C}_{F}^{^\subseteq}$
be the poset on ground set $\mathcal{C}_{F}$ ordered by inclusion.
\begin{enumerate}[label=\bfseries Claim \roman*]
\item\label{claimi1:polrealzon} $\mathcal{C}_F$ is closed under intersections and unions, so $\mathcal{C}_{F}^{^\subseteq}$ is a distributive lattice. 
\par
Let $y\in F$. By definition, both $\emptyset$ and $[n]$ belong to $\mathcal{C}_F$. Let us now  
take $\upsigma,\uprho \in \mathcal{C}_F$ and let us assume that $\upsigma\not\subseteq \uprho, \uprho\not\subseteq \upsigma$.
Then, $\sum_{i\in \upsigma}y_i = |\upsigma| + |E(G[\upsigma])|$ and $\sum_{j\in\uprho}y_j = |\uprho|+|E(G[\uprho])|$, so:
\begin{align*}
|\upsigma\cup\uprho|+|E(G[\upsigma\cup\uprho])|&\leq \sum_{i\in\upsigma\cup \uprho} y_i \\
& = \sum_{i\in\upsigma}y_i+\sum_{j\in\uprho}y_j-\sum_{k\in\upsigma\cap\uprho}y_k \\
& \leq |\upsigma\cup\uprho| + |E(G[\upsigma])|+|E(G[\uprho])|\\
& -|E(G[\upsigma\cap\uprho])|.
\end{align*}
In particular, $|E(G[\upsigma\cap\uprho])|+|E(G[\upsigma\cup\uprho])|\leq |E(G[\upsigma])|+|E(G[\uprho])|$. However, per Lemma~\ref{lem:tech1}:
$$|E(G[\upsigma\cap\uprho])|+|E(G[\upsigma\cup\uprho])| = |E(G[\upsigma])|+|E(G[\uprho])|.$$ 
This implies that $\upsigma\cap\uprho\in\mathcal{C}_F$ and
$\upsigma\cup\uprho \in\mathcal{C}_F$. 
\item\label{claimi2:polrealzon} Let $\emptyset=\upsigma_0\subsetneq\upsigma_1\subsetneq\cdots\subsetneq \upsigma_k = [n]$ be a maximal chain
in $\mathcal{C}_{F}^{^\subseteq}$. Then, $G[\upsigma_{i}\backslash\upsigma_{i-1}]$ is connected for all $i\in[k]$. 
\par
Let $i\in[k]$ and suppose that $G[\upsigma_{i}\backslash\upsigma_{i-1}]$ is disconnected. Let $\uprho_1$ and $\uprho_2$ be two 
completely non-adjacent sets of $G[\upsigma_{i}\backslash\upsigma_{i-1}]$ such that $\uprho_1\cup\uprho_2 = \upsigma_{i}\backslash\upsigma_{i-1}$. 
Then Lemma~\ref{lem:tech1} shows: 
$$|E(G[\upsigma_{i-1}\cup\uprho_1])| + E(G[\upsigma_{i-1}\cup\uprho_2])| = |E(G[\upsigma_i])|+|E(G[\upsigma_{i-1}])|.$$ 
Also, for $y\in F$: 
\begin{align*}
&|\upsigma_{i-1}| + |\upsigma_{i}|+|E(G[\upsigma_{i-1}\cup\uprho_1])| +E(G[\upsigma_{i-1}\cup\uprho_2])| \\
&= |\upsigma_{i-1}\cup\uprho_1| + |E(G[\upsigma_{i-1}\cup\uprho_1])| +|\upsigma_{i-1}\cup\uprho_2|+E(G[\upsigma_{i-1}\cup\uprho_2])|\\
&\leq \sum_{j\in \upsigma_{i-1}\cup\uprho_1}y_j + \sum_{k\in\upsigma_{i-1}\cup\uprho_2}y_k = \sum_{j\in\upsigma_{i-1}}y_j+\sum_{k\in\upsigma_i}y_k\\
& = |\upsigma_{i-1}| + |\upsigma_{i}| + |E(G[\upsigma_{i-1}])|+|E(G[\upsigma_{i}])|.
\end{align*}
This implies that $\upsigma_{i-1}\cup\uprho_1\in\mathcal{C}_F$ and $\upsigma_{i-1}\cup\uprho_2\in\mathcal{C}_F$, contradicting the choice of a maximal
chain in $\mathcal{C}_{F}^{^\subseteq}$. 
\item\label{claimi3:polrealzon} For a chain as in \ref{claimi2:polrealzon}, suppose that there exist $l,m\in\upsigma_i\backslash\upsigma_{i-1}$ with $m\neq l$, $i\in[k]$.
If $\upsigma\in\mathcal{C}_F$, then either $\{m,l\}\cap\upsigma=\emptyset$ or $\{m,l\}\subseteq \upsigma$. 
\par
Suppose on the contrary that for some $\upsigma\in\mathcal{C}_F$, $m\in \upsigma$ but $l\not\in\upsigma$. Then, 
$(\upsigma\cap\upsigma_i)\cup\upsigma_{i-1}\in \mathcal{C}_F$ per~\ref{claimi1:polrealzon}, and $\upsigma_{i-1}\subsetneq
(\upsigma\cap\upsigma_i)\cup\upsigma_{i-1}\subsetneq \upsigma_i$, which contradicts the choice of maximal chain. 
\item\label{claimi4:polrealzon} Per~\ref{claimi3:polrealzon}, every $\upsigma\in\mathcal{C}_F$ is a disjoint union of elements of
the connected partition $\Upsigma := \left\{\upsigma_{i}\backslash\upsigma_{i-1}\right\}_{i\in[k]}$ of $G$. Consider the 
acyclic orientation $O^{\Upsigma}$ of
$G^{\Upsigma}=G^{\Upsigma}(\Upsigma,E^{\Upsigma})$ such that for $e=\{\upsigma_{i}\backslash\upsigma_{i-1},\upsigma_{j}\backslash\upsigma_{j-1}\}\in E^{\Upsigma}$ and $i<j$, 
$O^{\Upsigma}(e)=(\upsigma_{i}\backslash\upsigma_{i-1},\upsigma_{j}\backslash\upsigma_{j-1})$. Then, both $\Upsigma$ and $O^{\Upsigma}$ are
well-defined, {\it i.e.} independent of the choice of maximal chain of $\mathcal{C}_{F}^{^{\subseteq}}$.  
\par
That $\Upsigma$ is well-defined follows from~\ref{claimi1:polrealzon} and~\ref{claimi3:polrealzon}, since distributed lattices are graded. To prove that
$O^{\Upsigma}$ is also consistent, it suffices to check that if $\upsigma\subsetneq \uprho$ for some $\upsigma,\uprho\in\mathcal{C}_F$ such that 
$\upsigma$ and $\uprho\backslash\upsigma$ are adjacent in $G$, then a set $\tau\in\mathcal{C}_F$ with 
$N_G(\upsigma)\cap (\uprho\backslash\upsigma)\subseteq \tau$ must satisfy
that $\upsigma \cap N_G(\uprho\backslash\upsigma)\subseteq \tau$. 
On the contrary, if $\upsigma \cap N_G(\uprho\backslash\upsigma)\not\subseteq \tau$, then both 
$\tau\backslash\upsigma$  and $\upsigma\backslash \tau$ are non-empty and adjacent in $G$. However, since per~\ref{claimi1:polrealzon} 
$\upsigma\cup \tau \in \mathcal{C}_F$, we obtain a contradiction with Lemma~\ref{lem:tech1}. 
\item\label{claimi5:polrealzon} From~\ref{claimi4:polrealzon}, let $O$ be the \emph{p.a.o.} 
of $G$ obtained from $O^{\Upsigma}$. If $\upsigma\in\mathcal{C}_F$, then $\upsigma\in\underline{J_G}(O)$.
\par
This is essentially a corollary to the proof of~\ref{claimi4:polrealzon}. Consider a maximal chain of $\mathcal{C}_{F}^{^{\subseteq}}$ that
contains $\upsigma$. Then, clearly $\upsigma\in\underline{J_G}(O)$.  
\end{enumerate}
For the first part, it remains to prove that if $\upsigma\in\underline{J_G}(O)$, then $\upsigma\in\mathcal{C}_F$. 
This is easy to establish by considering a point $y\in F$. For $\uprho\in \Upsigma$, note that 
$\sum_{i\in \uprho}y_i = |\uprho|+|E(G[\uprho])|+|O[E]\cap ([n]\backslash\uprho \times\uprho)|$. But then: 
\begin{align*}
\sum_{i\in\upsigma}y_i &= \sum_{\uprho \in\Upsigma:\uprho\subseteq \upsigma}|\uprho|+|E(G[\uprho])|+|O[E]\cap ([n]\backslash\uprho \times\uprho)|\\
&=|\upsigma| + |E(G[\upsigma])|,
\end{align*} 
since $\upsigma\in\underline{J_G}(O)$. Hence, $\upsigma\in\mathcal{C}_F$ and $\mathcal{C}_F = \underline{J_G}(O)\in\underline{\mathscr{Z}_G}$.
\par
For the second part, let us take a $\underline{J_G}(O)\in\underline{\mathscr{Z}_G}$ for some
\emph{p.a.o.} $O$ of $G$, and we want to prove that $\underline{J_G}(O) = \mathcal{C}_F$ for some
(relatively-open) non-empty face $F$ of $\CZon{G}$. The first part gives us a clear hint of how to proceed here. 
Let us consider the point $y\in\R^{[n]}$ given by $y_i = |O[E]\cap ([n]\backslash\{i\}\times \{i\})| + 
\frac{1}{2}\cdot |O[E]\cap \left\{e\in E:i\in e\right\}| + 1$ for all $i\in[n]$. Since $O$ is a \emph{p.a.o.} of $G$, then
for $\upsigma\in \underline{J_G}(O)$, we have that $\sum_{i\in\upsigma}y_i = |\upsigma| + |E(G[\upsigma])|$, 
as $O[E]\cap \left\lgroup\left([n]\backslash\upsigma\times \upsigma\right)\cup\left\{\{i,j\}:i\in\upsigma, j\in[n]\backslash\upsigma\right\}\right\rgroup=\emptyset$.
On the contrary, if $\upsigma\not\in \underline{J_G}(O)$, the later set is non-empty and
$\sum_{i\in\upsigma}y_i > |\upsigma| + |E(G[\upsigma])|$. Therefore, $\underline{J_G}(O) = \mathcal{C}_F$ for some face $F$ of $\CZon{G}$.
\par
We have now established how $\underline{\mathscr{Z}_G}$ corresponds to the set of (relatively-open) non-empty faces of 
$\CZon{G}$. Naturally, $\preceq_z$ corresponds to face containment and $\dim_z$ to affine dimension, and the
correctness of these two is an immediate consequence of our correspondence and 
of basic properties of the inequality description of a polytope.
\qed\end{proof2}
\begin{cor}\label{cor:zon}
Let $G=G([n],E)$ be a simple graph with graphical zonotope $\Zon{G}$ and degree vector $\mathbf{d}_G$, where
(see Definition~\ref{defn:convhull} for notation):
$$\Zon{G}:= \displaystyle\sum_{\{i,j\}\in
E}\left[e_i-e_j,e_j-e_i\right].$$ 
Then:
\begin{align}\label{eq:zon}
\CZon{G} = \frac{1}{2}\cdot\Zon{G}+\frac{1}{2}\cdot\mathbf{d}_G+\mathbf{1}.
\end{align}
\end{cor}
\begin{proof2}
As it is known from~\cite{iriartegraphs}, the vertices of $\Zon{G}$ are given by all points of the form 
$x_O =\left(\indeg{i}{O} - \outdeg{i}{O}\right)_{i\in[n]}$, where $O$ is an acyclic orientation
of $G$. Hence, translating $\frac{1}{2}\cdot\Zon{G}$ by $\frac{1}{2}\cdot\mathbf{d}_G + \mathbf{1}$, we obtain that the vertices of
the new polytope are given by all vectors of the form $y_O = \left(\indeg{i}{O}+1\right)_{i\in[n]}$, where $O$ is an acyclic orientation of $G$, but these are precisely
the vertices of $\CZon{G}$. 
\qed\end{proof2}
\begin{defn}\label{defn:transzon}
Let $G=G([n],E)$ be a simple graph with graphical zonotope $\Zon{G}$. From Corollary~\ref{cor:zon}, we will call 
the polytope $\CZon{G}$ the \emph{clean graphical zonotope of $G$}.
\end{defn}

\subsection{One More Degree of Freedom.}\label{subsec:onemoredeg}

\begin{defn}\label{defn:complex1}
Let $G=G([n],E)$ be a connected graph. Define $\mathscr{Y}_G = \left(\underline{\mathscr{Y}_G},\preceq_y,\dim_y\right)$ to be the abstract cell complex with underlying set of cells
$\underline{\mathscr{Y}_G}$, order relation $\preceq_y$, and dimension map $\dim_y$, given by:
\begin{enumerate}
\item $\underline{\mathscr{Y}_G} = \underline{\mathscr{X}_{G}}\cup \left(\midset{[n]}\right)$, where:
$$\underline{\mathscr{X}_{G}}:= \left\{(\upsigma,O):\text{ $\emptyset\neq\upsigma\subseteq [n]$, and $O$ is a \emph{p.a.o.} of 
$G[\upsigma]$}\right\}.$$
\item For $A,B \in \underline{\mathscr{Y}_G}$ with $A\neq B$, we have that $A\preceq_y B$ if and only if one of the following holds:
\begin{enumerate}
\item\label{defn:complex1:c21}  If $A,B\in \midset{[n]}$, then $A\subseteq B$.
\item\label{defn:complex1:c22}  If $A\in \midset{[n]}$ and $B=(\upsigma, O)\in\underline{\mathscr{X}_{G}}$, then  $A\subseteq [n]\backslash \upsigma$. 
\item\label{defn:complex1:c23}  If $A=(\upsigma_0,O_0),B=(\upsigma_1, O_1)\in\underline{\mathscr{X}_{G}}$, then 
$\underline{J_{G[\upsigma_1]}}(O_1)\subseteq \underline{J_{G[\upsigma_0]}}(O_0)$. 
\end{enumerate}
\item For $A\in \underline{\mathscr{Y}_G}$: 
\begin{enumerate}
\item If $A\in \midset{[n]}$, then $\dim_y(A) = |A|-1$.
\item\label{defn:complex1:c33}  If $A = (\upsigma,O)\in\underline{\mathscr{X}_{G}}$, then
$\dim_y(A) = |[n]\backslash \upsigma| + \dim_{G[\upsigma]}(O)$. 
\end{enumerate}
\end{enumerate}
\end{defn}
\begin{defn}\label{defn:convhull}
Let $S$ and $T$ be non-empty subsets of $\R^{[n]}$. The \emph{join of $S$ and $T$} is the
set: 
$$\left[S, T\right] := \{x\in \R^{[n]}: x = \upalpha s+(1-\upalpha) t\text{ for some $\upalpha\in[0,1]$, $s\in S$ and $t\in T$}\}.$$
The \emph{strict join of $S$ and $T$} is the set: 
$$\left(S, T\right) := \{x\in \R^{[n]}: x = \upalpha s+(1-\upalpha) t\text{ for some $\upalpha\in(0,1)$, $s\in S$ and $t\in T$}\}.$$
\end{defn}

\begin{figure}[hb]
\begin{tabular}{cc}
\begin{subfigure}{.5\textwidth}
  \centering
  \includegraphics[width=1.05\linewidth]{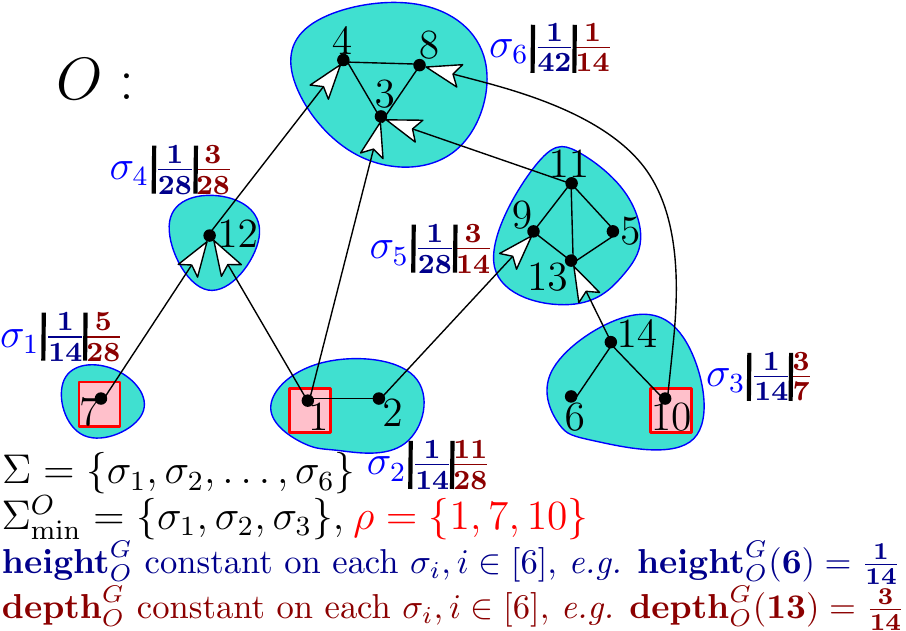}
  \caption{}
  \label{fig:tech2}
\end{subfigure}
&
\begin{subfigure}{.5\textwidth}
  \centering
  \includegraphics[width=1.05\linewidth]{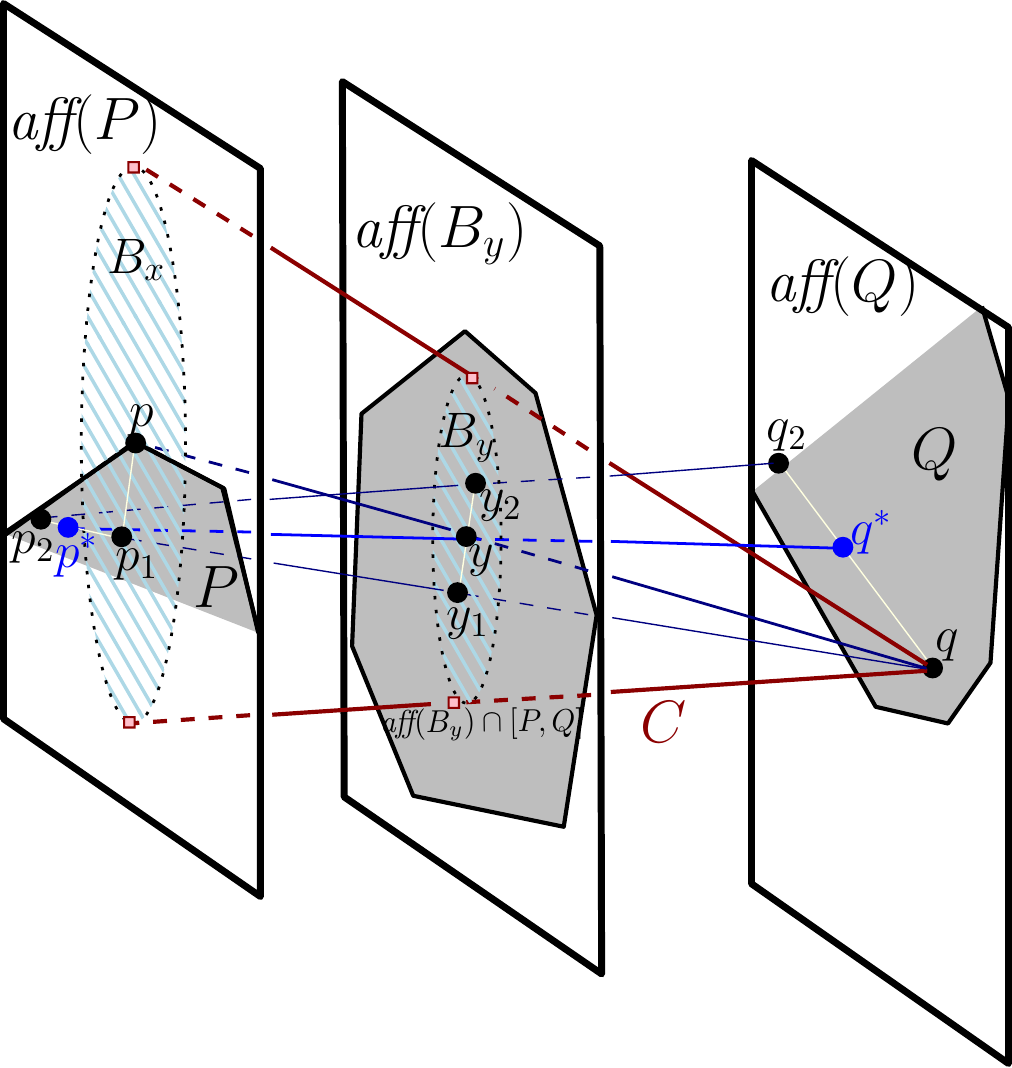}
  \caption{}
  \label{fig:tech1}
\end{subfigure}
\end{tabular}

\caption{Visual aids/guides to the proofs of Proposition~\ref{prop:tech3} ({\bf A}) and Proposition~\ref{prop:tech2}.\ref{prop:tech2:c1} ({\bf B}). 
{\bf A} also offers an example for Definition~\ref{defn:handd}.}\label{fig:tech}
\end{figure}

\begin{prop}\label{prop:tech2}
Let $P$ and $Q$ be $(n-1)$-dimensional polytopes in $\R^{n}$ such that
$\aff{P}$ and $\aff{Q}$ are parallel and disjoint affine hyperplanes. Consider
an open segment $(x,z)$ with $x\in P$, $z\in Q$, and let $y\in \hull{P}{Q}$. Then, the following are true:

\begin{enumerate}[label=\bfseries \roman*]
\item\label{prop:tech2:c1} $y \in \inter{\hull{P}{Q}}$ if and only if there exist $p^{\ast} \in \relint{P}$ and $q^{\ast}\in Q$ such that
$y\in (p^{\ast},q^{\ast})$, if and only if there exist $p^{\ast\ast} \in \relint{P}$ and $q^{\ast\ast}\in \relint{Q}$ such that
$y\in (p^{\ast\ast},q^{\ast\ast})$. 
\item\label{prop:tech2:c2} $(x,z)\subseteq\bound{\hull{P}{Q}}$ if and only if for every $p\in\relint{P}$ and $\upvarepsilon > 0$,  
$z+\upvarepsilon(x-p) \not\in Q$. On the contrary, $(x,z)\subseteq\inter{\hull{P}{Q}}$ if and only if there exists $p\in\relint{P}$ and $\upvarepsilon > 0$,  
such that $z+\upvarepsilon(x-p) \in Q$. 
\item\label{prop:tech2:c3} Let $\uppi_{\aff{P}}:\R^{[n]}\rightarrow\R^{[n]}$ be the projection operator onto the affine hyperplane containing $P$. If 
$\uppi_{\aff{P}}[(x,z)]\cap \ \relint{P}\cap \ \relint{\uppi_{\aff{P}}[Q]}\neq \emptyset$, then $(x,z) \subseteq \inter{\hull{P}{Q}}$. 
\end{enumerate} 
\end{prop}
\begin{proof2}
We will obtain these results in order.
\begin{enumerate}[label=\bfseries \roman*]
\item (See also Figure \ref{fig:tech1}) We prove the ``only if'' direction for both statements. Suppose that $y \in \inter{\hull{P}{Q}}$ and
let $p\in P$ and $q\in Q$ be such that $y\in(p,q)$. Let us assume  
that $p\in \bound{P}$. Take an open $(n-1)$-dimensional ball $B_y\subseteq\inter{\hull{P}{Q}}$ centered at $y$ such that $\aff{B_y}$ is parallel to 
$\aff{P}$ and $\aff{Q}$. Let $C$ be the positive open cone generated by $B_y - q$, and consider 
the affine open cone $q + C$. Then, $B_x := (q + C)\cap \aff{P}$ is an open $(n-1)$-dimensional ball
in $\aff{P}$ such that $p \in\relint{B_x}$. Hence, since $P$ is also $(n-1)$-dimensional, there exists 
some $p_1\in\relint{P}\cap B_x$. Now, let $y_1 = (p_1,q)\cap B_y \in \inter{\hull{P}{Q}}$. Since 
$y_2:= y+(y-y_1) \in B_y\subseteq \inter{\hull{P}{Q}}$, there exist $p_2\in P$ and $q_2\in Q$ such that
$y_2 = (p_2,q_2)\cap B_y$. But then, there exist $p^{\ast}\in(p_1,p_2)\subseteq \relint{P}$ and $q^{\ast}\in (q,q_2)\subseteq Q$
such that $y = (p^{\ast},q^{\ast})\cap B_y$, as we wanted. If $q^{\ast}\in\bound{Q}$, we can now repeat an analogous
construction starting from $q^{\ast}$ and $p^{\ast}$ to find $p^{\ast\ast} \in \relint{P}$ and $q^{\ast\ast}\in \relint{Q}$ such that
$y\in (p^{\ast\ast},q^{\ast\ast})$. 
\item This is a consequence of~\ref{prop:tech2:c1}, and not easy to prove without it. We prove the second statement, which is equivalent
to the first. For the ``if'' direction,
suppose that for some $p\in\relint{P}$ and $\upvarepsilon > 0$, $z+\upvarepsilon(x-p) \in Q$. Take some $y\in(x,z)$ and
consider the line containing both $z+\upvarepsilon(x-p)$ and $y$. For a sufficiently small $\upvarepsilon$, this line
intersects $\aff{P}$ in some $p_1\in\relint{P}$. But then, for a small open ball $B_{p_1}\subseteq\relint{P}$ 
centered at $p_1$ and with $\aff{B_{p_1}}=\aff{P}$, the open set $\shull{B_{p_1}}{z}$ contains $y$ and lies
completely inside $\inter{\hull{P}{Q}}$, so $y\in\inter{\hull{P}{Q}}$.
\par 
For the ``only if'' direction, suppose that $(x,z)\subseteq\inter{\hull{P}{Q}}$ and take $y\in (x,z)$. If $x\in\relint{P}$, then we are done
since $Q$ is also $(n-1)$-dimensional. If $x\in\bound{P}$, from~\ref{prop:tech2:c1}, take
$p\in\relint{P}$, $p\neq x$, and $q\in Q$ with $y\in(p,q)$. But then, $z + \upvarepsilon(y-p) = q \in Q$ for some
$\upvarepsilon > 0$. 
\item Take $p\in\uppi_{\aff{P}}[(x,z)]\cap \ \relint{P}\cap \ \relint{\uppi_{\aff{P}}[Q]}$ and let $\mathbf{p}\neq \mathbf{0}$ be 
a normal to $\aff{P}$. Then, for some $y\in(x,z)$ and real number $\upalpha\neq 0$, 
$y\in(p,p+\upalpha\mathbf{p})$ and $p+\upalpha\mathbf{p}\in\relint{Q}$, 
so~\ref{prop:tech2:c1} shows that $y\in\inter{\hull{P}{Q}}$. Clearly then $(x,z)\subseteq\inter{\hull{P}{Q}}$.
\end{enumerate}
\end{proof2}

\begin{defn}\label{defn:handd}
Let $G=G([n],E)$ be a simple graph, and let $O$ be a \emph{p.a.o.} of $G$ with
connected partition $\Upsigma$ and acyclic orientation $O^{\Upsigma}$ of $G^{\Upsigma}$.
Let us write $\Upsigma_{\min}^{O}$ for the set of elements of $\Upsigma$ that are minimal in $(\Upsigma,\leq_{O^{\Upsigma}})$, and for
$i\in[n]$ with $i\in\uprho\in\Upsigma$, let:
\begin{align*} 
I^{\wedge}_{G}(i,O^{\scriptscriptstyle \Upsigma})=&\{\upsigma\in\Upsigma:\text{$\upsigma\leq_{O^{\scriptscriptstyle \Upsigma}}\uprho$}\}\text{, and}\\
I^{\vee}_{G}(i,O)=&\{j\in[n]:\text{$j\in\upsigma\in\Upsigma$ and $\upsigma\geq_{O^{\scriptscriptstyle \Upsigma}}\uprho$}\}.
\end{align*}
With this notation, we now define certain functions associated to $O$ and $G$, called \emph{height} and \emph{depth}:
\begin{align*}
\height{O}{G},\depth{O}{G}&:[n]\rightarrow \QQ, \\
\height{O}{G}(i) &= \frac{1}{n\cdot\left\vert \Upsigma_{\min}^{O}\cap I^{\wedge}_{G}(i,O^{\scriptscriptstyle \Upsigma})\right\vert}, \\
\depth{O}{G}(i) &= \sum_{j\in I^{\vee}_{G}(i,O)}\height{O}{G}(j). 
\end{align*}
\end{defn}

\begin{exam}\label{exam:handd}
Figure \ref{fig:tech2} exemplifies Definition~\ref{defn:handd} on a particular graph
$G$ on vertex set $[14]=\{1,2,\dots,14\}$, with given \emph{p.a.o.} $O$. Since both
$\height{O}{G}$ and $\depth{O}{G}$ are constant within each element/block of the connected partition
$\Sigma=\left\{\sigma_1=\{7\},\sigma_2=\{1,2\},\sigma_3=\{6,10,14\},\dots,
\sigma_6=\{3,4,8\}\right\}$ associated to $O$, then we present only that common value
for each block in the figure. 
\end{exam}

\begin{prop}\label{prop:tech3}
In Definition~\ref{defn:handd}, let $\upsigma\in\underline{J_G}(O)$ and let $\uprho\subseteq[n]$ intersect
every element of $\Upsigma_{\min}^{O}$ in exactly one point and contain only minimal elements of $O$. Then:
$$1\geq \sum_{i\in\uprho\cap\upsigma}\depth{O}{G}(i) \geq \tfrac{|\upsigma|}{n}.$$
Moreover, if $G$ is connected, then $\sum_{i\in\uprho\cap\upsigma}\depth{O}{G}(i) > \tfrac{|\upsigma|}{n}$ if and only
if $\upsigma\neq[n]$, and whenever this holds, $\sum_{i\in\uprho\cap\upsigma}\depth{O}{G}(i) - \tfrac{|\upsigma|}{n}> \tfrac{1}{n^2}$.
\end{prop}
\begin{rem}\label{rem:tech3}
Figure~\ref{fig:tech2} shows one such choice of a set $\uprho$ in Proposition~\ref{prop:tech3} 
that works for Example~\ref{exam:handd} (in red). 
\end{rem}
\begin{proof2}
The verification is actually a simple double-counting argument using the fact that $\upsigma$ is an order ideal, so we omit it. When $G$ is connected,
if $\upsigma\neq[n]$, then there must exist $i\in[n]\backslash \upsigma$ that is strictly greater in $O$ than some element of $\upsigma$ 
(and hence strictly greater than some element of $\uprho$), again since $\upsigma$ is an order ideal. Clearly, we must have 
$\height{O}{G}(i)>\tfrac{1}{n^2}$.   
\end{proof2}

\begin{theo}\label{theo:acyc1}
Let $G=G([n], E)$ be a connected simple graph with abstract cell complex $\mathscr{Y}_G$ as in Definition~\ref{defn:complex1}.
For $N > 0, N\neq n+|E|$, consider the $(n-1)$-dimensional simplex $N\Updelta=\Conv{Ne_1, Ne_2,\dots, Ne_n}$ in $\R^{[n]}$. If
we let $\mathcal{Y}_G$ be the polytopal complex obtained from the join $ \hull{\CZon{G}}{N\Updelta}$ after removing
the (open) $n$-dimensional cell and the (relatively open) $(n-1)$-dimensional cell corresponding to
$N\Updelta$, then $\mathcal{Y}_G$ is a polytopal complex realization of $\mathscr{Y}_G$.
\end{theo}
\begin{proof2}
Let the faces of $\mathscr{Y}_G$ obtained from $\midset{[n]}$ correspond to the faces of $\bound{N\Updelta}$ 
in the natural way. Also, let the faces of $\underline{\mathscr{X}_{G}}$ of the form $([n],O)$ correspond to the faces of $\CZon{G}$
as in Theorem~\ref{theo:polrealzon}. The result is clearly true for the restriction to this two sub-complexes, so we will concentrate our efforts 
on the remaining cases.
\par
First, for the sake of having a lighter notation during the proof, we will let $\chat{\uprho}=[n]\backslash \uprho$ for any set $\uprho\subseteq[n]$. 
\par
A (relatively open) cell of $\mathcal{Y}_G\backslash(\CZon{G}\cup \bound{N\Updelta})$ can only be obtained as the strict join 
of a cell of $\bound{\CZon{G}}$ and a cell of $\bound{N\Updelta}$, so let us adopt some conventions to refer to this objects.  

\begin{conv}\label{conv:fands}
During the course of the proof, we will let $S$ (or $S_0$) denote a generic non-empty relatively open cell of $N\Updelta$ obtained from $\uprho\subseteq[n]$ (resp. 
$\uprho_0$), and $F$ (or $F_0$) a generic relatively open cell of $\CZon{G}$ with \emph{p.a.o.} $O$ of $G$, associated connected partition $\Upsigma$ of $G$, and
acyclic orientation $O^{\Upsigma}$ of $G^{\Upsigma}$ yielding $O$ (resp. $O_0,\Upsigma_0,O_{0}^{\Upsigma_0}$).
\end{conv}

We argue that we will be done if we can prove the following claim:

\begin{enumerate}[label=\bfseries Claim \roman*]
\item\label{theo:acyc1:c1} {\bf a)} $\shull{F}{S}$ is a cell of $\mathcal{Y}_G$ if and only if {\bf b)} $\uprho\neq [n]$ and
$\uprho$ is a non-empty union of
elements from the set $\left\{\upsigma\in\Upsigma:\text{$\upsigma$ is maximal in $(\Upsigma,\leq_{O^{\Upsigma}})$}\right\}$. When this equivalence is established, then
we will let $\shull{F}{S}$ correspond to the pair $(\chat{\uprho},
O\vert_{\scriptscriptstyle \chat{\uprho}})\in \underline{\mathscr{X}_{G}}$, where
$O\vert_{\scriptscriptstyle \chat{\uprho}}$ denotes the restriction of $O$ to 
$E(G[\chat{\uprho}])$. 
\end{enumerate}

Indeed, assume that~\ref{theo:acyc1:c1} holds. Then, under the stated correspondence
of ground sets of cells, all elements of $\underline{\mathscr{X}_{G}}$ are uniquely accounted for as cells of $\mathcal{Y}_G$. This is true 
for $\CZon{G}$ clearly,  and for the remaining cases since for
any choice of $\upsigma_{1}\subsetneq[n]$, $\upsigma_{1}\neq\emptyset$,  and of \emph{p.a.o.} $O_{1}$ of $G[\upsigma_{1}]$, we can always extend uniquely $O_{1}$ 
to a \emph{p.a.o.} of $G$ in which all the elements of $\chat{\upsigma_{1}}$ are maximal. 
\par
Secondly, we verify that $\preceq_y$ corresponds to face containment in $\mathcal{Y}_G$. 
Suppose that 
$\shull{F_0}{S_0}$ and $\shull{F}{S}$ are cells of $\mathcal{Y}_G$.
Then, 
$\overline{\shull{F_0}{S_0}}\subseteq \overline{\shull{F}{S}}$ if and only if 
$\overline{F_0}\subseteq \overline{F}$ and $\overline{S_0}\subseteq\overline{S}$, if and only if
$\underline{J_G}(O)\subseteq \underline{J_G}(O_0)$ and  $\uprho_0\subseteq \uprho$. 
Now, assuming~\ref{theo:acyc1:c1}, the last statement is true if and only if
$\underline{J_{G[{\scriptscriptstyle \chat{\uprho}}]}}(O\vert_{\scriptscriptstyle \chat{\uprho}})
\subseteq\underline{J_{G[{\scriptscriptstyle \chat{\uprho_0}}]}}(O_0\vert_{\scriptscriptstyle \chat{\uprho_0}})$: 
\begin{itemize}
\item[] The difficult part here is the
``if'' direction. Clearly, $\uprho_0\subseteq\uprho$. Since $\uprho$ is a union of elements of $\Upsigma$ that are maximal in $(\Upsigma,\leq_{O^{\Upsigma}})$,  
then $\underline{J_G}(O)\Big\backslash \underline{J_{G[{\scriptscriptstyle \chat{\uprho}}]}}(O\vert_{\scriptscriptstyle \chat{\uprho}})$
consists of ideals of $O$ whose intersection with $\uprho$ are non-empty unions of the connected components of $G[\uprho]$. But then, as
$\chat{\uprho}\in \underline{J_{G[{\scriptscriptstyle \chat{\uprho}}]}}(O\vert_{\scriptscriptstyle \chat{\uprho}})\subseteq\underline{J_{G[{\scriptscriptstyle \chat{\uprho_0}}]}}(O_0\vert_{\scriptscriptstyle \chat{\uprho_0}})
\subseteq  \underline{J_G}(O_0)$, these must also be ideals of $\underline{J_G}(O_0)$. 
\end{itemize}
The analogous verification pertaining to faces 
in $\underline{\mathscr{X}_{G}}$ of the form $([n],O)$, or
corresponding to Definition~\ref{defn:complex1:c21}-\ref{defn:complex1:c22}, is now a straightforward
application of the same ideas, so we omit it here.
\par
The correctness of $\dim_y$ will be established in~\ref{cc3}, so indeed if~\ref{theo:acyc1:c1} holds, the 
statement of the Theorem follows.
\par
Let us now begin with our proof of~\ref{theo:acyc1:c1}, which consists of
three main steps.   
\begin{enumerate}[label=\bfseries Claim i.\arabic*]
\item\label{cc1} Let $F$ and $S$ satisfy the conditions of~\ref{theo:acyc1:c1}{\bf .b)}. 
Then: 
$$\shull{F}{S}\subseteq \bound{\hull{\CZon{G}}{N\Updelta}}.$$
\par
Let $x\in F$ and $z\in S$. We must have that $O\neq O_{\scriptscriptstyle\emph{trivial}}$ here. Now, since $G$ is connected, there exists
$\upsigma\in\Upsigma$ that is minimal but not maximal in $(\Upsigma,\leq_{O^{\Upsigma}})$. Hence, $\upsigma\cap\uprho=\emptyset$ and
moreover, $\upsigma\in \underline{J_{G}}(O)$. But then, by the inequality description of $\CZon{G}$, for any
$p\in\relint{\CZon{G}}$, $\sum_{i\in \upsigma}p_i > |\upsigma|+|E(G[\upsigma])| = \sum_{j\in\upsigma}x_i$, and
$x-p$ must have a negative entry in $\upsigma$. Therefore, $z+\upvarepsilon(x-p)\not\in N\Updelta$ for all $\upvarepsilon>0$ and 
Proposition~\ref{prop:tech2}.\ref{prop:tech2:c2} shows that $(x,z)\subseteq\bound{\hull{\CZon{G}}{N\Updelta}}$. 

\item\label{cc2} Let $F_0$, $O_0$, $\Upsigma_0$, $O_{0}^{\Upsigma_0}$, $S_0$, $\uprho_0$ be as in Convention~\ref{conv:fands}. Then, there exist $F$, $O$, $\Upsigma$, 
$O^{\Upsigma}$, $S$, $\uprho$ also as
in Convention~\ref{conv:fands}, such that $\uprho$ is a union of elements of the set $\{\upsigma\in\Upsigma:\text{$\upsigma$ is maximal in $(\Upsigma,\leq_{O^{\Upsigma}})$}\}$ and
$(F_0,S_0)\subseteq (F,S)$. 
\par
(See Figures~\ref{fig:help1} and~\ref{fig:help2} for a particular example of the objects and setting considered during this proof) Let: 
$$\Upsigma_{0,\uprho_0}
:=\left\{\upsigma\in\Upsigma_0:\text{If $\varrho\in\Upsigma_0$ and $\varrho\leq_{O_{0}^{\scriptscriptstyle \Upsigma_0}}\upsigma$, then $\varrho\cap \uprho_0=\emptyset$}\right\}.$$ 
Then, define:
$$\upsigma_{\scriptscriptstyle 0} := \displaystyle\bigcup_{\upsigma\in \Upsigma_{0,\uprho_0}}\upsigma.$$ 
If $G[\chat{\upsigma_{\scriptscriptstyle 0}}] = G[\upsigma_1]+\dots+G[\upsigma_{k}]$ is the decomposition of 
$G[\chat{\upsigma_{\scriptscriptstyle 0}}]$ into its
connected components, we will let $\Upsigma = \Upsigma_{0,\uprho_0}\cup\{\upsigma_1,\dots,\upsigma_k\}$. We will use the acyclic orientation 
$O^{\Upsigma}$ of $G^{\Upsigma}$ obtained from the two conditions
{\bf 1)} $O^{\Upsigma}\vert_{\scriptscriptstyle \Upsigma_{0,\uprho_0}} = 
O_{0}^{\Upsigma_0}\vert_{\scriptscriptstyle \Upsigma_{0,\uprho_0}}$ and {\bf 2)} 
$\upsigma_1,\cdots,\upsigma_k$ are maximal in $(\Upsigma,\leq_{O^{\Upsigma}})$. The \emph{p.a.o.} $O$ is now obtained from $O^{\Upsigma}$, and let $F$ be associated to $O$ and 
$S$ be obtained from $\chat{\upsigma_{\scriptscriptstyle 0}}= \upsigma_1\cup\dots\cup\upsigma_k$. We now prove that $(F_0,S_0)\subseteq (F,S)$. 
\par
Since $\overline{(F_0,S_0)}\subseteq\overline{(F,S)}$, it is enough to find $x\in \overline{F_0}$ and $z\in \overline{S_0}$ such that $(x,z)\subseteq (F,S)$, so this is precisely 
what we will do. To begin,  
we note that for $i\in[k]$, the restriction $O_i := O_{0}\vert_{\scriptscriptstyle \upsigma_i}$ is a \emph{p.a.o.} of $G_i:=G[\upsigma_i]$, so we will let 
$\Upsigma_i$ be the connected partition of
$G_i$ and $O^{\Upsigma_i}_{i}$ the acyclic orientation of $G_{i}^{\Upsigma_i}$ associated to $O_i$; moreover, 
we note that $\uprho_0$ intersects every element of $\Upsigma_i$ minimal in $(\Upsigma_i,\leq_{O^{\Upsigma_i}_{i}})$. Hence, let us select $\varrho_{0}\subseteq\uprho_0$ 
so that for every $i\in[k]$, $\varrho_{0}$ intersects every element of $\Upsigma_i$ minimal in 
$(\Upsigma_i,\leq_{O_{i}^{\Upsigma_i}})$ in exactly one point and so that $\varrho_0\cap\upsigma_i$ contains only
minimal elements in $O_i$. Now, take any $x\in F_0$ and let:
$$z=\tfrac{N}{k}\sum_{i\in[k]}\sum_{j\in\varrho_0\cap\upsigma_i}\depth{O_i}{G_i}(j)\cdot e_j \in \overline{S_0}.$$ 
We will make use of the technique
of Proposition~\ref{prop:tech2}.\ref{prop:tech2:c2} to prove that $(x,z)\in(F,S)$, so for that we need to consider a point in $S$, which we select as:
$$s = \tfrac{N}{k}\sum_{i\in[k]}\tfrac{1}{|\upsigma_i|}\sum_{j\in\upsigma_i}e_j \in S.$$
For $i\in[k]$, if we consider a $\varrho_i\in\underline{J_{G_i}}(O_i)$ with $\varrho_i\neq\upsigma_i$, Proposition~\ref{prop:tech3} gives us:
\begin{align*}
\sum_{j\in\varrho_i}(z-s)_j &= \tfrac{N}{k}\left(\sum_{j\in \varrho_i\cap\varrho_0}\depth{O_i}{G_i}(j)\right)-\tfrac{N}{k}\cdot\tfrac{|\varrho_i|}{|\upsigma_i|}\\
&> \tfrac{N}{k}\left(\tfrac{|\varrho_i|}{|\upsigma_i|} + \tfrac{1}{|\upsigma_i|^2} - \tfrac{|\varrho_i|}{|\upsigma_i|}\right) = \tfrac{N}{k\cdot|\upsigma_i|^2}\\
&>0.
\end{align*}
Hence, for a sufficiently small $\upvarepsilon>0$, $x+\upvarepsilon(z-s) \in F$, so for each $y\in(x,z)$ we can find $x'\in F$ and $s'\in S$ such that
$y\in(x',z')$. That implies $(F_0,S_0)\subseteq (F,S)$. 

\begin{figure}[ht]
\begin{tabular}{cc}
\begin{subfigure}{.5\textwidth}
  \centering
  \includegraphics[width=1\linewidth]{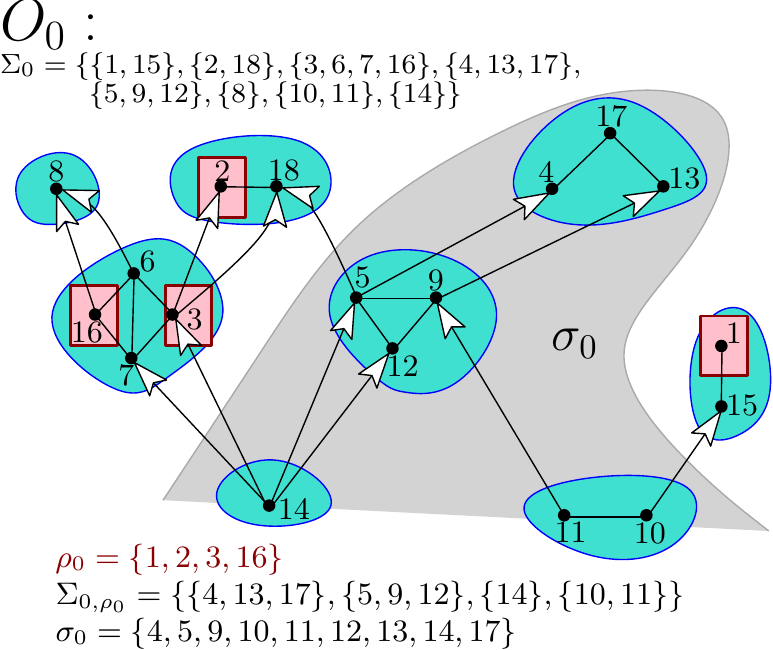}
  \caption{}
  \label{fig:help1}
\end{subfigure}
&
\begin{subfigure}{.5\textwidth}
  \centering
  \includegraphics[width=1\linewidth]{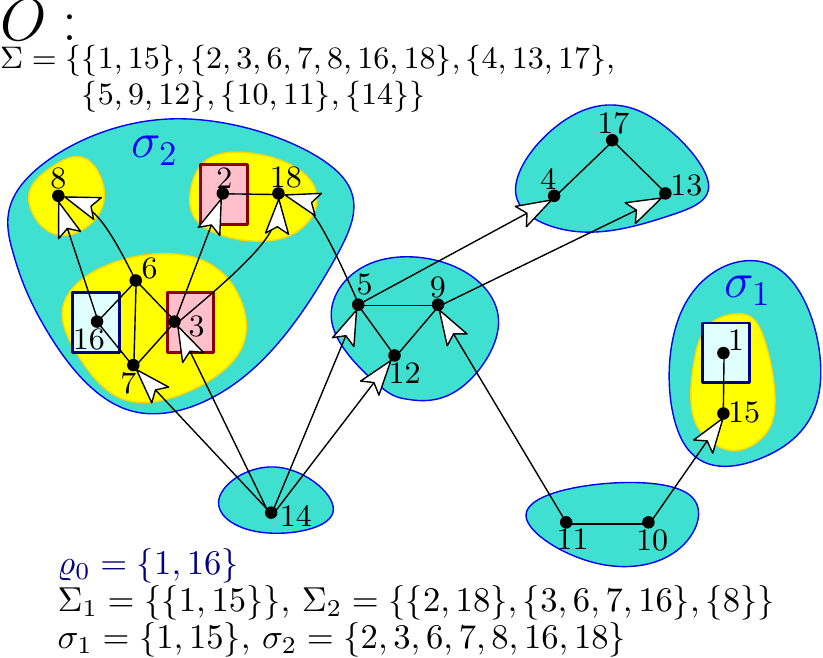}
  \caption{}
  \label{fig:help2}
\end{subfigure}
\end{tabular}

\caption{An example to the proof of \ref{cc2}}\label{fig:help} in Theorem~\ref{theo:acyc1}. 
\end{figure}

\item\label{cc3} Let both $F,S$ and $F_0,S_0$ satisfy the conditions of~\ref{theo:acyc1:c1}{\bf .b)}. Then, $(F,S)\cap(F_0,S_0)\neq\emptyset$ if and only
if $F=F_0$ and $S=S_0$. Moreover, $(F,S)$ is a face of $\mathcal{Y}_G$ and 
$\dim_{\scriptscriptstyle\emph{aff}} \left\langle(F,S)\right\rangle = |\uprho| + \dim_{G[\chat{\uprho}]}(O\vert_{\scriptscriptstyle \chat{\uprho}})$ 
(similarly for $(F_0,S_0)$).  
\par
Let $\upalpha\in(0,1)$ and consider the polytope $P_{\upalpha} = \{x\in\R^{[n]}:\sum_{i\in[n]}x_i=\upalpha(n+|E|)+(1-\upalpha)N\}\cap \hull{\CZon{G}}{N\Updelta}$.  
Every $x\in P_{\upalpha}$ satisfies the inequalities $\sum_{i\in[n]}x_i = (1-\upalpha)N + \upalpha(n+|E|)$ and
$\sum_{i\in\upsigma}x_i\geq \upalpha(|\upsigma|+|E(G[\upsigma])|)$ for all $\upsigma\subsetneq[n]$, $\upsigma\neq\emptyset$. 
Per~\ref{cc1} and~\ref{cc2}, the set $(F,S)\cap P_{\upalpha}$ can be characterized by the condition that it contains all the points $x\in P_{\upalpha}$ 
which, among those inequalities, satisfy the and only the following equalities:
\begin{align}
\sum_{i\in[n]}x_i &= (1-\upalpha)N + \upalpha(n+|E|)\text{ and }\label{eq:affball1}\\
\sum_{i\in\upsigma}x_i &= \upalpha(|\upsigma|+|E(G[\upsigma])|)\label{eq:affball2},\\
&\text{for all $\upsigma\in\underline{J_{G[\chat{\uprho}]}}(O\vert_{\scriptscriptstyle \chat{\uprho}})$, $\upsigma\neq\emptyset$}.\nonumber
\end{align}
This observation proves the first statement. 
\par
For the second statement, we assume without loss of generality that $N>n+|E|$ and select
generic coefficients $\upbeta_{\upsigma}\in \RP$ with $\upsigma\in{\scriptstyle\underline{J_{G[\chat{\uprho}]}}(O\vert_{\scriptscriptstyle \chat{\uprho}})\backslash}{\scriptscriptstyle\{\emptyset\}}$, 
such that: 
$$\sum_{\upsigma\in\underline{J_{G[\chat{\uprho}]}}(O\vert_{\scriptscriptstyle \chat{\uprho}})\backslash{\scriptscriptstyle\{\emptyset\}}}\upbeta_{\upsigma}(|\upsigma|+|E(G[\upsigma])|)=N-(n+|E|).$$ 
The linear functional, 
\begin{align}\label{eq:linfunc}
f:=\sum_{i\in[n]}e_{i}^{\ast}+\sum_{\upsigma\in\underline{J_{G[\chat{\uprho}]}}(O\vert_{\scriptscriptstyle \chat{\uprho}})\backslash{\scriptscriptstyle\{\emptyset\}}}
\upbeta_{\upsigma}\cdot\sum_{j\in\upsigma}e_{j}^{\ast},
\end{align}
satisfies that, for $x\in P_{\upalpha}$, 
$$f(x)\geq (1-\upalpha)N + \upalpha(n+|E|)+\upalpha\left(N-(n+|E|)\right)= N.$$ 
By the proof of the first claim, this inequality is tight if and only if $x\in\overline{(X,S)\cap P_{\upalpha}}=\overline{(X,S)}\cap P_{\upalpha}$. 
Moreover, since this minimum is independent of $\upalpha$, the
linear functional $f$ is minimized in $\hull{\CZon{G}}{N\Updelta}$ exactly at $\overline{(X,S)}$. If $N<n+|E|$, we must select negative coefficients and consider
instead the maximum of the linear functional in question, analogously. 
\par
For the third statement, we simply note that an open ball in the affine space determined by all $x\in\R^{[n]}$ satisfying Equalities~\ref{eq:affball1}-\ref{eq:affball2} can be
easily (but tediously) found inside $(F,S)$. Hence, $\dim_{\scriptscriptstyle\emph{aff}} \left\langle(F,S)\right\rangle  = |\uprho| + \dim_{G[\chat{\uprho}]}(O\vert_{\scriptscriptstyle \chat{\uprho}})$.
\end{enumerate}
\par
\qed\end{proof2}

\begin{defn}\label{defn:complex2}
Let $G=G([n],E)$ be a connected simple graph. Let $\mathscr{X}_{G}^{\ast}=(\underline{\mathscr{X}_{G}},\preceq_x,\dim_x)$ be the abstract cell complex dual to
$\mathscr{X}_{G}$ in Definition~\ref{defn:complex1}. Hence, for all $(\upsigma_0,O_0),(\upsigma_1, O_1)\in\underline{\mathscr{X}_{G}}$:
\begin{enumerate} 
\item\label{defn:complex2:c1} $(\upsigma_0,O_0)\preceq_x (\upsigma_1, O_1)$ if and only if $\underline{J_{G[\upsigma_0]}}(O_0)\subseteq \underline{J_{G[\upsigma_1]}}(O_1)$, and 
\item $\dim_x (\upsigma_0,O_0) = |\upsigma_0| - 1 -  \dim_{G[\upsigma_0]}(O_0)$.
\end{enumerate}
\end{defn}

\begin{theo}\label{theo:acyc2}
Let $G=G([n], E)$ be a connected simple graph with abstract cell complex $\mathscr{X}_{G}^{\ast}$ as in Definition~\ref{defn:complex2}.
Then, the polytopal complex $\mathcal{X}_{G}$ obtained from all faces of the intersection $\Arr{G}\cap \Updelta$ inside $\R^{[n]}$ is a polytopal
complex realization of $\mathscr{X}_{G}^{\ast}$, where $\Arr{G}$ is the \emph{graphical arrangement of $G$} and $\Updelta=\Conv{e_1, e_2,\dots, e_n}$:
$$\Arr{G}:=\{x\in\R^{[n]}:x_i-x_j=0\text{ , $\forall$ $\{i,j\}\in E$}\}.$$
\end{theo}
\begin{proof2}
From Theorems~\ref{theo:polrealzon}-\ref{theo:acyc1}, and letting $N\rightarrow\infty$ in Equation~\ref{eq:linfunc}, 
we know that the relatively open cone $C_{\scriptscriptstyle (\upsigma,O)}^{+}$ in the totally non-negative part of the normal fan of
the polytope $\hull{\CZon{G}}{N\Updelta}$ that corresponds to a cell $(\upsigma,O)\in\underline{\mathscr{X}_{G}}$, is given by:
$$C_{\scriptscriptstyle (\upsigma,O)}^{+} = 
\emph{span}_{\RP}\left\{\sum_{i\in\uprho}e_i:\uprho\in\underline{J_{G[\upsigma]}}(O)\backslash{\scriptscriptstyle \{\emptyset\}} \right\}.$$
Hence, since the affine dimension of the corresponding dual cell in $\mathcal{Y}_G$ is $|[n]\backslash \upsigma|+\dim_{G[\upsigma]}(O)$, then
$\dim_{\scriptscriptstyle\emph{aff}}\left\langle C_{\scriptscriptstyle (\upsigma,O)}^{+}\right\rangle 
=n-|[n]\backslash \upsigma|+\dim_{G[\upsigma]}(O)=|\upsigma|-\dim_{G[\upsigma]}(O)$ and
so $\dim_{\scriptscriptstyle\emph{aff}}\left\langle C_{\scriptscriptstyle (\upsigma,O)}^{+}
\cap\Updelta\right\rangle = |\upsigma|-1-\dim_{G[\upsigma]}(O)$, since
$C_{\scriptscriptstyle (\upsigma,O)}^{+}\subseteq \emph{span}_{\RNN}\{e_1,\dots,e_n\}$. Tangentially, we can also express 
$C_{\scriptscriptstyle (\upsigma,O)}^{+}$ more compactly 
by means of its positive basis as: 
$C_{\scriptscriptstyle (\upsigma,O)}^{+} 
= \emph{span}_{\RP}\left\{\sum_{i\in\uprho}e_i:\text{$\uprho\in\underline{J_{G[\upsigma]}}(O)\backslash{\scriptscriptstyle \{\emptyset\}}$ 
and $G[\uprho]$ is connected}\right\}$. 
\par
Now, the intersection, 
$$\Arr{G[\upsigma]}^{+}=\Arr{G}\cap\{x\in\R^{[n]}:\text{$x_i=0$ if $i\in [n]\backslash\upsigma$, and $x_j>0$ if $j\in\upsigma$}\}\text{ is,}$$ 
as suggested by our choice of 
notation, equal to the totally positive part of the graphical arrangement of $G[\upsigma]$, regarding here $\R^{\upsigma}$ as a subspace of $\R^{[n]}$.
Per Theorem~\ref{theo:polrealzon}, since $\Arr{G[\upsigma]}$ is precisely the normal fan of $\CZon{G[\upsigma]}$, and
$\Arr{G[\upsigma]}^{+}$ the totally positive part of this fan, we know that the relatively open cones of $\Arr{G[\upsigma]}^{+}$ correspond
to the \emph{p.a.o.}'s of $G[\upsigma]$. From the description of the cells of $\CZon{G[\upsigma]}$, the cone $C_{\scriptscriptstyle (\upsigma,O)}^{+}$
is exactly the cone in $\Arr{G[\upsigma]}^{+}$ normal to the cell of $\CZon{G[\upsigma]}$ corresponding to $O$. This establishes the correspondence
between cells of $\Arr{G}\cap \Updelta$ and elements of $\underline{\mathscr{X}_{G}}$, since we can go both ways in this discussion. 
\par
Using the same lens to regard cells of $\Arr{G}\cap \Updelta$, the correctness of Definition~\ref{defn:complex2}.\ref{defn:complex2:c1} now follows from
the analogous verification done in Theorem~\ref{theo:acyc1}, by a standard result on normal fans of polytopes, namely, the duality of face containment.   
\par
\qed\end{proof2}

\section{Two ideals for acyclic orientations.}\label{sec:twoideals}

\begin{defn}\label{defn:notation}
Let $G=G([n],E)$ be a simple graph.
\begin{enumerate} 
\item For an orientation $O$ of $G$ and 
for every $i\in[n]$, let:
\begin{align*}
\indeg{G}{O}(i)&:= \left\vert\left\{(j,i)\in O[E]:j\in[n]\right\}\right\vert,\\
\outdeg{G}{O}(i)&:= \left\vert\left\{(i,j)\in O[E]:j\in[n]\right\}\right\vert,\\
\nooutdeg{G}{O}(i)&:=\left\vert\left\{e\in O[E]:\text{either $e=(j,i)$ or $e=\{i,j\}$, $j\in[n]$}\right\}\right\vert,
\end{align*} 
where we denote the respective associated vectors in $\R^{[n]}$ as $\indegvec{G}{O}$, $\outdegvec{G}{O}$, and
$\nooutdegvec{G}{O}$. 

\item For $\sigma\subseteq[n]$ with $\sigma\neq\emptyset$, define $\mathbf{1}_\sigma:=\sum_{i\in\sigma}e_i\in\R^{[n]}$, further writing 
$\mathbf{1}:=\mathbf{1}_{[n]}$. Let now, for every $i\in[n]$:
\begin{align*}
\degin{G}{\sigma}(i)&:= \left\{\begin{array}{ll} \left\vert\left\{\{i,j\}\in E:j\in\sigma\right\}\right\vert&\text{if $i\in\sigma$,}\\ 0&\text{otherwise,}\end{array}\right.\\
\degout{G}{\sigma}(i)&:=\left\{\begin{array}{ll} \left\vert\left\{\{i,j\}\in E:j\in[n]\backslash\sigma\right\}\right\vert &\text{if $i\in\sigma$,}\\0&\text{otherwise,}\end{array}\right.\\
\end{align*} 
and denote the respective associated vectors of $\R^{[n]}$ as
$\deginvec{G}{\sigma}$ and $\degoutvec{G}{\sigma}$.
\end{enumerate}
\end{defn}

\begin{rem}\label{rem:notationalg}
During this section, we will follow the notation and definitions of~\cite{millersturmfels}, Chapters~1,4,5,6 and 8, in particular,
those pertaining to \emph{labelled polytopal cell complexes}. We refer the reader to this standard reference on the subject
for further details. Some key conventions worth mentioning here are:
\begin{enumerate}
\item The letter $\kfield$ will denote an infinite field.
\item For $\mathbf{a}:=(a_1,a_2,\dots,a_n)\in\NN^{[n]}$, 
$\mathbf{\mathfrak{m}}^{\mathbf{a}}:=\left\langle x_{i}^{a_i}:i\in[n]\right\rangle$ is the ideal of
$\kfield[x_1,\dots,x_n]$ associated to $\mathbf{a}$.
 \end{enumerate}
\end{rem}

\begin{defn}\label{defn:idealacy}
Let $G=G([n],E)$ be a connected simple graph. The \emph{ideal $A_G$ of acyclic orientations of $G$ }is the monomial ideal of 
$\kfield[x_1,\dots,x_n]$ minimally generated as: 
$$A_G:=\left\langle \mathbf{x}^{\indegvec{G}{O}+\mathbf{1}}=\prod_{i\in[n]}x_{i}^{\indeg{G}{O}(i)+1}: \text{$O$ is an acyclic orientation of $G$}\right\rangle.$$
\end{defn}

\begin{defn}\label{defn:idealtree}
Let $G=G([n],E)$ be a connected simple graph. The \emph{tree ideal $T_G$  of $G$ }is the monomial ideal of 
$\kfield[x_1,\dots,x_n]$ minimally generated as: 
$$T_G:=\left\langle \mathbf{x}^{\degoutvec{G}{\upsigma}+\mathbf{1}_{\upsigma}}
=\prod_{i\in \upsigma}x_{i}^{\degout{G}{\upsigma}(i)+1}: \text{$\upsigma\in\upset{[n]}$ and $G[\upsigma]$ is connected}\right\rangle.$$
\end{defn}

\begin{defn}\label{defn:alexander}
Given two vectors $\mathbf{a},\mathbf{b}\in \NN^{[n]}$ with $\mathbf{b}\preceq\mathbf{a}$ ($b_i\leq a_i$ for all $i\in[n]$), let
$\mathbf{a}\backslash\mathbf{b}$ be the vector whose $i$-th coordinate is:
$$a_i\backslash b_i = \left\{\begin{matrix} a_i + 1 - b_i &\text{if $b_i\geq 1$,}\\ 0 &\text{if $b_i=0$.}\end{matrix}\right.$$
If $I$ is a monomial ideal whose minimal generators all divide $\mathbf{x}^{\mathbf{a}}$, then the \emph{Alexander dual of $I$ with respect to $\mathbf{a}$} is:
$$I^{[\mathbf{a}]}:=\bigcap\{\mathbf{\mathfrak{m}}^{\mathbf{a}\backslash\mathbf{b}}:\text{$\mathbf{x}^{\mathbf{b}}$ is a minimal generator of $I$}\}.$$
\end{defn}

\begin{theo}\label{theo:alexander}
Let $G=G([n],E)$ be a simple connected graph. Then, the ideals $A_G$ and $T_G$ of Definitions~\ref{defn:idealacy}-\ref{defn:idealtree} are Alexander
dual to each other with respect to $\degvec{G}+\mathbf{1}$, so $A_{G}^{[\degvec{G}+\mathbf{1}]} = T_G$ and $T_{G}^{[\degvec{G}+\mathbf{1}]} = A_G$. 
\end{theo}
\begin{proof2}
It is enough to prove one of these two equalities, so we will prove that $A_{G}^{[\degvec{G}+\mathbf{1}]} = T_G$. Take some $\upsigma\in\upset{[n]}$ such that 
$G[\upsigma]$ is connected and consider the minimal generator of $T_G$ given by $\mathbf{x}^{\degoutvec{G}{\upsigma}+\mathbf{1}_{\upsigma}}=
\prod_{i\in \upsigma}x_{i}^{\degout{G}{\upsigma}(i)+1}$. We will verify that $\mathbf{x}^{\degoutvec{G}{\upsigma}
+\mathbf{1}_{\upsigma}}\in \mathbf{\mathfrak{m}}^{\left(\degvec{G}+\mathbf{1}\right)\backslash\mathbf{b}}$
for every minimal generator $\mathbf{x}^{\mathbf{b}}$ of $A_G$. Select an acyclic orientation $O$ of $G$ and let 
$\mathbf{x}^{\indegvec{G}{O}+\mathbf{1}}=\prod_{i\in[n]}x_{i}^{\indeg{G}{O}(i)+1}$ be the minimal generator of $A_G$ associated to $O$. If we 
take $m\in\upsigma$ to be maximal in $([n],\leq_O)$ among all elements of $\upsigma$, so that $i\geq_{O} m$ and $i\in\upsigma$ imply $i=m$, then 
$\left\lgroup\degnovec{G}(m) + 1\right\rgroup\big\backslash \left\lgroup\indeg{G}{O}(m) + 1\right\rgroup = \outdeg{G}{O}(m)+1
\leq |N_G(m)\backslash \upsigma|+1=\degout{G}{\upsigma}(m)+1$. 
Hence, 
\begin{align*}
\mathbf{x}^{\degoutvec{G}{\upsigma}+\mathbf{1}_{\upsigma}}&\in\left\langle x^{\degout{G}{\upsigma}(m)+1}_m\right\rangle\subseteq
 \left\langle x^{\left(\degnovec{G}(m) + 1\right)\big\backslash \left(\indeg{G}{O}(m) + 1\right)}_m\right\rangle\\
 &\subseteq
 \mathbf{\mathfrak{m}}^{\left(\degvec{G}+\mathbf{1}\right)\backslash\left(\indegvec{G}{O}+\mathbf{1}\right)}.
 \end{align*}
This proves that $T_G\subseteq A_{G}^{[\degvec{G}+\mathbf{1}]}$.
\par
Now, consider a monomial $\mathbf{x}^{\mathbf{b}}\not\in T_G$ with $\mathbf{0}\prec \mathbf{b}$ (so $b_i>0$ for some $i\in[n]$). 
Then, for every $\upsigma\in\upset{[n]}$ there exists $i\in\upsigma$ such that
$b_i<\degout{G}{\upsigma}(i)+1$, noting here that the condition on $G[\upsigma]$ being connected can be dropped. Hence, consider a bijective labeling $f:[n]\rightarrow [n]$ of the vertices of
$G$ such that $b_{f^{-1}(i)}<\degout{G}{f^{-1}[1,i]}({\scriptstyle f^{-1}(i)})+1$ for all $i\in[n]$. If we let $O$ be the acyclic orientation of $G$ such that for every $e=\{i,j\}\in E$,
$O(e)=(i,j)$ if and only if $f(i)<f(j)$, then for all $i\in[n]$, $b_{f^{-1}(i)}<\degout{G}{f^{-1}[1,i]}({\scriptstyle f^{-1}(i)})+1=\outdeg{G}{O}({\scriptstyle f^{-1}(i)})+1= 
\left\lgroup\degnovec{G}({\scriptstyle f^{-1}(i)}) + 1\right\rgroup\big\backslash \left\lgroup\indeg{G}{O}({\scriptstyle f^{-1}(i)}) + 1\right\rgroup $, or $\mathbf{x}^{\mathbf{b}}\not\in
\mathbf{\mathfrak{m}}^{\left(\degvec{G}+\mathbf{1}\right)\backslash\left(\indegvec{G}{O}+\mathbf{1}\right)}$. This shows that
$\mathbf{x}^{\mathbf{b}}\not\in A_{G}^{[\degvec{G}+\mathbf{1}]}$, therefore $T_G = A_{G}^{[\degvec{G}+\mathbf{1}]}$. 
\par
\qed\end{proof2}

\begin{cor}\label{cor:acycideal}
Let $G=G([n],E)$ be a simple connected graph. Then:
$$A_G = \bigcap\left\{\mathbf{\mathfrak{m}}^{\deginvec{G}{\upsigma}+\mathbf{1}_{\upsigma}}:\text{$\upsigma\in\upset{[n]}$ and $G[\upsigma]$ is connected}\right\},$$
is the irreducible decomposition of $A_G$. Also:
$$T_G = \bigcap\left\{\mathbf{\mathfrak{m}}^{\outdegvec{O}{G}+\mathbf{1}}:\text{$O$ is an acyclic orientation of $G$}\right\},$$
is the irreducible decomposition of $T_G$. 
\end{cor}

\begin{defn}\label{defn:complexlabs}
For a simple connected graph $G=G([n],E)$, consider the polytopal complexes
$\CZon{G}$, $\mathcal{Y}_G$ and $\mathcal{X}_G$, which respectively realize the abstract cell 
complexes $\mathscr{Z}_G$, $\mathscr{Y}_G$ and $\mathscr{X}_{G}^{\ast}$ of Definitions~\ref{defn:complexzon},~\ref{defn:complex1} 
and~\ref{defn:complex2}. We will let $Z_G=(\CZon{G},\ell_z)$, $Y_G=(\mathcal{Y}_G,\ell_y)$ and $X_G=(\mathcal{X}_G,\ell_x)$ 
be the $\NN^{[n]}$-labelled cell complexes with underlying polytopal complexes given by 
$\CZon{G}$, $\mathcal{Y}_G$ and $\mathcal{X}_G$, respectively, and face labelling functions $\ell_z$, $\ell_y$, $\ell_x$, defined according to:
\begin{enumerate}
\item $Z_G$: For a face $F$ of $\CZon{G}$ corresponding to $O\in\underline{\mathscr{Z}_G}$:
\begin{align*}
\ell_z(F)_i=\nooutdeg{G}{O}(i) + 1\text{, $i\in[n]$.} 
\end{align*}
\item $Y_G$: 
\begin{enumerate}
\item For a face $F$ of $\mathcal{Y}_G$ corresponding to $(\upsigma,O)\in\underline{\mathscr{X}_G}\subseteq \underline{\mathscr{Y}_G}$:
\begin{align*}
\ell_y(F)_i = \left\{\begin{array}{ll}
\nooutdeg{G[\upsigma]}{O}(i) + 1&\text{if $i\in\upsigma$,}\\
\degnovec{G}(i)+2&\text{otherwise.}
\end{array}\right.
\end{align*}
\item For a face $F$ of $\mathcal{Y}_G$ corresponding to $\upsigma \in \midset{[n]}\subseteq \underline{\mathscr{Y}_G}$:
\begin{align*}
\ell_y(F)_i=\left\{\begin{array}{ll}\degnovec{G}(i)+2&\text{if $i\in\upsigma$,}\\
0&\text{otherwise.}\end{array}\right.
\end{align*}
\end{enumerate}
\item $X_G$: For a face $F$ of $\mathcal{X}_G$ corresponding to $(\upsigma,O)\in\underline{\mathscr{X}_{G}^{\ast}}$:
\begin{align*}
\ell_x(F)_i= \left\{\begin{array}{ll}\outdeg{G[\upsigma]}{O}(i)+\degout{G}{\upsigma}(i)+1&\text{if $i\in\upsigma$,}\\
0&\text{otherwise.}\end{array}\right.
\end{align*}
\end{enumerate}
\end{defn}

\begin{lem}\label{lem:lcmz}
Let $G=G([n],E)$ be a simple connected graph. Then, for any face $F$ of $Z_G$ with vertices $v_1,\dots,v_k$, we have that:
$$\mathbf{x}^{\ell_z(F)} = \emph{LCM}\left\{\mathbf{x}^{\ell_z(v_i)}\right\}_{i\in[k]},$$
where $\emph{LCM}$ stands for ``least common multiple''. 
\end{lem}
\begin{proof2}
Let $F$ be a face of $Z_G$ with corresponding \emph{p.a.o.} $O$ of $G$ and connected partition $\Upsigma$. Every acyclic orientation of
$G$ that corresponds to a vertex of $F$ is obtained by {\bf 1)} selecting an acyclic orientation for each of the $G[\upsigma]$ with
$\upsigma\in \Upsigma$, and then by {\bf 2)} combining those $|\Upsigma|$ acyclic orientations with $O[E]\cap \overleftrightarrow{E}$. For a fixed vertex
$i\in\upsigma$ with $\upsigma\in\Upsigma$, it is possible to select an acyclic orientation of $G[\upsigma]$ in which $i$ is maximal and then to extend this to an acyclic
orientation of $G$ that refines $O$, so if vertex $v_j$ of $F$ corresponds to one such orientation, then $\ell_z(v_j)_i = \nooutdeg{G}{O}(i)+1$. 
On the other hand, clearly $\ell_z(v_j)_i\leq \nooutdeg{G}{O}(i)+1$ for all vertices $v_j$ of $F$. Hence, 
$\mathbf{x}^{\ell_z(F)} = \emph{LCM}\left\{\mathbf{x}^{\ell_z(v_i)}\right\}_{i\in[k]}$.   
\par
\qed\end{proof2}

\begin{cor}\label{cor:lcmy}
Similarly, for $G$ as in Lemma~\ref{lem:lcmz} and for any face $F$ of $Y_G$ with vertices $v_1,\dots,v_k$, we have that:
$$\mathbf{x}^{\ell_y(F)} = \emph{LCM}\left\{\mathbf{x}^{\ell_y(v_i)}\right\}_{i\in[k]},$$
where $\emph{LCM}$ stands for ``least common multiple''.  
\end{cor}
\begin{proof2}
If $F$ is a face of $Y_G$ inside the simplex $N\Updelta$, then this is immediate. If
$F$ corresponds to some $(\upsigma, O)$, then this is a consequence of the proof of Lemma~\ref{lem:lcmz}, since the vertices of $F$ are
all the $N\cdot e_i$ with $i\in[n]\backslash \upsigma$, and all the vertices of $Y_G$ that correspond to acyclic orientations of $G$ whose restrictions to
$G[\upsigma]$ refine $O$ and in which all edges of $G$ connecting $\upsigma$ with $[n]\backslash\upsigma$ are directed out of $\upsigma$.  
\par
\qed\end{proof2}

\begin{prop}\label{prop:resofY}
Let $G=G([n],E)$ be a simple connected graph. The cellular free complex $\mathcal{F}_{Y_G}$ supported on $Y_G$ is a minimal free
resolution of the artinian quotient 
$\kfield[x_1,\dots,x_n]\big\slash \left(A_G + \mathbf{\mathfrak{m}}^{\degvec{G}+2} \right)$. 
\end{prop}
\begin{proof2}
Without loss of generality, we assume here that $N>n+|E|$. 
From standard results in topological combinatorics it is easy to see that for $\mathbf{b}\in\NN^{[n]}$, the closed faces of $Y_G$ that are contained
in the closed cone $C_{\preceq\mathbf{b}}=\{\mathbf{v}\in\R^{[n]}:\mathbf{v}\preceq \mathbf{b}\}$ form a contractible polytopal complex, whenever this cone contains at least one
face of $Y_G$. Now, suppose that $\mathbf{b}$ satisfies that $b_i\leq \degnovec{G}(i)+1$ for all $i\in[n]$. Then, the complex of faces of 
$Y_G$ in the cone $C_{\preceq\mathbf{b}}$ coincides with $Y_{G,\preceq \mathbf{b}}$, so the later is contractible and acyclic if non-empty. On the contrary, let $U_{\mathbf{b}}$
be the set of all $i$ such that $b_i\geq \degnovec{G}(i)+2$, and let $D_{\mathbf{b}}=[n]\backslash U_{\mathbf{b}}$. Consider the vector
$\mathbf{a} \in\R^{[n]}$ such that:
$$a_i =  \left\{\begin{array}{ll}N&\text{if $i\in U_{\mathbf{b}}$,}\\
b_i&\text{if $i\in D_{\mathbf{b}}$.}\end{array}\right.$$   
Then, the set of faces of $Y_G$ in the cone $C_{\preceq \mathbf{a}}$ coincides with $Y_{G,\preceq\mathbf{b}}$, so again the later is contractible and acyclic if non-empty. 
This shows that $F_{Y_G}$ supports a cellular resolution of $\kfield[x_1,\dots,x_n]\big\slash \left(A_G + \mathbf{\mathfrak{m}}^{\degvec{G}+2} \right)$.
\par
To prove that this resolution is minimal, it suffices to check that whenever $F_0$ and $F_1$ are closed faces of $Y_G$ such that
$F_0\subsetneq F_1$, then $\ell_y(F_0)\prec \ell_y(F_1)$. There are three cases to study:
\begin{enumerate}
\item $F_0$ and $F_1$ correspond respectively to $\upsigma_0,\upsigma_1\in\midset{[n]}\subseteq \underline{\mathscr{Y}_G}$: 
\\*
Then, $\upsigma_0\subsetneq\upsigma_1$ and for $i\in\upsigma_1\backslash\upsigma_0$, $\ell_y(F_0)_i=0<\degnovec{G}(i)+2=\ell_y(F_1)_i$. 
\item $F_0$ corresponds to $\upsigma_0\in\midset{[n]}\subseteq \underline{\mathscr{Y}_G}$ and $F_1$ to 
$(\upsigma_1,O_1)\in\underline{\mathscr{X}_{G}^{\ast}}$:
\\*
Then, $\upsigma_0\subseteq [n]\backslash \upsigma_1$ and for $i\in\upsigma_1$, $\ell_y(F_0)_i= 0<1\leq \nooutdeg{G[\upsigma_1]}{O_1}(i) + 1 = \ell_y(F_1)_i$.
\item $F_0$ and $F_1$ correspond respectively to $(\upsigma_0,O_0),(\upsigma_1,O_1)\in\underline{\mathscr{X}_{G}^{\ast}}$:
\\*
Therefore, $\underline{J_{G[\upsigma_1]}}(O_1)\subsetneq \underline{J_{G[\upsigma_0]}}(O_0)$, so {\bf 1)} if $\upsigma_1\subsetneq\upsigma_0$, then
for $i\in\upsigma_0\backslash\upsigma_1$, $\ell_y(F_0) = \nooutdeg{G[\upsigma_0]}{O_0}(i) + 1 \leq\degnovec{G}(i)+1<\degnovec{G}(i)+2 = \ell_y(F_1)$;
and {\bf 2)} if $\upsigma = \upsigma_0 = \upsigma_1$, then letting $\Upsigma_0$ and $\Upsigma_1$ be the connected partitions of $G[\upsigma]$ corresponding respectively to
$O_0$ and $O_1$, we observe that $\Upsigma_0$ is a strict refinement of $\Upsigma_1$, so there exist $\uprho_0\in\Upsigma_0$ and $\uprho_1\in\Upsigma_1$ such that
$\uprho_0\subsetneq\uprho_1$ and such that for some $i\in\uprho_1\backslash\uprho_0$, we have that 
$\ell_y(F_0)_i=\nooutdeg{G[\upsigma]}{O_0}(i) + 1 <  \nooutdeg{G[\upsigma]}{O_1}(i)+1=\ell_y(F_1)_i$, since $G[\uprho_1]$ is connected (so there is an edge
directed out of $i$ in $O_0$ which was not directed in $O_1$).  
\end{enumerate}
\par
\qed\end{proof2}

\begin{prop}\label{prop:resofZ}
For $G$ as in Proposition~\ref{prop:resofY}, the cellular free complex $\mathcal{F}_{Z_G} = \mathcal{F}_{Y_G,\preceq\degvec{G}+\mathbf{1}}$ supported on $Z_G$ 
gives a minimal free
resolution of the quotient ring: 
$\kfield[x_1,\dots,x_n]\big\slash A_G$. 
\end{prop}
\begin{proof2}
This follows from the proof of Proposition~\ref{prop:resofY}, since $Z_G = Y_{G,\preceq \degvec{G}+\mathbf{1}}$.
\par
\qed\end{proof2}

\begin{cor}\label{cor:coresofX}
For $G$ as in Proposition~\ref{prop:resofY}, let $Y_{G}^{\emph{col}}= \degvec{G}+\mathbf{2}-Y_G$. Then, the cocellular free complex 
$\mathcal{F}^{Y_{G,\preceq\mathbf{d}_G+\mathbf{1}}^{\emph{col}}}$ supported on $Y_{G}^{\emph{col}}$ is a minimal 
cocellular resolution of the monomial ideal
$T_G$. 
\end{cor}

\begin{prop}\label{prop:resofX}
For $G$ as in Proposition~\ref{prop:resofY}, the cellular free complex 
$\mathcal{F}_{X_G}$ supported on $X_G$ is a minimal 
cellular resolution of the monomial ideal $T_G$. 
\end{prop}
\begin{proof2}
This is now a consequence of Corollary~\ref{cor:coresofX}, since the underlying polytopal complex of $Y_{G, \preceq\mathbf{d}_G+\mathbf{1}}^{\emph{col}}$ 
is combinatorially dual to the underlying complex of $X_G$, and cells from both complexes dual to each other have equal labels: If a face $F_y$ of $\mathcal{Y}_G$ and
a face $F_x$ of $\mathcal{X}_G$
both correspond to $(\upsigma,O)\in\underline{\mathscr{X}_G}$, then,
\begin{align*}
\degnovec{G}(i)+2-\ell_y(F_y)_i&=\left\{\begin{array}{ll}
\degnovec{G}(i)+2-\left(\nooutdeg{G[\upsigma]}{O}(i) + 1\right)&\text{if $i\in\upsigma$,}\\
\degnovec{G}(i)+2-\left(\degnovec{G}(i)+2\right)&\text{otherwise,}
\end{array}\right.\\
&= \left\{\begin{array}{ll}\outdeg{G[\upsigma]}{O}(i)+\degout{G}{\upsigma}(i)+1&\text{if $i\in\upsigma$,}\\
0&\text{otherwise.}\end{array}\right.\\
&=\ell_x(F_x).
\end{align*} 
\par
\qed\end{proof2}
The following is, in reality, a well-known result about Betti numbers of monomial quotients with a given cellular resolution, 
and not a definition. We present it here as a definition given
its immediate connection to the topology of cellular complexes, clearly central for the results of this section. 
\begin{defn}\label{defn:bettisturm}
If $\mathcal{F}_{X}$ is a cellular resolution of the monomial quotient
$S\slash I$, then the Betti numbers of $I$ are the numbers calculated, for all $i\geq 1$, as:
$$\upbeta_{i,\mathbf{b}}(I)=\dim_{\kfield}\widetilde{H}_{i-1}(X_{\prec\mathbf{b}};\kfield),$$
where $\widetilde{H}_{\ast}$ stands for the reduced homology functor.
\end{defn}

\begin{lem}\label{lem:facenumY}
For a simple connected graph $G=G([n],E)$, the Betti numbers of the ideals $A_G$ and
$T_G$ satisfy that, for all $i\geq 0$: 
\begin{align*}
\sum_{\mathbf{b}\in\NN^{[n]}}\upbeta_{i,\mathbf{b}}(A_G)&=\text{$\boldsymbol{\varhash}$ \emph{p.a.o's} of $G$ on $n-i$ connected parts,}\\
\sum_{\mathbf{b}\in\NN^{[n]}}\upbeta_{i,\mathbf{b}}(T_G)&=\text{$\boldsymbol{\varhash}$ of pairs $(O,\upsigma)\boldsymbol{:}O$ is a \emph{p.a.o.} of $G[\upsigma]$ on 
$i+1$ connected parts.}
\end{align*}
 \end{lem}
 \begin{proof2}
These results are clear from our choice of minimal cellular resolutions for these ideals, since $i$-th syzygies of each ideal correspond to
$i$-dimensional faces of the respective geometrical complex. 
 \par
 \qed\end{proof2}
 
 \section{Non-crossing trees.}\label{sec:noncrosstrees}
 
 In this section we investigate, for a simple graph $G=G([n],E)$, a useful and novel unifying relation between the standard monomials of $T_G$,
 the rooted spanning forests of $G$, and the maximal chains of the poset of non-crossing partitions. 
 We show that,
 arguably, the phenomenology that binds these objects together and which has been hitherto discovered in the literature, 
is largely due to the existence of a simple canonical way to represent rooted spanning forests of a graph on vertex set $[n]$ 
as non-crossing spanning trees. 
\par
An analogous extension of the theory presented here to a more general poset of non-crossing partitions associated to $G$, and
the consideration of the equally arbitrary non-nesting trees and their connection to the \emph{Catalan arrangement}, will not
be discussed here, and will be the subject of a future writing by the author. 
 
\begin{figure}[hb]
\begin{tabular}{llll}
\begin{subfigure}{.25\textwidth}
  \centering
   \caption{}
  \includegraphics[width=0.9\linewidth]{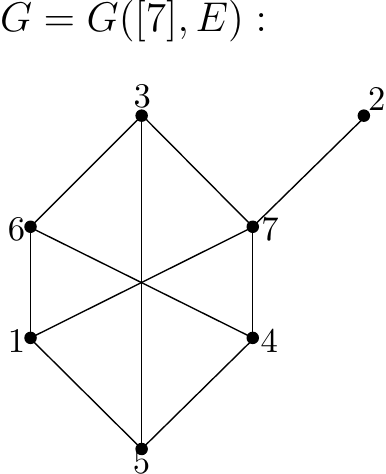}
  \label{fig:nonct1}
\end{subfigure}
&
\begin{subfigure}{.25\textwidth}
  \centering
   \caption{}
  \includegraphics[width=0.9\linewidth]{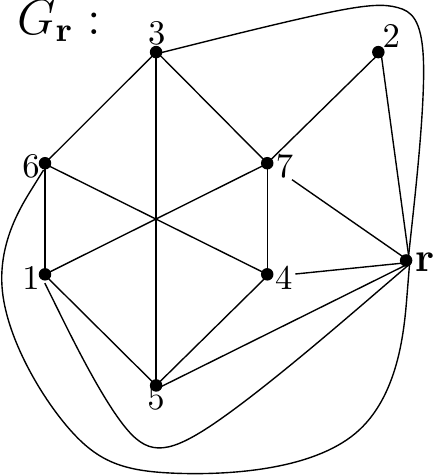}
  \label{fig:nonct2}
\end{subfigure}
&
\multicolumn{2}{c}{
\begin{subfigure}{.5\textwidth}
  \centering
   \caption{}
  \includegraphics[width=1\linewidth]{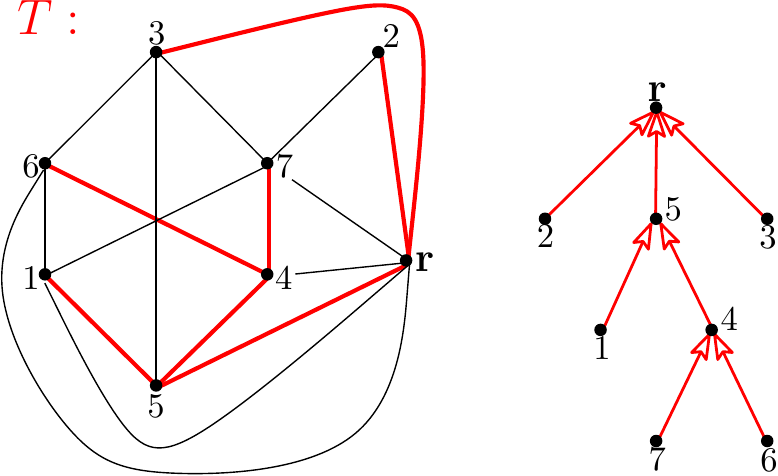}
  \label{fig:nonct3}
\end{subfigure}
}
\\
\multicolumn{3}{c}{
\begin{subfigure}{.8\textwidth}
  \centering
   \caption{}
  \includegraphics[width=1.1\linewidth]{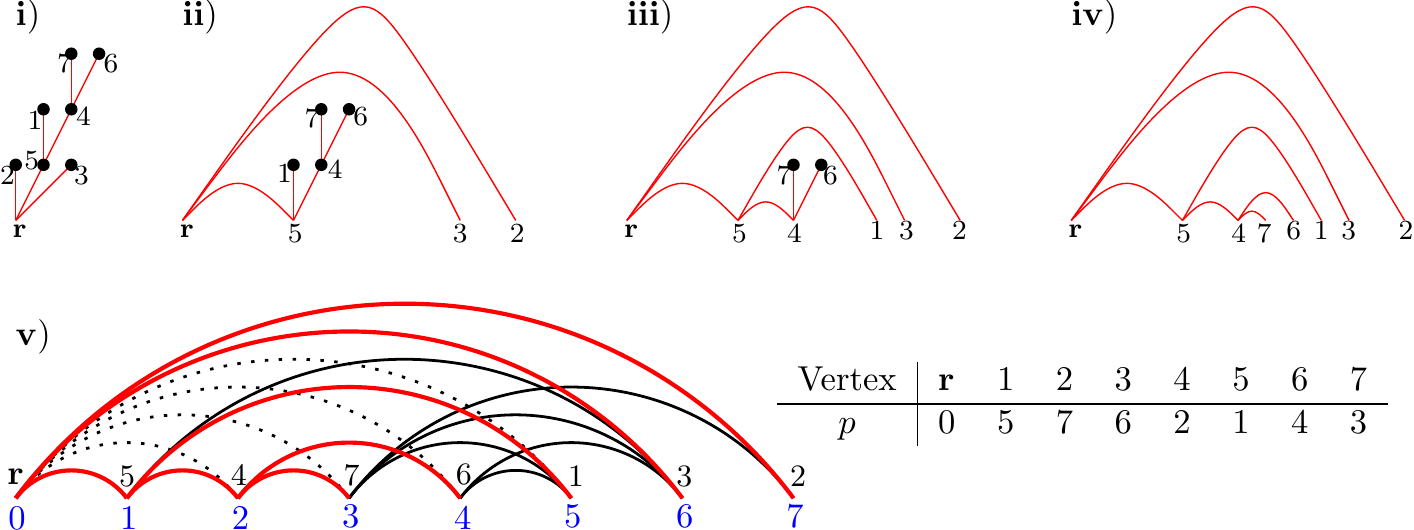}
  \label{fig:nonct4}
\end{subfigure}
}
&
\begin{subfigure}{.3\textwidth}
  \raggedleft
  \caption{}
  \includegraphics[width=0.7\linewidth]{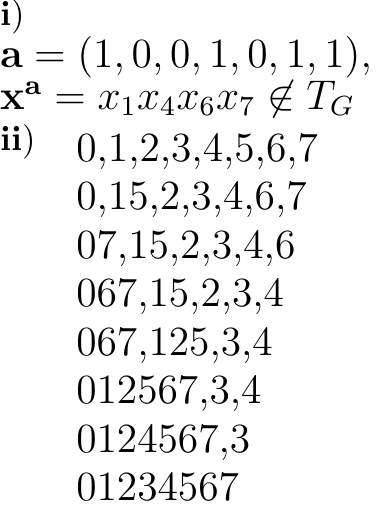}
  \label{fig:nonct5}
\end{subfigure}
\end{tabular}

\caption{Fully worked example illustrating the central dogma of Section~\ref{sec:noncrosstrees}.  
Theorems~\ref{theo:bijstdmontrees} and~\ref{theo:chainsofncp} are dwelled on in tables {\bf E.i} and 
{\bf E.ii}, respectively.}
\label{fig:nonct}
\end{figure}

 \begin{defn}\label{defn:rootgraph}
 For a simple graph $G=G(V,E)$, we will let 
 $G_{\mathbf{r}}$ denote the graph on vertex set 
 $V\sqcup\{\mathbf{r}\}$ and with edge set 
 $E\sqcup\{\{\mathbf{r},v\}:v\in V\}$, so $G_{\mathbf{r}}$ is the graph obtained from
 $G$ by adding a new vertex $\mathbf{r}$ and connecting it to all other vertices in $G$
 ({\it e.g.} Figures \ref{fig:nonct1}-\ref{fig:nonct2}).  
 \end{defn}
    
 \begin{defn}\label{defn:depiction}
A \emph{planar depiction $(D,p)$ of a finite acyclic di-graph $T=T(V,E)$} is a
finite union of closed curves $D\subseteq\R^{2}$ and a bijection $p:V\rightarrow \{0,1,2,\dots,|V|-1\}$ 
(called a \emph{depiction function}) such that: 
 \begin{enumerate}[label=\bfseries \arabic*)]
 \item $p$ is order-reversing, so if $e\in E$ and $e=(u,v)$, then $p(v)<p(u)$. 
 \item There exist strictly increasing and continuous real functions $f$ and $g$ such that
$f(0)=g(0)=0$,
 and $D$ is the image under $(f,g):\R^2\rightarrow\R^2$ of the following union of semicircles:
$$\bigcup_{(u,v)\in E}\left\{(x,y)\in\R^2:y=+\sqrt{\left(\tfrac{p(u)-p(v)}{2}\right)^2-\left(x-\tfrac{p(u)+p(v)}{2}\right)^2}\right\}.$$
 \end{enumerate}   
 A planar depiction $(D,p)$ of $T$ is said to be \emph{non-crossing} if for all $(x,y)\in D$ with $y>0$, a sufficiently small
 neighborhood of $(x,y)$ in $D$ is homeomorphic to the real line. 
 \end{defn}   
 
  \begin{lem}\label{lem:noncross}
 In Definition~\ref{defn:depiction}, the property of being a \emph{non-crossing planar depiction} is independent of the choice
 of functions $f$ and $g$, and only depends on $p$ and $T$. In other words, any two planar depictions $(D_1,p)$ and $(D_2,p)$ of $T$ are either
 both non-crossing or both crossing. 
 \end{lem}
 
\begin{figure}[ht]
\begin{tabular}{cc}
\begin{subfigure}{.3\textwidth}
  \centering
  \includegraphics[width=1.1\linewidth]{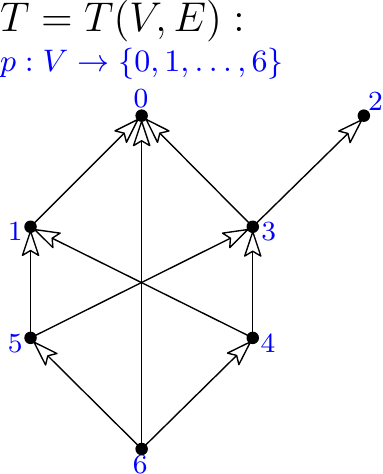}
  \caption{}
  \label{fig:depict1}
\end{subfigure}
&
\begin{subfigure}{.7\textwidth}
  \centering
  \includegraphics[width=1.1\linewidth]{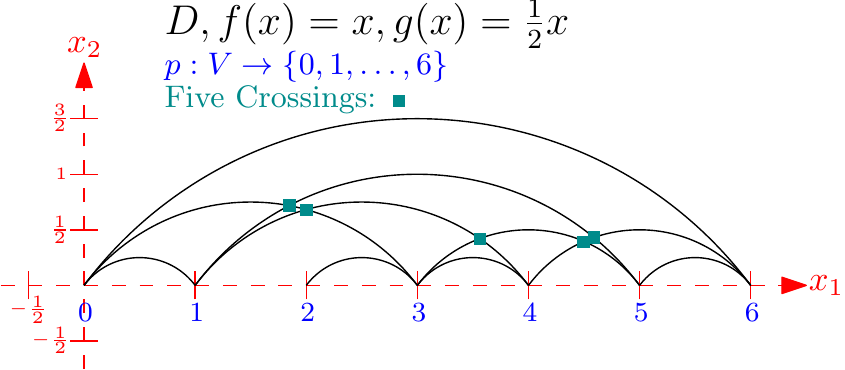}
  \caption{}
  \label{fig:depict2}
\end{subfigure}
\end{tabular}

\caption{Example of a \emph{planar depiction}, according to Definition~\ref{defn:depiction}.}
\label{fig:depict}
\end{figure}

\begin{exam}\label{exam:depiciton}
Figure~\ref{fig:depict1} shows a particular acyclic
directed graph $T=T(V,E)$ with $|V|=7$, and a choice of depiction function
$p:V\rightarrow \{0,1,2,3,4,5,6\}$ (in blue). With this choice of $p$, Figure~\ref{fig:depict2} then presents
the set $D$ obtained by taking $f(x)=x$ and $g(x)=\tfrac{1}{2}x$ in Definition~\ref{defn:depiction}. There are five \emph{crossings} in $D$, 
each marked with a square; these crossings are the points $(x,y)\in D$, $y>0$, that are locally non-homeomorphic to the real line.
\end{exam}
 
  \begin{defn}\label{defn:noncrossingtree}
A \emph{non-crossing tree} is a non-crossing planar depiction of a rooted tree $T=T(V,E)$. Vaguely, 
$T$ is obtained from an acyclic connected simple graph on vertex set $V$ by orienting all
of its edges towards a distinguished vertex of $T$, called the root of $T$ ({\it e.g.} Figure~\ref{fig:nonct3}). 
\end{defn}   
\begin{rem}\label{rem:noncrossingtree}
In Definition~\ref{defn:noncrossingtree}, for one such non-crossing tree $(D,p)$ of $T$, 
if $\mathbf{r}$ is the root of $T$, then
necessarily $p(\mathbf{r})=0$. 
\end{rem}

\begin{theo}\label{theo:noncross}
Let $G=G([n],E)$ be a simple graph, and consider a spanning tree $T$ of $G_{\mathbf{r}}$ rooted at
$\mathbf{r}$. Then, there exists a unique depiction function $p$ as in Definition~\ref{defn:depiction} such that:
\begin{enumerate}[label=\bfseries \roman*]
\item\label{noncross:c1} For all edges $(i,k)$ and $(j,k)$ of $T$, $p(i)>p(j)$ if $i<j$.
\item\label{noncross:c2} Any planar depiction $(D,p)$ of $T$ is a non-crossing tree.  
\end{enumerate}
\end{theo}
\begin{proof2}
For any two $i,j\in[n]$ with $i\neq j$, consider the directed paths from $i$ and $j$ to the root $\mathbf{r}$ of $T$. These paths meet initially at a unique vertex
$r_{ij}$ of $T$. Let us say that 
$i\prec_T j$ if either {\bf 1)} $r_{ij}=i$ or if {\bf 2)} there exist edges $(i_j,r_{ij})$ in the path from $i$ to $r_{ij}$ and $(j_i,r_{ij})$ in the path from
$j$ to $r_{ij}$ such that $i_j>j_i$. 
\par
Firstly, we verify that the relation $\preceq_T$ is a total order on the set $[n]$ of vertices of $G$. This is true since for $i\prec_T j\prec_T k$ with
$i,j,k\in[n]$:
\begin{enumerate}[label=\bfseries \alph*.]
\item If $r_{ij} = i$, then either $r_{ik}=i$ or $r_{ik}=r_{jk}$ and in the later case $i_k =j_k>k_j = k_i$.
\item If $r_{jk} = j$, then $r_{ij}=r_{ik}$ and $i_k = i_j > j_i = k_i$.
\item If $r_{ij}\neq i$, $r_{jk}\neq j$ and $r_{ij}=r_{jk}$, then $i_k = i_j > j_i = j_k > k_j = k_i$.
\item If $r_{ij}\neq i$, $r_{jk}\neq j$ and $r_{ij}\prec_T r_{jk}$, then $i_k = i_j>j_i = k_i$.
\item If $r_{ij}\neq i$, $r_{jk}\neq j$ and $r_{jk}\prec_T r_{ij}$, then $i_k = j_k > k_j = k_i$.
\end{enumerate}
Let $f:[n]\rightarrow [n]$ be the unique linear extension of this chain poset $([n],\preceq_T)$ and define
$p$ by requiring that $p(\mathbf{r}) = 0$ and $p(i)=f(i)$ for all $i\in[n]$. Clearly then $p$ satisfies 
Condition~\ref{noncross:c1}. 
\par
We now want to check that any depiction 
$(D,p)$ of $T$ is non-crossing. Suppose on the contrary that one such depiction is crossing. If that is the case, then
there exist edges $(j,i)$ and $(m,k)$ in $T$ such that $p(i)<p(k)<p(j)<p(m)$, and hence 
$j_m=j_k<k_j=m_j<j_m$, a contradiction. This proves~\ref{noncross:c2}. 
\par
To prove that $p$ is the unique bijection $[n]\cup\{\mathbf{r}\}\rightarrow\{0,1,2,\dots,n\}$
satisfying~\ref{noncross:c1}-\ref{noncross:c2}, let us suppose that another depiction function $q$ works as well. 
Since $q$ is order-reversing, then for any $i,j\in[n]$ with $i\neq j$ and $r_{ij}=i$, we must have that $q(i)<q(j)$.
If instead $r_{ij}\neq i,j$ and $i_j>j_i$, then Condition~\ref{noncross:c1} and transitivity imply that $q(i_j)<q(j_i)<q(j)$, 
and then Condition~\ref{noncross:c2} shows that $q(i_j)<q(i)<q(j_i)<q(j)$ since in any planar depiction of $T$
using $q$, the depiction of the path from $i$ to $r_{ij}$ (or to $i_j$) does not cross
the depiction of the path from $j$ to $r_{ij}$ (or to $j_i$). Hence, $q(i)<q(j)$. This shows that
$q=p$ from {\bf 1)} and {\bf 2)} above. 
\par
\qed\end{proof2}

\begin{exam}\label{exam:noncross}
Figures~\ref{fig:nonct1}-\ref{fig:nonct5} offer an example of the unique depiction function $p$ of Theorem~\ref{theo:noncross}. For the graph
$G=G([7],E)$ of Figure~\ref{fig:nonct1}, we calculate $G_{\mathbf{r}}$ in Figure~\ref{fig:nonct2}. We then select a particular
spanning tree $T$ of $G_{\mathbf{r}}$ (Figure~\ref{fig:nonct3}, in red, left diagram) and root it at $\mathbf{r}$ (Figure~\ref{fig:nonct3}, right diagram). 
Next, we present an inductive construction of the depiction function $p$ of Theorem~\ref{theo:noncross} associated to $T$.
Figure~\ref{fig:nonct4}.{\bf i-v} exhibits an inductive calculation from $T$ of a certain special diagram $D$ (in red), and the final output of this calculation is fully 
illustrated in \ref{fig:nonct4}.{\bf v}. This final diagram \ref{fig:nonct4}.{\bf v} shows a non-crossing tree from which $p$ can be instantly read off (table). 
At every step of the construction, we aim to respect
both Conditions~\ref{noncross:c1} and \ref{noncross:c2} of Theorem~\ref{theo:noncross}, and this is seen to imply
the uniqueness of $p$ for this example. In fact, it is not difficult to observe that the analogous inductive process can be 
readily applied to any other example, from
which Theorem~\ref{theo:noncross} follows. 
\end{exam}

 \subsection{Standard monomials of $T_G$.}\label{subsec:standardmon}

\begin{defn}\label{defn:bijstdmontrees}
Let $G=G([n],E)$ be a simple graph and let $\mathbf{x}^{\mathbf{a}}$ be a standard monomial
of the ideal $T_G$. From $\mathbf{a}$, let us define a bijection 
$f_{\mathbf{a}}:\{0,1,\dots,n\}\rightarrow[n]\sqcup\{\mathbf{r}\}$ and an 
$\mathbf{r}$-rooted spanning tree 
$T_{\mathbf{a}}$ of $G_{\mathbf{r}}$ recursively as follows:
\begin{enumerate}
\item The edge set $E(T_{\mathbf{a}})$ of $T_{\mathbf{a}}$ will be constructed one edge at a time. 
Similarly, a set $K$ will contain at each step the set of values in $\{0,1,\dots,n\}$
for which $f_{\mathbf{a}}$ has already been defined.
\item Initially, set $E(T_{\mathbf{a}})=\emptyset$, $f_{\mathbf{a}}(0)=\mathbf{r}$, $i=1$, and $K = \{0\}$. 
\par
Since $f_{\mathbf{a}}(k)$ has been defined for all $k\in K$, let us also denote this partially-defined function by 
$f_{\mathbf{a}}$ (which should not cause any confusion). 
\\*
\noindent\rule{5cm}{0.4pt}{\bf $_{\textstyle \emph{Step} \ i:}$}\noindent\rule{5cm}{0.4pt}

\item\label{algostart} Let $(k,j)$ be the lexicographically-maximal pair among all pairs such that: 
\begin{enumerate}[label=\bfseries \alph*)]
\item $k\in K$,
\item $j\in [n]\backslash f_{\mathbf{a}}[K]$, and 
\item\label{algostartp1} for $\{l_0<\dots <l_m\} = f^{-1}_{\mathbf{a}}\left[N_{G_{\mathbf{r}}}(j)\cap f_{\mathbf{a}}[K]\right]$, we have $k=l_{a_j}$.
 \end{enumerate}
\item From this pair $(k,j)$, set $f_{\mathbf{a}}(i)=j$ and $E(T_{\mathbf{a}})=E(T_{\mathbf{a}})\cup\{\left(j,f_{\mathbf{a}}(k)\right)\}$.
\item $K=K\cup\{i\}$. 
\item $i=i+1$.
\item Go back to~\ref{algostart} if $i\leq n$, otherwise stop. 
\end{enumerate}
\end{defn}

\begin{prop}\label{prop:bijstdmontrees1}
In Definition~\ref{defn:bijstdmontrees}, both $f_{\mathbf{a}}$ and $T_{\mathbf{a}}$ are well-defined. Furthermore, 
if we set $p_{\mathbf{a}}=f^{-1}_{\mathbf{a}}$, then $p_{\mathbf{a}}$ is the unique function of Theorem~\ref{theo:noncross} such that
any planar depiction $(D,p_{\mathbf{a}})$ of $T_{\mathbf{a}}$ is a non-crossing tree. 
\end{prop}
\begin{proof2}
If the condition of Definition~\ref{defn:bijstdmontrees}.\ref{algostart}.\ref{algostartp1} can be attained at each step of the recursion, 
that is, if for all $i\in[n]$ we are able to find
at least one such pair of $k$ and $j$ 
for which $k=l_{a_j}$, then it is clear that $f_{\mathbf{a}}$ is a bijection and 
$T_{\mathbf{a}}$ (with edge set $E(T_{\mathbf{a}})$) is a spanning $\mathbf{r}$-rooted tree
of $G_{\mathbf{r}}$. It then follows easily that $p_{\mathbf{a}}$ is order-reversing. Now suppose that we are at the $i$-th step of the recursion, $i\leq n$,
 so that $K=\{0,1,\dots,i-1\}$.  
Since for $\emptyset \neq \upsigma = [n]\backslash f_{\mathbf{a}}[K]$ we have that
$\mathbf{x}^{\degoutvec{G}{\upsigma}+\mathbf{1}_\upsigma}\in T_G$,  
then there must exist at least one $j\in\upsigma$ such that $a_j \leq \degout{G}{\upsigma}(j)$. Therefore,
if we write $\{l_0<\dots <l_m\} = f^{-1}_{\mathbf{a}}\left[N_{G_{\mathbf{r}}}(j)\cap f_{\mathbf{a}}[K]\right]$ and observe that 
in fact $m = \degout{G}{\upsigma}(j)$, it follows that $k=l_{a_j}$ is defined correctly for this choice of $j$. 
\par
Let us now establish the non-crossing condition given the choice of depiction function
$p_{\mathbf{a}}=f^{-1}_{\mathbf{a}}$. Notably, the recursive definition of $f_{\mathbf{a}}$ is tailored
at making this verification rather simple. Indeed, suppose that there exists a first step of the recursion, say the $i$-th step, $i\leq n$, where
a pair of crossing curves will be formed in any depiction $(D,p_{\mathbf{a}})$ of $T_{\mathbf{a}}$, and let $(k,j)$ be the lexicographically-maximal pair found in this step. 
Let also $(k_0,j_0)$ be the optimal pair found at the $i_0$-th step with $i_0<i$, such that the curves representing the edges
$(j_0,f_{\mathbf{a}}(k_0))$ and $(j,f_{\mathbf{a}}(k))$ cross
in all $p_{\mathbf{a}}$-depictions of $T_{\mathbf{a}}$. Then, $k_0<k<i_0<i$. This implies that the pair $(k,j)$ is lexicographically-larger
than $(k_0,j_0)$ and that, during the $i_0$-th step, the condition of Definition~\ref{defn:bijstdmontrees}.\ref{algostart}.\ref{algostartp1} is also attained
for $(k,j)$, so that $k=l_{a_j}$. Contradiction. 
\par
It remains to prove that $p_{\mathbf{a}}$ satisfies Condition~\ref{noncross:c1} of Theorem~\ref{theo:noncross}, but this follows 
immediately from the choice of lexicographically-maximal pairs at each step of the recursion.
\par
\qed\end{proof2}

\begin{defn}\label{defn:inverse}
Let $G=G([n],E)$ be a simple graph, $T$ an $\mathbf{r}$-rooted spanning tree of $G_{\mathbf{r}}$, and
$p$ the unique depiction function of Theorem~\ref{theo:noncross} associated to $T$. Let us associate
with $T$ a vector $\mathbf{b}{\scriptstyle(T)}\in \NN^{[n]}$ in the following way:
\begin{align*}
\text{For all $i\in[n]$}& \text{  and unique directed edge $(i,i_{\mathbf{r}})$ in $T$, let }\\
&\text{$b{\scriptstyle(T)}_i=\left\vert \left\{j\in N_{G_{\mathbf{r}}}(i):p(j)<p(i_{\mathbf{r}})\right\}\right\vert$.}
\end{align*}
\end{defn}

\begin{prop}\label{prop:bijstdmontrees2}
In Definition~\ref{defn:inverse}, the monomial $\mathbf{x}^{\mathbf{b}{\scriptscriptstyle(T)}}$ is a standard monomial of the ideal
$T_G$. 
\end{prop}
\begin{proof2}
Consider the bijective function $f:[n]\rightarrow[n]$ given by $f(i)=n+1-p(i)$ for all $i\in[n]$. Clearly then
$b{\scriptstyle(T)}_{f^{-1}(i)}\leq \degout{G}{f^{-1}[1,i]}({\scriptstyle f^{-1}(i)})<\degout{G}{f^{-1}[1,i]}({\scriptstyle f^{-1}(i)})+1$, and we are 
exactly in the situation
of the second part of the proof of Theorem~\ref{theo:alexander}, so we obtain that 
$\mathbf{x}^{\mathbf{b}{\scriptscriptstyle(T)}}\not\in A_{G}^{[\degvec{G}+\mathbf{1}]}=T_G$. 
\par
\qed\end{proof2}

\begin{theo}\label{theo:bijstdmontrees}
Let $G=G([n],E)$ be a simple graph, $\mathbf{x}^{\mathbf{a}}$ a standard monomial of
$T_G$, and $T$ an $\mathbf{r}$-rooted spanning tree of $G_{\mathbf{r}}$. 
Then, using the notation and functions from Definitions~\ref{defn:bijstdmontrees}-\ref{defn:inverse} and
Proposition~\ref{prop:bijstdmontrees1},
we have that $\mathbf{b}{\scriptstyle (T_{\mathbf{a}})}=\mathbf{a}$ and $T_{\mathbf{b}{\scriptscriptstyle(T)}}=T$.
Hence, the non-crossing trees obtained from the spanning trees of $G_{\mathbf{r}}$ interpolate in a bijection
between rooted spanning forests of $G$ and standard monomials of $T_G$, in such a way that every non-crossing tree naturally corresponds
to a uniquely determined object from each of these two sets of combinatorial objects associated to $G$.  
\end{theo}
\begin{proof2}
This is now a straightforward application of the recursive definition of $f_{\mathbf{a}}$ 
(or of $f_{\mathbf{b}{\scriptscriptstyle(T)}}$). For the first equality, let us 
suppose that during the $i$-th step of the recursion to define $f_{\mathbf{a}}$, so $K=\{0,1,\dots,i-1\}$ and $i\leq n$, we find a lexicographically-maximal pair
$(k,j)$ with $k=l_{a_j}$, where $\{l_0<\dots <l_m\} = f^{-1}_{\mathbf{a}}\left[N_{G_{\mathbf{r}}}(j)\cap f_{\mathbf{a}}[K]\right]$. Then: 
\begin{align*}
b{\scriptstyle(T_{\mathbf{a}})}_j &= 
\left\vert \left\{\ell\in N_{G_{\mathbf{r}}}(j):p_{\mathbf{a}}(\ell)<p_{\mathbf{a}}(j_{\mathbf{r}})\right\}\right\vert &\text{($(j,j_{\mathbf{r}})\in E(T_{\mathbf{a}})$.)}\\
&=\left\vert \left\{\ell\in N_{G_{\mathbf{r}}}(j):f^{-1}_{\mathbf{a}}(\ell)<f^{-1}_{\mathbf{a}}(f_{\mathbf{a}}{\scriptstyle(k)})\right\}\right\vert&\\
&=\left\vert \left\{\ell\in N_{G_{\mathbf{r}}}(j):f^{-1}_{\mathbf{a}}(\ell)< k \right\}\right\vert&\\
&=\left\vert \left\{\ell\in N_{G_{\mathbf{r}}}(j):f^{-1}_{\mathbf{a}}(\ell)< l_{a_j} \right\}\right\vert &\\
&=\left\vert \left\{\ell\in N_{G_{\mathbf{r}}}(j)\cap f_{\mathbf{a}}[K]:f^{-1}_{\mathbf{a}}(\ell)< l_{a_j}\leq i-1 \right\}\right\vert = a_j.&
\end{align*}
This proves the first equality. 
\par
For the second equality, we use induction on $N$ to prove that $f_{\mathbf{b}{\scriptscriptstyle(T)}}(N)=p^{-1}(N)$ for all
$N\in\{0,1,\dots,n\}$, and then to argue that during step $N\geq 1$ of the recursion to define $f_{\mathbf{b}{\scriptscriptstyle(T)}}$,
$N\leq n$, the
edge that will be added to the set $E\left(T_{\mathbf{b}{\scriptscriptstyle(T)}}\right)$ is an edge of $T$. 
Initially, when $N=0$, we have
$f_{\mathbf{b}{\scriptscriptstyle(T)}}(0)=p^{-1}(0)=\mathbf{r}$ and $E\left(T_{\mathbf{b}{\scriptscriptstyle(T)}}\right)=\emptyset$. 
Suppose that the result is true for all $N<i$, $i\in[n]$, and let us consider the
$i$-th step of the recursion, so that $K=\{0,1,\dots,i-1\}$. By induction, 
if $j\in [n]\backslash f_{\mathbf{b}{\scriptscriptstyle(T)}}[K]$ and
$\{l_0<\dots <l_m\} = f^{-1}_{\mathbf{b}{\scriptscriptstyle(T)}}\left[N_{G_{\mathbf{r}}}(j)\cap f_{\mathbf{b}{\scriptscriptstyle(T)}}[K]\right]$, 
since $f_{\mathbf{b}{\scriptscriptstyle(T)}}(k)=p^{-1}(k)$ for all $k\in K$, we have that when $\mathbf{b}{\scriptscriptstyle(T)}_j\leq m$: 
\begin{align*}
l_{\mathbf{b}{\scriptscriptstyle(T)}_j}&=l_{\left\vert \left\{\ell\in N_{G_{\mathbf{r}}}(j):p(\ell)<p(j_{\mathbf{r}})\right\}\right\vert}&\text{($(j,j_{\mathbf{r}})\in E(T)$, definition of
$\mathbf{b}{\scriptstyle(T)}$)}\\
&=p(j_{\mathbf{r}}) & \text{(definition of $l_{\ast}$ and induction)}
\end{align*}
Hence, the choice of lexicographically-maximal pair $(k,j)$ necessarily corresponds to an edge of $T$, that is, 
$(j,f_{\mathbf{b}{\scriptscriptstyle(T)}}(k))\in E(T)$. Letting $s:= p^{-1}(i)$ and $(s,s_{\mathbf{r}})\in E(T)$, 
that maximal pair selected from $T$ is
easily seen to be $\left(p(s_{\mathbf{r}}),s\right)$, again by the inductive step and the conditions
satisfied by $p$ (and $T$) from Theorem~\ref{theo:noncross}.   
\par
\qed\end{proof2}

\begin{exam}\label{exam:bijstdmontrees}
Figure~\ref{fig:nonct5}.{\bf i} presents the standard monomial of $T_G$ that corresponds to the spanning $T$ tree of
$G_{\mathbf{r}}$ in Example~\ref{exam:noncross}. For example, to calculate
$(\mathbf{a})_4=a_4$, we find cusp $4$ (in black) in Figure~\ref{fig:nonct4}.{\bf v}. To the left
of cusp $4$ in this diagram, there is exactly one adjacent cusp to $4$ through a red arc. This is cusp 
$5$ (in black), so we say that $5=4_{\mathbf{r}}$. 
There is exactly \textcolor{red}{{\bf one}} cusp in the diagram strictly to the left of $5$ that is adjacent to $4$, that is
$\mathbf{r}$. Therefore, $a_4 =$ \textcolor{red}{ $1$}, as in Definition~\ref{defn:inverse}. 
\end{exam}

\begin{prop}\label{prop:acyctrees}
Let $G=G([n],E)$ be a simple graph. Then, there exists a bijection between the following sets:
\begin{enumerate} 
\item The set of acyclic orientations
of $G$. 
\item The set of $\mathbf{r}$-rooted spanning trees $T$ of
$G_{\mathbf{r}}$ such that if $p$ is the depiction function for $T$ of Theorem~\ref{theo:noncross}, then  
for all $(i,i_{\mathbf{r}})\in E(T)$ and $j\in[n]$ with 
$p(i_{\mathbf{r}})<p(j)< p(i)$, we have that $\{i,j\}\not\in E$.
\end{enumerate} 
Moreover, if $T$ (with depiction function $p$) corresponds to an acyclic orientation $O$ of $G$ under this bijection, then the function $f:[n]\rightarrow[n]$ given by
$f(m)=n+1-p(m)$ for all $m\in[n]$ is a linear extension of $O$, and for $(i,i_{\mathbf{r}})\in E(T)$ with
$i,i_{\mathbf{r}}\in[n]$, $i_{\mathbf{r}}$ covers
$i$ in $O$.  
\end{prop}
\begin{proof2}
Let us first show that the maximal (by divisibility) standard monomials of $T_G$ are in bijection with
the acyclic orientations of $G$. Let $\mathbf{a}\in \NN^{[n]}$ be such that $\mathbf{x}^{\mathbf{a}}\not\in T_G$ but 
$\mathbf{x}^{\mathbf{a}+\mathbf{e}_i}\in T_G$ for all $i\in[n]$. From the Alexander duality of $A_G$ and $T_G$, consider
an acyclic orientation $O$ of $G$ such that $a_i\leq \outdeg{G}{O}(i)$ for all $i\in[n]$. Since $a_i+1\geq \outdeg{G}{O}(i)+1$
for all $i$, then it must be the case that $a_i = \outdeg{G}{O}(i)$, so $\mathbf{a}=\outdegvec{G}{O}$. It is well-known and not difficult
to prove that the out-degree (or in-degree) sequences uniquely determine the acyclic orientations of a simple graph, so this establishes
that the maximal standard monomials of $T_G$ are in bijection with the (out-degree sequences of the) acyclic orientations of $G$.  
\par
Now, given an $\mathbf{r}$-rooted spanning tree $T$ of $G_{\mathbf{r}}$ with depiction function $p$ as in Theorem~\ref{theo:noncross}, let us
define an orientation $O$ (not necessarily a \emph{p.a.o.}) of $G$ associated to $T$. For all $e=\{i,j\}\in E$, let:
\begin{align*}
O(e)=\left\{\begin{array}{ll}(i,j)&\text{if $p(j)\leq p(i_{\mathbf{r}})$, where $(i,i_{\mathbf{r}})\in E(T)$,}\\ e 
&\text{otherwise.}\end{array}\right.
\end{align*}
Consider the out-degree sequence $\outdegvec{G}{O}$ 
associated to the orientation $O$, {\it i.e.} $\outdeg{G}{O}(i)=\left\vert \left\{j\in[n]:(i,j)\in O[E]\right\}\right\vert$ for all 
$i\in[n]$. We then note that $b{\scriptstyle(T)}_i = \outdeg{G}{O}(i)$ for all $i$, so
$\mathbf{b}{\scriptstyle(T)}=\outdegvec{G}{O}$. However, the out-degree sequence $\outdegvec{G}{O}$ corresponds
to an acyclic orientation of $G$ if and only if $T$ satisfies that for all $(i,i_{\mathbf{r}})\in E(T)$ and $j\in[n]$ with 
$p(i_{\mathbf{r}})<p(j)< p(i)$, we have that $\{i,j\}\not\in E$, since we require that all edges
of $E$ get oriented (or get mapped to directed edges) through $O$. This proves the main statement. 
\par
That $f$ is a linear extension when $O$ is an acyclic orientation follows since then, for $(i,j)\in O[E]$, necessarily
$p(j)\leq p(i_{\mathbf{r}})<p(i)$ by the Definition of $O$ from $T$ and $p$; likewise if $(i,i_{\mathbf{r}})\in E(T)$ with
$i,i_{\mathbf{r}}\in[n]$, 
then $i_{\mathbf{r}}$ covers $i$ in $O$ since $p(i_{\mathbf{r}})\geq p(j)$ for all $(i,j)\in E(T)$ and $p$ is order-reversing. 
\par
\qed\end{proof2}

\begin{figure}[ht]
  \centering
  \includegraphics[width=0.6\linewidth]{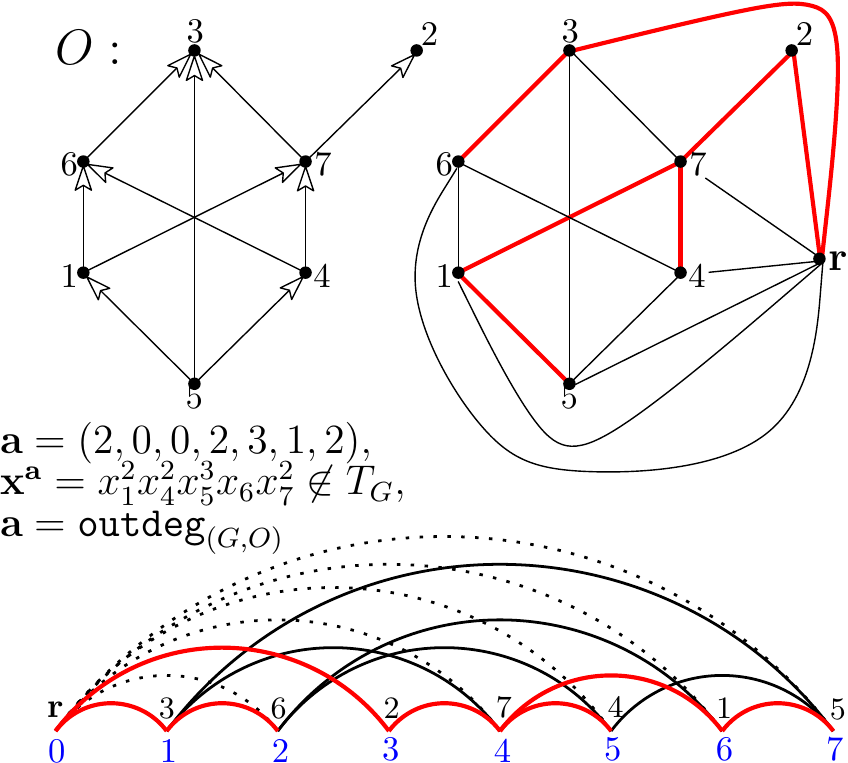}
  \caption{Example of the bijection of Proposition~\ref{prop:acyctrees}. The selected spanning tree
  of $G_{\mathbf{r}}$ (in red) corresponds to the acyclic orientation $O$ of $G$ presented.}
  \label{fig:acyct}
\end{figure}

\begin{exam}\label{exam:acyctrees}
Figure~\ref{fig:acyct} illustrates both the statement and proof of Proposition~\ref{prop:acyctrees}. Firstly,
we show an acyclic orientation $O$ of a graph $G=G([7],E)$ (Fig.~\ref{fig:acyct}, left). Then, we select a particular special 
spanning tree of $G_{\mathbf{r}}$
(Fig.~\ref{fig:acyct}, in red), and 
calculate the non-crossing tree representation of this spanning tree (Fig.~\ref{fig:acyct}, below). 
Arcs of this lower diagram represent edges of $G_{\mathbf{r}}$. To each cusp $i$ (in black) of the diagram with $i\in[7]=\{1,2,\dots,7\}$, there 
is a unique adjacent red arc to the left, and we let $i_{\mathbf{r}}$ (in black) be the other 
cusp adjacent to the same red arc, {\it e.g.} for $i=5$ we have $1=5_\mathbf{r}$. Let us orient from right to left every arc of the diagram 
adjacent to cusp $i$ if the other cusp adjacent to the arc is either $i_{\mathbf{r}}$ or lies to the left
of $i_{\mathbf{r}}$, {\it e.g.} the arcs from $5$ to $1$, $5$ to $4$, $5$ to $3$, and $5$ to $\mathbf{r}$, get all oriented 
from right to left. Doing this for all $i$,
we obtain an orientation of (some of the arcs of) the diagram, and hence an orientation of $G_{\mathbf{r}}$. In our example, this orientation 
yields an acyclic orientation of $G_{\mathbf{r}}$, and all edges are assigned an orientation; however, this might not be the case for 
several other choices of 
spanning tree of $G_{\mathbf{r}}$! Moreover, the restriction of this acyclic orientation to the edges of $G$ is precisely $O$, 
and this is
the bijection of Proposition~\ref{prop:acyctrees}. 
\end{exam}

\subsection{Non-crossing partitions.}\label{subsec:noncrosspart}

\begin{defn}\label{defn:ncp}
A \emph{non-crossing partition of the totally ordered set $[0,n]=\{0,1,\dots,n\}$} is a set partition $\uppi$ of $[0,n]$ in which
every block is non-empty and such that there does not exist integers $i<j<k<l$ and
blocks $B\neq B'$ of $\uppi$ with $i,k\in B$ and $j,l\in B'$. 
\par
The set of all non-crossing partitions of $[0,n]$ ordered by refinement ($\preceq_{\scriptscriptstyle\emph{ref}}$) forms a graded lattice
of length $n$, and we will
denote this \emph{lattice of non-crossing partitions of $[0,n]$} by $\noncross{[0,n]}$.
\end{defn}

\begin{defn}\label{defn:chainsofncp}
Consider a maximal chain $C=\{\uppi_0\precref \uppi_1\precref
\dots\precref \uppi_n\}$ of
$\noncross{[0,n]}$. For each $i\in[n]$, there exists a unique element $\bar{i}\in[n]$ such that 
$\bar{i}$ is the minimal element of its block in $\uppi_{i-1}$ but $\bar{i}$ is not the minimal element of its block in
$\uppi_{i}$. Let then $\underline{i}\neq \bar{i}$ be the element of the block of $\bar{i}$ in $\uppi_i$ that immediately
precedes $\bar{i}$. 
\par
With this notation, we define a bijection $p_{\scriptscriptstyle C}:[n]\sqcup\{\mathbf{r}\}\rightarrow [0,n]$ and an $\mathbf{r}$-rooted tree
$T_{\scriptscriptstyle C}$ of $\left(K_{[n]}\right)_{\mathbf{r}}\cong K_{[n]\sqcup\{\mathbf{r}\}}$ associated to the chain $C$ in the following way:
\begin{align*}
p_{\scriptscriptstyle C}(\mathbf{r}) &= 0,\tag{$p_{\scriptscriptstyle C}$:}\\
p_{\scriptscriptstyle C}(i) &= \bar{i}\text{, for all $i\in[n]$.}\\
E(T_{\scriptscriptstyle C}) &= \left\{\left(p_{\scriptscriptstyle C}^{-1}({\scriptstyle\bar{i}}),p_{\scriptscriptstyle C}^{-1}({\scriptstyle\underline{i}})\right): 
i\in[n]\right\}.\tag{$T_{\scriptscriptstyle C}:$}
\end{align*}
\end{defn}

\begin{prop}\label{prop:chainsofncp1}
In Definition~\ref{defn:chainsofncp}, both $p_{\scriptscriptstyle C}$ and $T_{\scriptscriptstyle C}$ are well-defined and moreover,
$p_{\scriptscriptstyle C}$ is the function of Theorem~\ref{theo:noncross} such that any planar depiction $(D,p_{\scriptscriptstyle C})$
of $T_{\scriptscriptstyle C}$ is a non-crossing tree. 
\end{prop}
\begin{proof2}
That $p_{\scriptscriptstyle C}$ is well-defined is a consequence of the fact that taking the union of two disjoint blocks in a partition of
$[0,n]$ will make exactly one minimal element of these blocks non-minimal in the newly formed block. Hence, in a maximal chain
of $\noncross{[0,n]}$, every non-zero element of $[0,n]$ stops being minimal in its own block at exactly one cover relation in the chain, and every cover relation
in the chain gives rise to one such element. That $T_{\scriptscriptstyle C}$ is an $\mathbf{r}$-rooted spanning tree of
$\left(K_{[n]}\right)_{\mathbf{r}}$ comes from observing that, since $p_{\scriptscriptstyle C}$ is well-defined, the di-graph 
$p_{\scriptscriptstyle C}\circ T_{\scriptscriptstyle C}$ on vertex set 
$[0,n]$ and edge set $(\bar{i},\underline{i})$ for all $i\in[n]$, is a $0$-rooted spanning tree of $K_{[0,n]}$. This is true because
for every $i\in[n]$, there exists exactly one edge in $p_{\scriptscriptstyle C}\circ T_{\scriptscriptstyle C}$ of the form $(i,j)$ with $j<i$, and these are all the edges of
$p_{\scriptscriptstyle C}\circ T_{\scriptscriptstyle C}$.  
\par
To verify that $p_{\scriptscriptstyle C}$ and $T_{\scriptscriptstyle C}$ satisfy Condition~\ref{noncross:c1} of Theorem~\ref{theo:noncross}, suppose on the contrary that there are
edges $(i,k),(j,k)\in E(T_{\scriptscriptstyle C})$ with $i<j$ and $p_{\scriptscriptstyle C}(i)<p_{\scriptscriptstyle C}(j)$. This means that $\bar{i}$ was minimal in its block in $\uppi_{i-1}$ but not in
$\uppi_i$, and that both $\bar{i}$ and $\underline{i}=p_{\scriptscriptstyle C}(k)$ lied in the same block of $\uppi_i$. Similarly, $\bar{j}$ was minimal in $\uppi_{j-1}$
but not in $\uppi_{j}$, where it was immediately preceded by $\underline{j}=p_{\scriptscriptstyle C}(k)=\underline{i}$. Since $j>i$, all three $\bar{i}$, $\bar{j}$ and
$\underline{j}$ belonged
to the same block of $\uppi_j$, but $\underline{j}=\underline{i}<\bar{i}=p_{\scriptscriptstyle C}(i)<p_{\scriptscriptstyle C}(j)=\bar{j}$ shows that
$\underline{j}$ does not immediately precede $\bar{j}$ in $\uppi_j$, contradiction.  
\par
To verify the non-crossing condition, note that if there is a crossing in a depiction $(D,p_{\scriptscriptstyle C})$ of $T_{\scriptscriptstyle C}$, then there is a smallest $i\in[n]$ 
such that there exists $j<i$ with either $\underline{j}<\underline{i}<\bar{j}<\bar{i}$ or $\underline{i}<\underline{j}<\bar{i}<\bar{j}$. In both cases,
we observe that $\{\underline{i},\bar{i}\}$ and $\{\underline{j},\bar{j}\}$ belong to different blocks of $\uppi_i$. 
But then, these two blocks must cross in
$\uppi_i$, clearly. This is a contradiction. 
\par
\qed\end{proof2}

\begin{defn}\label{defn:chainsofncp2}
Let $G=G([n],E)$ be a simple graph, and let $T$ be an $\mathbf{r}$-rooted spanning tree of
$G_{\mathbf{r}}$. Suppose that $p$ is the depiction function of Theorem~\ref{theo:noncross} such that
any depiction $(D,p)$ of $T$ is non-crossing. 
\par
From $T$ and $p$, let us form a chain 
$C_T=\{\uppi_0\precref \uppi_1\precref\dots\precref \uppi_n\}$ of partitions of the set $[0,n]$ in the following way:
\begin{enumerate}
\item Let $\uppi_0 = \left\{\{0\},\{1\},\dots,\{n-1\},\{n\}\right\}$, and
\item for each $i\in[n]$, let $\uppi_i$ be obtained from $\uppi_{i-1}$ by taking the union of the block that contains
$p(i)$ and the block that contains $p(i_{\mathbf{r}})$, where $(i,i_{\mathbf{r}})$ is an edge of $T$. 
\end{enumerate}
\end{defn}

\begin{prop}\label{prop:chainsofncp2}
In Definition~\ref{defn:chainsofncp2}, $C_T$ is well-defined and moreover, it is a maximal chain of partitions in
$\noncross{[0,n]}$.  
\end{prop}
\begin{proof2}
That $C_T$ is a well-defined (maximal) chain of partitions of $[0.n]$ is a consequence of $p$ being a bijection
$[n]\cup\mathbf{r}\rightarrow [0,n]$ and of $T$ being a spanning tree of $G_{\mathbf{r}}$: We can think of the procedure of 
Definition~\ref{defn:chainsofncp2} as that of beginning with an independent set of vertices $[n]\cup\mathbf{r}$,
and then adding one edge of $T$ at a time until we form $T$, keeping track at each step of the connected components of the graph
so far formed (and mapping those connected components through $p$); there are $n$ such steps and at each step we add a different
edge of $T$. In fact, since $T$ is rooted and $p$ is order-reversing, if for some $i\in[n]$ we consider the edges 
$(1,1_{\mathbf{r}}),\dots,(i,i_{\mathbf{r}})$ of $T$ that have been added up to the $i$-th step in this process (so that the graph
in consideration is a rooted forest), we see that if
two numbers $k<l$ in the set $[0,n]$ belong to the same block $B$ of $\uppi_i$, then either $(p^{-1}(l'),p^{-1}(k))$ is an edge of
$T$ for some $l'\in B$ with $k<l'\leq l$ and $p^{-1}(l')\leq i$ , 
or there exist $k',l'\in B$ with $k'<k<l'\leq l$ such that $(p^{-1}(l'),p^{-1}(k'))$ is an edge of $T$ and
$p^{-1}(l')\leq i$. 
\par
Suppose now that some of the partitions in $C_T$ are crossing, and let us assume that $i$ is minimal such that
$\uppi_i$ is crossing. Hence, the block $B_i$ in $\uppi_i$ that contains both $p(i)$ and $p(i_{\mathbf{r}})$ crosses with 
another block $B_j$ of $\uppi_i$, so there exist two consecutive elements $i_1<i_2$ of $B_i$ and two
consecutive elements $j_1<j_2$ of $B_j$ such that either {\bf a)} $i_1<j_1<i_2<j_2$ or {\bf b)} $j_1<i_1<j_2<i_2$.
In $\uppi_{i-1}$, $i_1$ and $i_2$ belong to different blocks $B_{i_1}$ and $B_{i_2}$ respectively, 
and $B_i = B_{i_1}\sqcup B_{i_2}$. Moreover, since $i$ was chosen minimally, 
if {\bf a)} holds above then $B_{i_2}\subseteq (j_1,j_2)$ and $B_{i_1}\cap (j_1,j_2)=\emptyset$, and if
{\bf b)} holds then $B_{i_1}\subseteq (j_1,j_2)$ and $B_{i_2}\cap (j_1,j_2)=\emptyset$. As
$p$ is order-reversing, so $p(i_{\mathbf{r}})<p(i)$, we see that $p(i_{\mathbf{r}})\in B_{i_1}$ and $p(i)\in B_{i_2}$,
and then that $i_1<p(i)$ and $p(i_{\mathbf{r}})<i_2$. These last two inequalities imply that
$p(i_{\mathbf{r}})\leq i_1<i_2\leq p(i)$. 
Also, since $p$ satisfies Condition~\ref{noncross:c1} of Theorem~\ref{theo:noncross}, we observe that
necessarily $i_1=p(i_{\mathbf{r}})$. Otherwise, as both $i_1$ and $p(i_{\mathbf{r}})$ belong to the same block $B_{i_1}$ of $\uppi_{i-1}$ and
$p(i_{\mathbf{r}})\leq i_1$, then
either $(p^{-1}(l),i_{\mathbf{r}})$ is an edge of $T$ for some
$l\in B_{i_1}$ with $p(i_{\mathbf{r}})<l\leq i_1$  and $p^{-1}(l)<i$ (which cannot hold since $i_1<p(i)$), or there exist $k,l\in B_{i_1}$ with
$k<p(i_{\mathbf{r}})<l\leq i_1<p(i)$ such that $(p^{-1}(l),p^{-1}(k))$ is an edge of $T$ (which cannot hold because
that edge crosses $(i,i_{\mathbf{r}})$ in any depiction $(D,p)$ of $T$). 
More easily, since $i_2\leq p(i)$ and there are no edges of the form $(i,l)$ in $T$ except for
$(i,i_{\mathbf{r}})$, we must in fact have that $i_2=p(i)$. 
It is now clear that if {\bf a)} or {\bf b)} holds above with $i_1=p(i_{\mathbf{r}})$ and
$i_2=p(i)$, then in any depiction $(D,p)$ of $T$ we may find an edge of $T$ that crosses
$(i,i_{\mathbf{r}})$, which is impossible. 
\par
\qed\end{proof2}

\begin{theo}\label{theo:chainsofncp}
Let $K_{[n]}$ be the complete graph on $[n]$, $T$ be an $\mathbf{r}$-rooted spanning tree of
$\left(K_{[n]}\right)_{\mathbf{r}}\cong K_{[n]\sqcup\{\mathbf{r}\}}$, and 
$C=\{\uppi_0\precref \uppi_1\precref\dots\precref \uppi_n\}$ a maximal chain of 
$\noncross{[0,n]}$. Then, using the notation and functions of Definitions~\ref{defn:chainsofncp}-\ref{defn:chainsofncp2},
we have that $T_{\left(C_{\scriptscriptstyle T}\right)}=T$ and $C_{\left(T_C\right)}=C$. Hence, the non-crossing trees obtained from the spanning trees of 
$\left(K_{[n]}\right)_{\mathbf{r}}$ interpolate 
in a bijection between rooted spanning forests of $K_{[n]}$ and maximal chains of the non-crossing partitions
lattice $\noncross{[0,n]}$: Every non-crossing tree corresponds bijectively to an element of each of these two combinatorial sets. 
\end{theo}
\begin{proof2}
 This is clear from the proofs of Propositions~\ref{prop:chainsofncp1}-\ref{prop:chainsofncp2} through the following simple
observations.
\par
Firstly, the edges of $T_{\scriptscriptstyle C}$ correspond to the cover relations in $C$ so that 
 an edge $(i,i_{\mathbf{r}})$ with $i\in[n]$ exists in $T_{\scriptscriptstyle C}$ for every minimal element $p_{\scriptscriptstyle C}(i)$ in its block of $\uppi_{i-1}$ that stops being minimal
 in its block of $\uppi_{i}$; the number $i_{\mathbf{r}}$ is then recollected by requiring that $p_{\scriptscriptstyle C}(i_{\mathbf{r}})$ is the immediate
 predecessor of $p_{\scriptscriptstyle C}(i)$ in the newly formed block of $\uppi_i$. Nextly, for all $i\in[n]$, the $i$-th cover relation in $C_{\left(T_C\right)}$ corresponds to taking the union
 of the block that contains $p_{\scriptscriptstyle C}(i)$ and $p_{\scriptscriptstyle C}(i_{\mathbf{r}})$. Therefore, $C = C_{\left(T_C\right)}$.
 \par
 Secondly, the $i$-th cover relation in $C_T$, $i\in[n]$, corresponds to taking the union of the (disjoint) blocks that contain $p(i)$ and $p(i_{\mathbf{r}})$, where 
 $(i,i_{\mathbf{r}})$ is an edge of $T$ (and from the second part of the proof of Proposition~\ref{prop:chainsofncp2}, 
 $p(i)$ was minimal in its initial block and $p(i_{\mathbf{r}})$ immediately precedes $p(i)$ in the newly formed block). But then, the edges of $T_{\left(C_T\right)}$ are given 
by all the $(i,i_{\mathbf{r}})$. Hence, $T_{\left(C_T\right)}=T$.   
\par
\qed\end{proof2}

\begin{exam}\label{exam:chainsofncp}
Table~\ref{fig:nonct5}.{\bf ii} shows an example of the bijection of Theorem~\ref{theo:chainsofncp}, presenting the
maximal chain of $\noncross{[0,7]}$ corresponding to the spanning tree $T$ of $G_{\mathbf{r}}$ of Example~\ref{exam:noncross} 
(top to bottom of table, blocks separated by commas). Let us discuss how this list of non-crossing partitions can be calculated from
 Figure~\ref{fig:nonct4}.{\bf v}. We will  inductively define a set of graphs $G_0,G_1,\dots,G_7$, each on vertex set $[0,7]=\{0,1,\dots,7\}$ and with
 edge sets $E_0,E_1,\dots,E_7$, respectively.
 Initially, $G_0$ has no edges, so $E_0=\emptyset$. Suppose then that we have
 defined $G_{i-1}$ and $E_{i-1}$ with $i\leq 7$, and that we want to define  
 $G_i$ and $E_i$. We find cusp 
$i$ (in black) in Figure~\ref{fig:nonct4}.{\bf v} and note that, to this cusp, there is exactly one red arc adjacent to the left. This arc is also adjacent to cusp
$i_{\mathbf{r}}$ (in black). Let us then read off the blue labellings of cusps $i$ and $i_{\mathbf{r}}$ in Figure~\ref{fig:nonct4}.{\bf v}, and say that these are
$p(i)$ and $p(i_{\mathbf{r}})$. Then, writing $e:=\{p(i),p(i_{\mathbf{r}})\}$, we let $E_i=E_{i-1}\cup\{e\}$ and update $G_i$ accordingly. We stop when $G_7$ is defined. Notably, 
$G_7$ is 
a spanning tree. Non-crossing partitions of 
Table~\ref{fig:nonct5}.{\bf ii} are then, in order, given by the connected components of the spanning forests $G_0,G_1,\dots,G_7$.  
\end{exam}

\begin{cor}[Germain Kreweras]\label{cor:numchainsncp}
The number of maximal chains in $\noncross{[0,n]}$ is $(n+1)^{n-1}$. 
\end{cor}

\begin{cor}\label{cor:decomptrees}
We have that:
$$(n+1)^{n-1}=\sum_{\{B_1,\dots,B_m\}\in\noncross{[n]}}\binom{n}{|B_1|!,|B_2|!,\dots,|B_m|!}.$$
Therefore, using Speicher's exponential formula for $\noncross{[n]}$~{\bfseries [}\cite{speicher}{\bfseries ]}, we obtain the classic result:
$$\sum_{n=1}^{\infty}n^{n-1}\frac{x^{n}}{n!}=\left(\frac{x}{e^x}\right)^{\langle -1\rangle}.$$
\end{cor}
\begin{proof2}
For each $\{B_1,\dots,B_m\}\in\noncross{[n]}$, where $b_1<b_2<\dots<b_m$ are respectively the minimal elements of $B_1$, $B_2$, $\dots$, $B_m$, 
and for each bijection $f:[n]\rightarrow[n]$ such that $f$
is strictly decreasing on each block $B_i$, $i\in[m]$, we can define an $\mathbf{r}$-rooted spanning 
tree $T$ of $\left(K_{[n]}\right)_{\mathbf{r}}$ by taking $E(T) = \left\{(f(i),\mathbf{r}):i\in B_1\right\}\cup\left\{(f(i),f(b_k-1)):i\in B_k\text{ with }k>1\right\}$. 
If we let $p(\mathbf{r}) = 0$ and $p(i)=f^{-1}(i)$ for all $i\in[n]$ above, we can readily check that $p$ is the depiction
function of Theorem~\ref{theo:noncross} associated to $T$. 
\par
Conversely, given an $\mathbf{r}$-rooted spanning tree $T$ with depiction function $p$ as in Theorem~\ref{theo:noncross}, the partition
$\displaystyle\sqcup_{k\in [n]}\left\{p(i)\in[n]:(i,k)\in T\right\}$ is an element of
$\noncross{[n]}$. 
\par
Hence, since given a partition $\{B_1,\dots,B_m\}\in\noncross{[n]}$, there are
${\scriptstyle \binom{n}{|B_1|!,|B_2|!,\dots,|B_m|!}}$ choices for $f$ above, the result follows. 
\par
\qed\end{proof2}

 \section{Applications.}\label{sec:applications}
 
 \subsection{Random Acyclic Orientations of a Simple Graph: Markov Chains.}\label{subsec:randwalk}
 
 \begin{defn}\label{defn:randwalk}
 Let $G=G(V,E)$ be a connected (finite) simple graph. A \emph{simple random walk on $G$} is a Markov chain $(v_t)_{t=0,1,2,\dots}$, obtained by selecting
an initial vertex $v_0\in V$, and then for all $t\geq 1$, selecting $v_t\in V$ from a uniform distribution
 on the set $N_G(v_{t-1})$. If $P$ is the Markov transition matrix for a simple random walk on $G$, then for 
 $u,v\in V$:
 \begin{align*}
 (P)_{uv}=p_{uv}=\left\{\begin{array}{ll}\tfrac{1}{\degnovec{G}(u)}&\text{if $v\in N_G(u)$,}\\
 0&\text{otherwise.}\end{array}\right.
 \end{align*}
 \end{defn}
 
 \begin{theo}\label{theo:randwalk}
 The Markov chain of Definition~\ref{defn:randwalk} is always \emph{irreducible}. Furthermore, it is
 \emph{aperiodic} if and only if $G$ is not bipartite. 
 \par
If for all $w\in V$, we let $\uppi_w := \tfrac{\degnovec{G}(w)}{2|E|}$, 
then for any pair of vertices $u,v\in V$, we have that: 
 $$\uppi_v p_{vu} = \uppi_u p_{uv}.$$
Consequently, since $\sum_{v\in V}\uppi_v =1$, random walks on $G$ are \emph{reversible} and 
they have a unique stationary distribution given by $\boldsymbol{\uppi}=(\uppi_v)_{v\in V}$, so that:
\begin{equation}\label{eq:statdist}
\uppi_v = \lim_{N\rightarrow \infty}\tfrac{1}{N}\sum_{t=1}^{N}\proba{v_t = v}\text{, for all $v\in V$}.
\end{equation}  
Moreover, if $G$ is not bipartite, then:
\begin{equation}\label{eq:statdistaper}
\uppi_v = \lim_{t\rightarrow \infty}\proba{v_t = v}\text{, for all $v\in V$}.
\end{equation}
 \end{theo}
 
 \begin{defn}[Card-Shuffling Markov Chain, see also~\cite{athanasiadisdiaconis}]\label{defn:cardshuffle}
 Let $G=G(V,E)$ be a simple graph with $|V|=n\geq 3$, and select an arbitrary bijective 
 map $f_0:V\rightarrow [n]$ (regarded as a labelling of $V$). Let us consider a sequence $(f_t)_{t=0,1,2,\dots}$ of
 bijective maps $V\rightarrow [n]$ such that for $t\geq 1$,
 $f_t$ is obtained from $f_{t-1}$ through the following random process:
 Let $v_t\in V$ be chosen uniformly at random, and let, 
$$f_t(v)=\left\{\begin{array}{ll} n&\text{if $v = v_t$,}\\
f_{t-1}(v)-1&\text{if $f_{t-1}(v)>f_{t-1}(v_t)$,}\\
f_{t-1}(v)&\text{otherwise.}\end{array}\right.$$
 Consider now the sequence of acyclic orientations $(O_t)_{t=0,1,2,\dots}$ of $G$ induced by the labellings $(f_t)_{t=0,1,2,\dots}$, so that
 for all $e=\{u,v\}\in E$ and $t\geq 0$, we have that $O_t(e)=(u,v)$ if and only if $f_t(u)<f_t(v)$. The sequence
 $(O_t)_{t=0,1,2,\dots}$ is called the \emph{Card-shuffling (CS) Markov chain on the set of acyclic orientations of $G$}.  
 \par
 Equivalently, we can define this Markov chain by selecting an arbitrary acyclic orientation 
 $O_0$ of $G$, and then for each $t\geq 1$, letting $O_t$ be obtained from $O_{t-1}$ by selecting
 $v_t\in V$ uniformly at random and taking, for all $e\in E$:
 $$O_t(e)=\left\{\begin{array}{ll}O_{t-1}(e)&\text{if $v_t\not\in e$,}\\
 (v,v_t)&\text{if $e=\{v,v_t\}$.} \end{array}\right.$$ 
 \end{defn}

\begin{theo}\label{theo:cardshuffle}
The Card-Shuffling Markov chain of $G$ in Definition~\ref{defn:cardshuffle} is a well-defined, irreducible and aperiodic Markov chain
on state space equal to the set of all acyclic orientations of $G$; its unique 
stationary distribution $\boldsymbol{\uppi^{\scriptscriptstyle \emph{CS}}}$ is given by:
$$\uppi^{\scriptscriptstyle \emph{CS}}_O = \tfrac{e(O)}{n!}, \text{ for all acyclic orientations $O$ of $G$,}$$
where $e(\cdot)$ denotes the number of linear extensions of the induced poset $(V,\leq_O)$. 
\end{theo}
\begin{proof2}
If we consider instead the Markov chain $(f_t)_{t=0,1,2,\dots}$, whose set of states is the set of all bijections
$V\rightarrow[n]$, it is not difficult to observe that this Markov chain is irreducible and aperiodic (see below),
and hence that it has a unique stationary distribution $\boldsymbol{\uppi}$ satisfying
Equations~\ref{eq:statdistaper}. By the symmetry of the set of all bijective labelings
$V\rightarrow[n]$, or simply by direct inspection of the stationary equations for this Markov chain 
(since every state can be accessed in one step from exactly
$n$ different states and each one of these transitions occurs with probability $\tfrac{1}{n}$), we obtain that $\uppi_f = \tfrac{1}{n!}$ for all 
bijective maps $f:V\rightarrow[n]$.  
Hence, by the aforementioned construction of the Card-Shuffling (CS) Markov chain of $G$ from bijective labellings of $V$, we must have
that this CS chain is also irreducible (since each labelling is accessible from every other labelling, hence each acyclic orientation
from every other acyclic orientation), aperiodic 
(since both $\proba{f_t = f_{t-1}\vert f_{t-1}}>0$ and $\proba{O_t = O_{t-1}\vert O_{t-1}}>0$ for all $t\geq 1$), 
and has a unique stationary distribution $\boldsymbol{\uppi^{\scriptscriptstyle \emph{CS}}}$, necessarily then given 
by $\uppi^{\scriptscriptstyle \emph{CS}}_O = \tfrac{e(O)}{n!}$ for every acyclic orientation $O$ of $G$, from
Equations~\ref{eq:statdist}.   
\par
\qed\end{proof2}

 \begin{defn}[Edge-Label-Reversal Stochastic Process]\label{defn:edgelabrev}
 Let $G=G(V,E)$ be a connected simple graph with $|V|=n$, and select an arbitrary bijective 
 map $f_0:V\rightarrow [n]$ (regarded as a labelling of $V$). Let us consider a sequence $(f_t)_{t=0,1,2,\dots}$ of
 bijective maps $V\rightarrow [n]$ such that for $t\geq 1$,
 $f_t$ is obtained from $f_{t-1}$ through the following random process:
 Let $e_t=\{u_t,v_t\}\in E$ be chosen uniformly at random from this set, and let, 
$$f_t(v)=\left\{\begin{array}{ll}f_{t-1}(u_t) &\text{if $v = v_t$,}\\
f_{t-1}(v_t)&\text{if $v=u_t$,}\\
f_{t-1}(v)&\text{otherwise.}\end{array}\right.$$
 Consider now the sequence of acyclic orientations $(O_t)_{t=0,1,2,\dots}$ of $G$ induced by the labellings $(f_t)_{t=0,1,2,\dots}$, so that
 for all $e=\{u,v\}\in E$ and $t\geq 0$, we have that $O_t(e)=(u,v)$ if and only if $f_t(u)<f_t(v)$. The sequence
 $(O_t)_{t=0,1,2,\dots}$ is called the \emph{Edge-Label-Reversal (ELR) stochastic process on the set of acyclic orientations of $G$}.  
 \end{defn}
 
 \begin{theo}\label{theo:edgelabrev}
The Edge-Label-Reversal stochastic process of $G$ in Definition~\ref{defn:edgelabrev} satisfies that,
for every acyclic orientation $O$ of $G$:
\begin{align*}
\left(\boldsymbol{\uppi^{\scriptscriptstyle \emph{ELR}}}\right)_O=\uppi^{\scriptscriptstyle \emph{ELR}}_O &:= \lim_{N\rightarrow \infty}\tfrac{1}{N}\sum_{t=1}^{N}\proba{O_t=O}=\tfrac{e(O)}{n!},
\end{align*}
where $e(O)$ denotes the number of linear extensions of the induced poset $(V,\leq_O)$, and this
result holds independently of the initial choice of $O_0$.
\end{theo}
\begin{proof2}
Consider the simple graph $H$ on vertex set equal to the set of all bijective maps $V\rightarrow[n]$, and where two
maps $f$ and $g$ are connected by an edge if and only if there exists $\{u,v\}\in E$ such that $f(u)=g(v)$,
$f(v)=g(u)$, and $f(w)=g(w)$ for all $w\in V\backslash\{u,v\}$. Since $G$ is connected, a standard result
in the algebraic theory of the symmetric group shows that $H$ is connected, {\it e.g.} consider a spanning tree $T$ of $G$; then,
any permutation in $\mathfrak{S}_{V}$ can be written as a product of transpositions of the form $(u \ v)$ with $\{u,v\}\in E(T)$. Moreover,
by considering the parity of permutations in $\mathfrak{S}_{V}$, we observe that $H$ is bipartite. Now, the sequence 
$(f_t)_{t=0,1,2,\dots}$ of Definition~\ref{defn:edgelabrev} is precisely a simple random walk on $H$, and the degree of each bijective map $f:V\rightarrow[n]$ in $H$
is clearly $|E|$, so the stationary distribution for this Markov chain is uniform. Necessarily then, the result follows from the
construction of $(O_t)_{t=0,1,2,\dots}$ and Equations~\ref{eq:statdist}.  
\par
\qed\end{proof2}
 
  \begin{defn}[Sliding-$(n+1)$ Stochastic Process]\label{defn:sliding}
 Let $G=G(V,E)$ be a connected simple graph with $|V|=n$, and consider the graph
 $G_{\mathbf{r}}$. Let us select an arbitrary bijective 
 map $f_0:V\sqcup\{\mathbf{r}\} \rightarrow [n+1]$, which we regard as a labelling of the vertices of $G_{\mathbf{r}}$, and
 define a sequence $(f_t)_{t=0,1,2,\dots}$ of
 bijective maps $V\sqcup\{\mathbf{r}\}\rightarrow [n+1]$ such that for $t\geq 1$,
 $f_t$ is obtained from $f_{t-1}$ through the following random process:
 Let $v_{t-1}\in V\sqcup\{\mathbf{r}\}$ be such that $f_{t-1}(v_{t-1})=n+1$, and select $v_t\in N_{G_{\mathbf{r}}}(v_{t-1})$ uniformly at random from this set. Then, 
$$f_t(v)=\left\{\begin{array}{ll}n+1 &\text{if $v = v_t$,}\\
f_{t-1}(v_t)&\text{if $v=v_{t-1}$,}\\
f_{t-1}(v)&\text{otherwise.}\end{array}\right.$$
 Consider now the sequence of acyclic orientations $(O_t)_{t=0,1,2,\dots}$ of $G_{\mathbf{r}}$ induced by the labellings $(f_t)_{t=0,1,2,\dots}$, so that
 for all $e=\{u,v\}\in E(G_{\mathbf{r}})$ and $t\geq 0$, we have that $O_t(e)=(u,v)$ if and only if $f_t(u)<f_t(v)$. The sequence
 $(O_t)_{t=0,1,2,\dots}$ is called the \emph{Sliding-$(n+1)$ (SL) stochastic process on the set of acyclic orientations of $G_{\mathbf{r}}$}.  
 \end{defn}
 
 \begin{theo}\label{theo:sliding}
The Sliding-$(n+1)$ stochastic process of $G_{\mathbf{r}}$ of Definition~\ref{defn:sliding} satisfies that,
if $S_{\mathbf{r}}$ is the set of all acyclic orientations of $G_{\mathbf{r}}$ whose unique maximal element is $\mathbf{r}$, then 
$\sum_{t=1}^{\infty} \proba{O_t\in S_{\mathbf{r}}}=\infty$ and for every $O\in S_{\mathbf{r}}$:
\begin{align*}
\left(\boldsymbol{\uppi^{\scriptscriptstyle \emph{SL}}}\right)_O=\uppi^{\scriptscriptstyle \emph{SL}}_O &:= 
\lim_{N\rightarrow \infty}\frac{\sum_{t=1}^{N}\proba{O_t=O}}{\sum_{t=1}^{N}\proba{O_t\in S_{\mathbf{r}}}}
=\frac{e(O\vert_{\scriptscriptstyle V})}{n!},
\end{align*}
where $O\vert_{\scriptscriptstyle V}$ is the restriction of $O$ to $E$ (hence an acyclic orientation of $G$) and 
$e(O\vert_{\scriptscriptstyle V})$ denotes the number of linear extensions of the induced poset $(V,\leq_{O\vert_{_V}})$. These
results hold independently of the initial choice of $O_0$.
\end{theo}
\begin{proof2}
Consider the simple graph $H$ on vertex set equal to the set of all bijective maps $V\sqcup\{\mathbf{r}\}\rightarrow[n+1]$, and where two
maps $f$ and $g$ are connected by an edge if and only if there exists $\{u,v\}\in E(G_{\mathbf{r}})$ such that $f(u)=g(v)=n+1$,
$f(v)=g(u)$, and $f(w)=g(w)$ for all $w\in V\backslash\{u,v\}$. If two bijective maps $f,g:V\sqcup\{\mathbf{r}\}\rightarrow[n+1]$ differ
only in one edge of $G_{\mathbf{r}}$, so that $f(u)=g(v)\neq n+1$ and $f(v)=g(u)\neq n+1$ for some $\{u,v\}\in E(G_{\mathbf{r}})$, but $f(w)=g(w)$ for all $w\in V\backslash\{u,v\}$, then
we can easily but somewhat tediously show that $f$ and $g$ belong to the same connected component of $H$, making use of the facts that vertex
$\mathbf{r}$ is adjacent to all other vertices of $G_{\mathbf{r}}$ and that $G$ is connected. But then, the proof
of Theorem~\ref{theo:edgelabrev} shows that $H$ is a connected graph. Now, the sequence 
$(f_t)_{t=0,1,2,\dots}$ of Definition~\ref{defn:sliding} is a simple random walk on $H$, and the degree of a bijective map $f:V\rightarrow[n]$ in $H$
is clearly $\degnovec{G_{\mathbf{r}}}(v_f)$, where $v_f\in V$ depends on $f$ and is such that $f(v_f)=n+1$, 
so the stationary distribution $\boldsymbol{\uppi}$ for this Markov chain satisfies that
$\uppi_{f}=c\cdot \degnovec{G_{\mathbf{r}}}(v_f)$, for some fixed normalization constant $c\in\RP$. The vertices of
$H$ that induce acyclic orientations of $G_{\mathbf{r}}$ from the set $S_{\mathbf{r}}$ are exactly the bijective maps $f:V\sqcup\{\mathbf{r}\}\rightarrow[n+1]$ such that
$f(\mathbf{r})=n+1$, and for these we have that $\uppi_f = c\cdot n$. 
The result then follows from the
construction of $(O_t)_{t=0,1,2,\dots}$ and from Equations~\ref{eq:statdist}.  
\par
\qed\end{proof2}
 
 \begin{defn}[Cover-Reversal Random Walk]\label{defn:coverrev}
 Let $G=G(V,E)$ be a simple graph with $|V|=n$, and select an arbitrary acyclic orientation 
 $O_0$ of $G$. Let us consider a sequence $(O_t)_{t=0,1,2,\dots}$ of
 acyclic orientations of $G$ such that for $t\geq 1$,
 $O_t$ is obtained from $O_{t-1}$ through the following random process:
 Let $(u,v)$ be selected uniformly at random from the set, 
 $$\emph{Cov}(O_{t-1}):=\left\{e\in O_{t-1}[E]:\text{$e$ represents a cover relation
 in $(V,\leq_{O_{t-1}})$}\right\},$$
 and for all $e\in E$, let,
 $$O_t(e) = \left\{\begin{array}{ll}(v,u)&\text{if $e=\{u,v\}$,}\\
 O_{t-1}(e)&\text{otherwise.} \end{array}\right.$$
The sequence
 $(O_t)_{t=0,1,2,\dots}$ is called the \emph{Cover-Reversal (CR) random walk on the set of acyclic orientations of $G$}.  
 \end{defn}
 
\begin{theo}\label{theo:coverrev}
The Cover-Reversal random walk in $G$ of Definition~\ref{defn:coverrev} is a simple
$2$-period random walk on the \emph{$1$-skeleton} of the clean graphical 
zonotope $\CZon{G}$ of Theorem~\ref{theo:polrealzon} (hence, on a particular simple connected bipartite graph on vertex set equal to the set of all
acyclic orientations of $G$), and its
stationary  distribution $\boldsymbol{\uppi^{\scriptscriptstyle \emph{CR}}}$ satisfies that, for every 
acyclic orientation $O$ of $G$:
$$\uppi^{\scriptscriptstyle \emph{CR}}_O=c\cdot |\emph{Cov}(O)|,$$
where $c\in\RP$ is a normalization constant independent of $O$. 
\end{theo}
\begin{proof2}
From the proof of Theorem~\ref{theo:polrealzon}, the edges of $\CZon{G}$ are in bijection with the set of
all \emph{p.a.o.}'s $O$ of $G$ such that if $\Upsigma$ is the connected partition associated to $O$, then
$|\Upsigma| = n-1$. Hence, the edges of $\CZon{G}$ are in bijection with the set of all pairs of the form $(e,O)$, where
$e\in E$ and $O$ is an acyclic orientation of the graph $G\slash e$, obtained from $G$ by contraction of the edge $e$. 
The two vertices of $\CZon{G}$ adjacent to an edge corresponding to a $(e,O)$ with $e=\{u,v\}$ are, respectively, 
obtained from the acyclic orientations
$O_1$ and $O_2$ of $G$ such that $O_1(e)=(u,v)$, $O_2(e)=(v,u)$, and 
such that $O_1\vert_{\scriptscriptstyle E\backslash e}=O_2\vert_{\scriptscriptstyle E\backslash e}$ are naturally induced by $O$ 
({\it e.g.} see Definition~\ref{defn:pao1}). Necessarily then, both $(u,v)$ and $(v,u)$ correspond respectively to cover 
relations in the posets $(V,\leq_{O_1})$ and $(V,\leq_{O_2})$, since otherwise the orientation $O$ of $G\slash e$ would not be acyclic. 
\par
On the other hand, given an acyclic orientation $O_1$ of $G$ and an edge $(u,v)\in O_1[E]$ such that
$v$ covers $u$ in $(V,\leq_{O_1})$, then, reversing the orientation of (only) that edge in $O_1$ yields a new acyclic orientation
$O_2$ of $G$, so $(v,u)\in O_2[E]$. Otherwise, using a directed cycle formed by edges from $O_2[E]$, which must then include the edge
$(v,u)$, we observe that the relation $u\leq_{O_1} v$ is a consequence
of other order relations in $(V,\leq_{O_1})$ and $v$ does not cover $u$ there. This is a contradiction, and it 
furthermore implies that both $O_1$ and $O_2$ naturally induce a well-defined acyclic orientation
$O$ of $G\slash\{u,v\}$. 
\par
Hence, the Cover-Reversal random walk of $G$ corresponds to a simple random walk on the $1$-skeleton
of $\CZon{G}$ (or of $\Zon{G}$) and the result follows now from Theorem~\ref{theo:randwalk}, since
this graph is connected and bipartite, clearly. 
\par
\qed\end{proof2}
\begin{rem}\label{rem:coverrev}
Variants of the Cover-Reversal random walk on $G$, obtained for example by flipping biased coins at each step, 
can be used
to obtain stochastic processes that converge to a uniform distribution on the set of acyclic orientations of $G$. However,
these variants are clearly not very illuminating or efficient. 
\end{rem}
 
 \begin{defn}[Interval-Reversal Random Walk]\label{defn:interrev}
 Let $G=G(V,E)$ be a simple graph with $|V|=n$, and select an arbitrary acyclic orientation 
 $O_0$ of $G$. Let us consider a sequence $(O_t)_{t=0,1,2,\dots}$ of
 acyclic orientations of $G$ such that for $t\geq 1$,
 $O_t$ is obtained from $O_{t-1}$ through the following random process:
 Let $\{u,v\}\in E$ be selected uniformly at random from this set, with
 $(u,v)\in O_{t-1}[E]$, and
for all $e=\{x,y\}\in E$ with $(x,y)\in O_{t-1}[E]$, let,
 $$O_t(e) = \left\{\begin{array}{ll}(y,x)&\text{if $u\leq_{O_{t-1}}x<_{O_{t-1}} y\leq_{O_{t-1}} v$,}\\
 (x,y)=O_{t-1}(e)&\text{otherwise.} \end{array}\right.$$
The sequence
 $(O_t)_{t=0,1,2,\dots}$ is called the \emph{Interval-Reversal (IR) random walk on the set of acyclic orientations of $G$}.  
 \end{defn}

\begin{lem}\label{lem:interrev}
Let $G=G(V,E)$ be a simple graph and let $O$ be any given acyclic orientation of $G$. For an arbitrary edge $\{u,v\}\in E$, say with
$(u,v)\in O[E]$, let us 
define a new orientation $O_{\scriptscriptstyle \{u,v\}}$ of $G$ by requiring that, for all $e=\{x,y\}\in E$ with $(x,y)\in O[E]$:
$$O_{\scriptscriptstyle\{u,v\}}(e) = \left\{\begin{array}{ll}(y,x)&\text{if $u\leq_{O}x<_{O} y\leq_{O} v$,}\\
 (x,y)=O(e)&\text{otherwise.} \end{array}\right.$$
 Then, $O_{\scriptscriptstyle\{u,v\}}$ is also an acyclic orientation of $G$ and, furthermore, 
 $\left(O_{\scriptscriptstyle\{u,v\}}\right)_{\scriptscriptstyle\{u,v\}}=O$.
 \par
 Additionally, for any choice of $e_1,e_2\in E$, we have that $O_{e_1} = O_{e_2}$ if and only if
 $e_1=e_2$.  
\end{lem}
\begin{proof2}
Suppose on the contrary that $O_{\scriptscriptstyle\{u,v\}}$ is not an acyclic orientation of $G$. Then, 
there exists at least one directed cycle
$C\subseteq O_{\scriptscriptstyle\{u,v\}}[E]$ that has the following form: 
\begin{enumerate}[label =\bfseries ${\scriptstyle}$]
\item For $E^{\scriptscriptstyle O}_{\scriptscriptstyle (u,v)} = \left\{\{x,y\}\in E:u\leq_O x<_{O}y\leq_{O} v\right\}$,
\item there exists $k\in\PN$ and pairwise disjoint non-empty sets, 
$$P_1,Q_1,P_2,Q_2,\dots,P_k,Q_k\subseteq O_{\scriptscriptstyle\{u,v\}}[E]\text{ with $C = \bigcup_{i=1}^k \left(P_i\cup Q_i\right)$},$$
\item such that for all $i\in[k]$, 
\item $P_i=\left\{(p^{i}_{j-1},p^{i}_j)\right\}_{j=1,\dots,|P_i|}\subseteq 
O_{\scriptscriptstyle\{u,v\}}\left[E^{\scriptscriptstyle O}_{\scriptscriptstyle (u,v)} \right]$,
\item $Q_i=\left\{(q^{i}_{j-1},q^{i}_j)\right\}_{j=1,\dots,|Q_i|}\subseteq 
O_{\scriptscriptstyle\{u,v\}}[E]\backslash \left(O_{\scriptscriptstyle\{u,v\}}\left[E^{\scriptscriptstyle O}_{\scriptscriptstyle (u,v)} \right]\right)$,
\item $p^{i}_{|P_i|}=q^{i}_0$ and $q^{i}_{|Q_i|}=p^{i+1}_0$, where $p^{k+1}_0:= p^{1}_0$. 
\end{enumerate}
This is true simply because any directed cycle in $O_{\scriptscriptstyle\{u,v\}}[E]$ must necessarily involve edges from both
$E^{\scriptscriptstyle O}_{\scriptscriptstyle (u,v)}$ and $E\backslash E^{\scriptscriptstyle O}_{\scriptscriptstyle (u,v)}$. 
\par
Since {\bf 1)} $O$ and $O_{\scriptscriptstyle\{u,v\}}$ agree on $E\backslash E^{\scriptscriptstyle O}_{\scriptscriptstyle (u,v)}$, 
{\bf 2)} $u\leq_{O} p^{1}_{|P_1|}, p^{2}_0 \leq_{O} v$, and {\bf 3)} $q^{1}_0=p^{1}_{|P_1|}$, $q^{1}_{|Q_1|}=p^{2}_0$, then 
$u\leq_O q^{1}_0\leq_{O} q^{1}_{1}\leq_O\dots\leq_O q^{1}_{|Q_1|}\leq_O v$, so in particular $\{q^{1}_{0},q^{1}_{1}\}\in E^{\scriptscriptstyle O}_{\scriptscriptstyle (u,v)}$
by definition,
a contradiction with the construction of $C$. Hence, $O_{\scriptscriptstyle\{u,v\}}$ is also an acyclic orientation of $G$. 
\par
To prove that $\left(O_{\scriptscriptstyle\{u,v\}}\right)_{\scriptscriptstyle\{u,v\}}=O$, it suffices to check that if
for some $\{x,y\}\in E$ with $(y,x)\in O_{\scriptscriptstyle\{u,v\}}[E]$ we have that 
$v\leq_{O_{\scriptscriptstyle\{u,v\}}} y<_{O_{\scriptscriptstyle\{u,v\}}} x\leq_{O_{\scriptscriptstyle\{u,v\}}} u$,
then in fact $u\leq_{O} x<_{O}y\leq_{O} v$. Somewhat analogously with the previous argument, suppose on the contrary
that there exists some $\{x,y\}\in E$ with $(y,x)\in O_{\scriptscriptstyle\{u,v\}}[E]$ for which the condition fails to hold. 
Then, inside any directed path $P=\{(p_{j-1},p_j)\}_{j=1,\dots,|P|}\subseteq O_{\scriptscriptstyle\{u,v\}}[E]$ such that
$(y,x)\in P$, $p_0=v$, and $p_{|P|}=u$, there must exist a maximal (by containment) sub-path $Q=\{(q_{j-1},q_j)\}_{j=1,\dots,|Q|}\subseteq P$ such that
$(y,x)\in Q\subseteq O_{\scriptscriptstyle\{u,v\}}[E]\backslash \left(O_{\scriptscriptstyle\{u,v\}}\left[E^{\scriptscriptstyle O}_{\scriptscriptstyle (u,v)} \right]\right)$.
Necessarily then, $u\leq_O q_0 <_O q_{|Q|}\leq_{O} v$, so $u\leq_O q_0 \leq_O y <_O x\leq_O q_{|Q|}\leq_{O} v$, and hence 
$\{y,x\}\in E^{\scriptscriptstyle O}_{\scriptscriptstyle (u,v)} $. This is a contradiction. 
\par
The last statement is a simple consequence of observing that, for every choice of $\{u,v\}\in E$, $u$ and $v$ determine a unique
interval inside each of the posets $(V,\leq_{O})$, where $O$ is an acyclic orientation of $G$: A non-empty closed interval of a finite poset
is uniquely determined by its maximal and minimal elements.  
\par
\qed\end{proof2}

\begin{prop}\label{prop:interrev}
In Lemma~\ref{lem:interrev}, consider the simple graph $\emph{AO}_{G}^{\emph{inter}}$ on vertex set equal to the set of all acyclic orientations of $G$, and in which
two acyclic orientations $O_1$ and $O_2$ of $G$ are connected by an edge, if and only if there exists $\{u,v\}\in E$ such that
$\left(O_1\right)_{\scriptscriptstyle \{u,v\}} = O_2$. Then, $\emph{AO}_{G}^{\emph{inter}}$ is an $|E|$-regular connected graph. 
\end{prop}
\begin{proof2}
Firstly, let us note that $\emph{AO}_{G}^{\emph{inter}}$ is indeed a well-defined simple graph 
(so it does not have loops or multiple edges) per the three main statements of Lemma~\ref{lem:interrev}. 
Now, we point out that $\emph{AO}_{G}^{\emph{inter}}$ contains as a spanning sub-graph the $1$-skeleton of the (clean) graphical zonotope
$\CZon{G}$ since, colloquially, all \emph{cover-reversals} are also \emph{interval-reversals}. Hence, since the later graph has been
observed to be connected in the proof of Theorem~\ref{theo:coverrev}, then $\emph{AO}_{G}^{\emph{inter}}$ is also connected. 
Every vertex of this graph must have degree $|E|$, clearly.  
\par
\qed\end{proof2}

\begin{theo}\label{theo:interrev}
The Interval-Reversal random walk in $G$ of Definition~\ref{defn:interrev} is a simple
random walk on the graph $\emph{AO}_{G}^{\emph{inter}}$ of Proposition~\ref{prop:interrev} 
(hence, on a particular regular connected graph on vertex set equal to the set of all
acyclic orientations of $G$), and its
stationary  distribution $\boldsymbol{\uppi^{\scriptscriptstyle \emph{IR}}}$ satisfies that, for every 
acyclic orientation $O$ of $G$:
$$\uppi^{\scriptscriptstyle \emph{IR}}_O=\frac{1}{\vert \chi_G(-1)\vert},$$
where $\vert \chi_G(-1)\vert$ is the number of acyclic orientations of $G$~{\bfseries [}\cite{stanleyacyclic}{\bfseries ]}. 
\end{theo}
\begin{proof2}
That the Interval-Reversal random walk of $G$ corresponds to a simple random walk on $\emph{AO}_{G}^{\emph{inter}}$ is
a direct consequence of Lemma~\ref{lem:interrev}. That $\emph{AO}_{G}^{\emph{inter}}$ is connected and $|E|$-regular is the content of
Proposition~\ref{prop:interrev}, so we can now rely on Theorem~\ref{theo:randwalk} to
obtain the result. 
\par
\qed\end{proof2}

\begin{figure}[ht]
\begin{tabular}{cc}
\begin{subfigure}{.5\textwidth}
  \centering
  \includegraphics[width=0.9\linewidth]{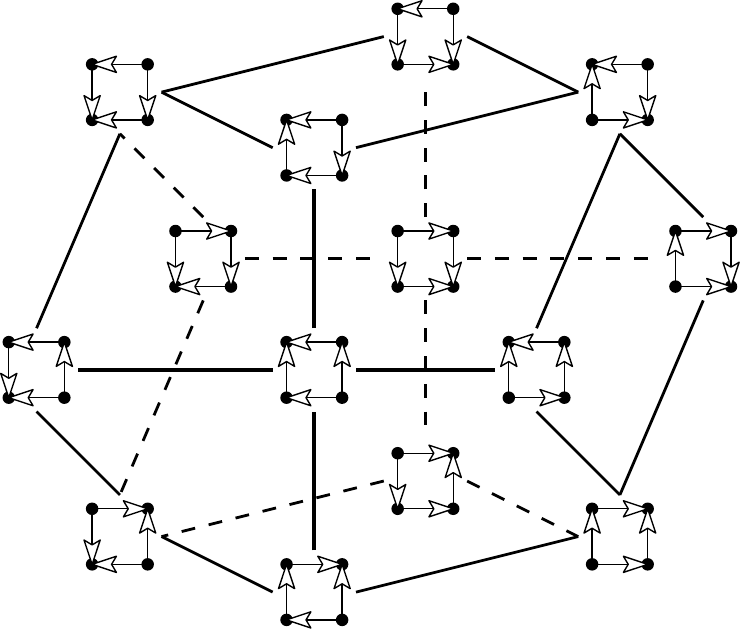}
  \caption{}
  \label{fig:stoc1}
\end{subfigure}
&
\begin{subfigure}{.5\textwidth}
  \centering
  \includegraphics[width=0.9\linewidth]{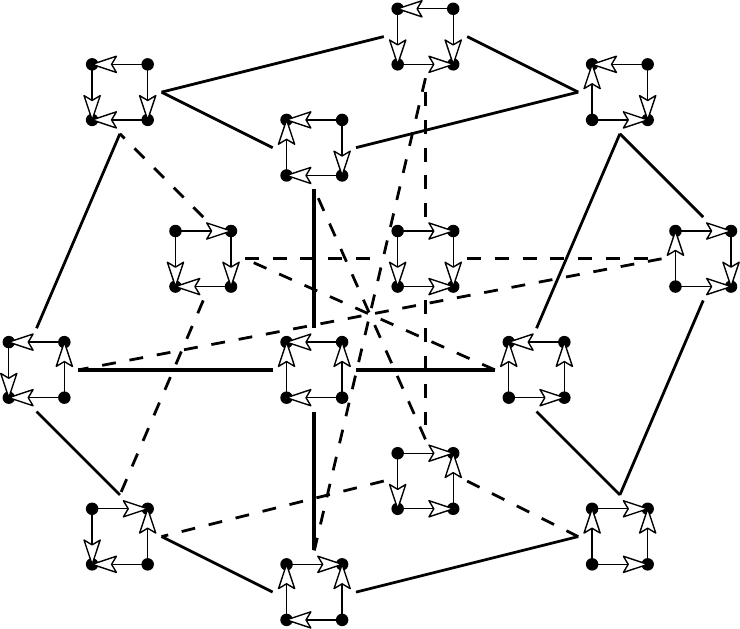}
  \caption{}
  \label{fig:stoc2}
\end{subfigure}
\end{tabular}

\caption{ Examples of Definitions~\ref{defn:coverrev} and \ref{defn:interrev} for the
$4$-cycle $C_4$. In \ref{fig:stoc1}, we present the $1$-skeleton of the graphical zonotope of $C_4$, a \emph{rhombic dodecahedron}, 
where the Cover-Reversal random walk runs; notably, it is not a regular
graph. If four diagonals are added to the graph as shown in \ref{fig:stoc2}, 
we obtain a $4$-regular graph, $\emph{AO}_{C_4}^{\emph{inter}}$ in Proposition~\ref{prop:interrev},
where the Interval-Reversal random walk runs.}
\label{fig:stoc}
\end{figure}

 \subsection{Acyclic Orientations of a Random Graph.}\label{subsec:randgraph}
 
 This short subsection is aimed at proving a surprising formula for the expected number of acyclic orientations
 of an Erd\"os-R\'enyi random graph from $\rgr{[n]}{p}$, with $p\in (0,1)$. This formula will follow from the results of Section~\ref{sec:noncrosstrees}, and more
 specifically from those of Subsection~\ref{subsec:standardmon}. 
 
 \begin{defn}\label{defn:park}
 Let $n\in \PN$. A \emph{parking function of $[n]$} is a vector $\mathbf{a} \in \NN^{[n]}$ such that for
 any $\upsigma\in\mathfrak{S}_{[n]}$ with
 $a_{\upsigma(1)}\leq a_{\upsigma(2)}\leq \dots\leq a_{\upsigma(n)}$, we have that
 $a_{\upsigma(i)}\leq i-1$ for all $i\in[n]$. The set of all parking functions of $[n]$ will be denoted by
 $\emph{Park}_{[n]}$. 
 \par
 For $\mathbf{a}\in \NN^{[n]}$, let us write $\emph{Area}(\mathbf{a}) := a_1+a_2+\dots+a_n$ and
 $\emph{supp}(\mathbf{a}):=\{i\in[n]:a_i>0\}$.
 \end{defn}
 
 \begin{theo}\label{theo:park}
 Let $n\in \PN$, $p\in(0,1)$, and $G\sim \rgr{[n]}{p}$. Write $q:=1-p$. If $\vert\chi_G(-1)\vert$ is the number of
 acyclic orientations of $G$, and we let $S:=\kfield[x_1,\dots,x_n]\slash T_{K_{[n]}}$, where as usual $K_{[n]}$ denotes the complete graph on 
 vertex set $[n]$, 
 then we have:
 \begin{align}
 \expe{\vert\chi_G(-1)\vert} &= q^{\binom{n}{2}}\cdot\sum_{\mathbf{a}\in 
 \emph{Park}_{[n]}}\left(\frac{1}{q}\right)^{\emph{Area}(\mathbf{a}) }p^{\vert\emph{supp}(\mathbf{a})\vert}.
 \end{align}
 \end{theo}
 \begin{proof2}
We make use of Proposition~\ref{prop:acyctrees}. In general, for any simple graph $H$ on vertex set $[n]$ (as $G$ here and the complete graph
$K_{[n]}$), we will let
$\mathbf{T}^{\mathbf{r}}H$ be the set of all $\mathbf{r}$-rooted spanning trees of $H_{\mathbf{r}}$. 
Now, for $T\in \mathbf{T}^{\mathbf{r}}K_{[n]}$, we will say
that $T$ is \emph{useful} if $T\in \mathbf{T}^{\mathbf{r}}G$ 
and its unique depiction function $p$ of Theorem~\ref{theo:noncross}, satisfies the conditions of
Proposition~\ref{prop:acyctrees}. Then:
\begin{align*}
\expe{\vert \chi_{G}(-1)\vert} &= \sum_{T\in \mathbf{T}^{\mathbf{r}}K_{[n]}}\proba{T\text{ is \emph{useful}}}&\\
&=\sum_{T\in \mathbf{T}^{\mathbf{r}}K_{[n]}}\left\lgroup\proba{T\in \mathbf{T}^{\mathbf{r}}G}\right\rgroup\cdot \left\lgroup\proba{\{i,j\}\not\in E(G_{\mathbf{r}})\text{ for all }i,j\in[n],\right.&\\
& \left. \ \ \ \ \ \ \ \ \ \ \ \ \ \ \ \ \ \ \ \ \ \ \ \ \ \ \ \ \ \ \ \ \ \ \ \ \ \ \ \ \ \ \ \text{$(i,i_{\mathbf{r}})\in E(T)$, $p(i_{\mathbf{r}})<p(j)< p(i)$}}\right\rgroup &\\
&&\\
&=\sum_{T\in \mathbf{T}^{\mathbf{r}}K_{[n]}} \frac{p^{n}}{p^{\degnovec{T}(\mathbf{r})}}\cdot\prod_{\substack{i\in[n]\\
(i,i_{\mathbf{r}})\in E(T)}}\proba{\{i,j\}\not\in E(G_{\mathbf{r}})\text{ for all $j\in[n]$, $p(i_{\mathbf{r}})<$}&\\
& \ \ \ \ \ \ \ \ \ \ \ \ \ \ \ \ \ \ \ \ \ \ \ \ \ \ \ \ \ \ \ \ \ \ \ \ \ \ \ \ \ \ \ \ \text{$p(j)< p(i)$}}&\\
&&\\
&=\sum_{T\in \mathbf{T}^{\mathbf{r}}K_{[n]}}\frac{p^n}{p^{\degnovec{T}(\mathbf{r})}}\prod_{\substack{i\in[n]\\(i,i_{\mathbf{r}})\in E(T)}}
q^{p(i)-1-p(i_{\mathbf{r}})}&\\
&=q^{\binom{n}{2}}\cdot\sum_{T\in \mathbf{T}^{\mathbf{r}}K_{[n]}}\left(\frac{1}{q}\right)^{\sum_{i\in[n]}p(i_\mathbf{r})}p^{\vert\{i\in[n]:a_i>0\}\vert}&\\
&= q^{\binom{n}{2}}\cdot\sum_{\mathbf{a}\in \emph{Park}_{[n]}}\left(\frac{1}{q}\right)^{\emph{Area}(\mathbf{a}) }p^{\vert\emph{supp}(\mathbf{a})\vert},&
\end{align*}
as we wanted.
\par
\qed\end{proof2}
 
 \subsection{$k$-Neighbor Bootstrap Percolation.}\label{subsec:kbootstrap}
 
\begin{defn}\label{defn:kbootstrap}
Let $G=G([n],E)$ be a finite simple graph, $k\in\PN$, and $A\subseteq [n]$. The \emph{$k$-neighbor bootstrap percolation on $G$ with initial set $A$}, is
the process $\{A_t\}_{t=0,1,2,\dots}$, where $A_0=A$ and $A_{t}=A_{t-1}\cup\{i\in [n]:\vert N_G(i)\cap A_{t-1}\vert\geq k\}$ for all $t\geq 1$.
The \emph{closure} of $A$ is the set $\emph{cl}(A):=\cup_{t\geq 0}A_t$, and we say that \emph{$A$ percolates in $G$} if $\emph{cl}(A)=[n]$.    
\end{defn}

\begin{quest}\label{quest:perc}
Given a graph $G$ as in Definition~\ref{defn:kbootstrap}, what is minimal size $|A|$ of $A\subseteq [n]$ such that
$A$ percolates in $G$? 
\end{quest}

\begin{defn}\label{defn:kbootstrapideal}
 For fixed $G$ and $k$ as in Definition~\ref{defn:kbootstrap}, let $\mathcal{C}{\scriptstyle(G,k)}:=
 \{\upsigma\subseteq[n]:\text{$\degout{G}{\upsigma}(i)<k$ for all $i\in\upsigma$}\}$. 
 The \emph{$k$-bootstrap percolation ideal $B_{\mathcal{C}{\scriptscriptstyle (G,k)}}$ of $G$} is the square-free monomial ideal of 
$\kfield[x_1,\dots,x_n]$ generated as: 
$$B_{\mathcal{C}{\scriptscriptstyle (G,k)}}=\left\langle \prod_{i\in\upsigma}x_i:\upsigma\in\mathcal{C}{\scriptstyle(G,k)} \right\rangle.$$ 
\end{defn}

\begin{prop}\label{prop:perc1}
In Definitions~\ref{defn:kbootstrap}-\ref{defn:kbootstrapideal}, the function that associates to each standard monomial 
$\mathbf{x}^{\mathbf{b}}\not\in B_{\mathcal{C}{\scriptscriptstyle (G,k)}}$, $\mathbf{b}\in\NN^{[n]}$, the set
of vertices $\{i\in[n]:b_i=0\}$ of $G$, restricts to a bijection between the set of all square-free standard monomials of 
$B_{\mathcal{C}{\scriptscriptstyle (G,k)}}$ and the set of all $A\subseteq[n]$ such that $A$ percolates in $G$. 
\par
Colloquially, the percolating sets of $G$ are in bijection with the supporting sets of standard monomials of the ideal 
$B_{\mathcal{C}{\scriptscriptstyle (G,k)}}$.  
\end{prop}
\begin{proof2}
Let $A\subseteq [n]$ be such that $\emph{cl}(A)\subsetneq [n]$, and consider the set 
$\upsigma:=[n]\backslash \emph{cl}(A)$. Necessarily, every element of $\upsigma$ must have fewer than $k$ 
neighbors inside $\emph{cl}(A)$, so $\degout{G}{\upsigma}(i)< k$, for all $i\in\upsigma$. This implies that
$\upsigma\in \mathcal{C}{\scriptstyle(G,k)}$, and $\mathbf{x}^{\mathbf{\upsigma}}:=\prod_{i\in\upsigma}x_i\in B_{\mathcal{C}{\scriptscriptstyle (G,k)}}$.
But then, since $\upsigma\subseteq \widehat{A}:=[n]\backslash A$, we have that $\mathbf{x}^{\mathbf{\upsigma}}\vert
\mathbf{x}^{\boldsymbol{\widehat{A}}}:=\prod_{i\in\widehat{A}}x_i$, so
$\mathbf{x}^{\boldsymbol{\widehat{A}}}\in B_{\mathcal{C}{\scriptscriptstyle (G,k)}}$ as well.
\par
On the contrary, if $\mathbf{x}^{\boldsymbol{\widehat{A}}}\in B_{\mathcal{C}{\scriptscriptstyle (G,k)}}$ for some
$A\subseteq[n]$, there must exist some $\upsigma\in \mathcal{C}{\scriptstyle(G,k)}$ such that
$\mathbf{x}^{\mathbf{\upsigma}}\vert
\mathbf{x}^{\boldsymbol{\widehat{A}}}$. Necessarily then, $\emph{cl}(A)\subseteq [n]\backslash\upsigma$, since it is never possible to
percolate the elements of $\upsigma$ during a $k$-bootstrap percolation on $G$ from an initial set disjoint from $\upsigma$, as $A$ here.
\par
\qed\end{proof2}

\nocite{*}
\bibliographystyle{abbrvnat}
\bibliography{dual}
\label{sec:biblio}

\end{document}